%% file: preprint_shortened.tex
\documentclass[a4paper,oneside,reqno]{article}
\usepackage[utf8]{inputenc}
\usepackage{lipsum}
\usepackage{standalone}
\usepackage{atbegshi}
\AtBeginDocument{\AtBeginShipoutNext{\AtBeginShipoutDiscard}}

\usepackage{float}

\newcommand{\cJnu}{\mathcal{J}_i}
\newcommand{\nep}{NEP}
\newcommand{\gnep}{\ensuremath{\mathrm{GNEP}}}
\newcommand{\vep}{\ensuremath{\mathrm{VEP}}}
\newcommand{\Vhat}{\ensuremath{\widehat V}}
\newcommand{\cBhat}{\ensuremath{\widehat{\calB}}}

\newcommand{\pengnep}{\eqref{eq:gnep:penalized}}
\newcommand{\penvep}{\eqref{eq:vep:penalized}}

\newcommand\fenics{FEniCS}

\usepackage{mystyle}
\usepackage{review}

\title{Path-Following Methods for Generalized Nash Equilibrium Problems}
\author{Steven-Marian Stengl\thanks{Weierstrass Institute, Mohrenstrasse 39, 10117 Berlin, Germany (\texttt{steven-marian.stengl@wias-berlin.de}) and Humboldt-Universit\"at zu Berlin, Unter den Linden 6, 10099 Berlin, 
Germany (\texttt{stengl@math.hu-berlin.de})}}
\date{\today}
\addtocounter{page}{-1}%
\begin{document}%
\thanks{This research was supported by the DFG under the grant HI 1466/10-1  associated to the project ``Generalized Nash Equilibrium Problems with Partial Differential Operators: Theory, Algorithms, and Risk Aversion'' within the priority programme SPP1962 ``Non-smooth and Complementarity-based Distributed Parameter Systems: Simulation and Hierarchical Optimization''.}
\maketitle%
\begin{abstract}
Building upon the results in \cite{bib:HintSurKaemmler}, generalized Nash equilibrium problems are considered, in which the feasible set of each player is influenced by the decisions of their competitors.
This is realized via the existence of one (or more) state constraint(s) establishing a link between the players.
Special emphasis is put on the situation of a state encoded in a possibly non-linear operator equation. First order optimality conditions under a constraint qualification are derived.
\msa{Aiming at a practically meaningful method, an}{} approximation scheme using a penalization technique leading to a sequence of (Nash) equilibrium problems without \msa{dependence of the}{} constraint set \msa{on the other players' strategies}{} is established.
An associated path-following strategy related \ms{to}{} a \mh{value}{} function is then proposed.
This happens at first on the most abstract level and is subsequently established to a narrower framework geared to the presence of partial differential equations in the constraint.
Our findings are illustrated \msa{with}{} examples \msa{having}{} \msa{distributed}{} and boundary controls---both involving semi-linear elliptic PDEs.
\end{abstract}
\section{Introduction}
In recent years a growing research effort was focused on \msa{the}{}
theoretical as well as numerical treatment of (generalized) Nash equilibrium
problems (\msa{abbr.:}{} (G)NEPs)\ms{.}{} For problems in finite dimensions a high level of sophistication \mh{has been}{}
reached \mh{there}{}. We refer to the \ms{overview articles \cite{bib:PangFukushima, bib:KanzowFacchinei}}{} \ms{as well as the}{} contributions \cite{bib:KanzowHeusingerOptReformulation,
bib:FacchineiFischerPiccialli,bib:DrevesHeusingerKanzowFukushima} and the references given therein.
For GNEPs within a function space setting
the level of research is \mh{far}{} less complete\ms{, but is nonetheless an active working field (cf. \cite{bib:BorziKanzow, bib:HintSurGNEP, bib:HintSurKaemmler, bib:KanzowKarlSteckDWachsmuth}).}{}\\
In this article we \mh{study}{} a generalization of the techniques developed
and used in \cite{bib:HintSurKaemmler} as well as \cite{bib:HintSurGNEP}. The
approach developed therein \mh{focussed}{} on a class of Nash games with tracking
functionals involving a partial differential equation (\msa{abbr.:}{} PDE\mh{)}{} constraint along with control as well as state
constraints. The latter condition establish\mh{e}s an influence of the players on their
feasible strategy choices of the players' strategy sets. The techniques have been
developed and successfully applied to a selection of linear PDEs (elliptic, parabolic
and hyperbolic). So far, non-linear PDEs have not been considered in this regard. One
of the main reasons is the difficulty to establish suitable conditions \mh{for}{} the
existence of Nash equilibria.\\
Within the scope of this work we extend the methods in \cite{bib:HintSurKaemmler} to
include a class of GNEPs with possibly non-linear operator equation\mh{s}. We
distinguish two different levels of abstraction compared to the setting therein: In the first step we are guided by
the typical form discussed in (abstract) optimization as well as \mh{in the}{} 
literature \mh{on Nash games}{} and derive a path-following scheme. The \mh{obtained}{} results are then applied to
a framework similar to the one in \cite{bib:HintSurKaemmler}\mh{,}{} specifically
tailored to the treatment of (non-)linear PDE constraints.\\
The remaining sections are organized as follows. In section 2 we introduce
notations and important definitions which are vital for the understanding of the
remainder of the work. Section 3 is devoted to the introduction of the concept of
(generalized) Nash equilibrium problems and the derivation of first order conditions for Nash
equilibria in form of Lagrange multiplier systems. In section 4  a penalization and path-following technique is used to derive asymptotic
convergence results. Our findings therein are utilized in Section 5 to the
specialized framework to facilitate the application to PDE related problems. In
Section 6 the results are illustrated by a selection of problems involving semi-linear elliptic PDEs with distributed as well as boundary controls.
\section{Notation and Preliminaries}
In the following, for a given Banach space $X$ denote by $X^*$ its \emph{topological dual space} and the associated dual pairing $\langle \ccdot, \ccdot \rangle_{X^*,X} : X^* \times X \rightarrow \R$ by $\langle x^*, x \rangle_{X^*,X} := x^*(x)$.
Oftentimes we simply denote $\langle \ccdot, \ccdot \rangle$, if the corresponding
spaces are clear from the context. Two elements $x^* \in X^*$ and $x \in X$ are called \emph{orthogonal}, if $\langle x^*, x \rangle = 0$ and we write $x^* \perp x$ or $x \perp x^*$. We write $x_n \to x$ for the strong convergence and $x_n \rightharpoonup x$ to denote the weak convergence. \\
A subset $C \subseteq X$ is called \emph{convex}, if for all $t \in (0,1)$ and $x_0,x_1 \in C$ holds $t x_1 + (1-t) x_0 \in C$. A set
$K \subseteq X$ is called a \emph{cone}, if for all $t \in \R$, $t \geq 0$ and 
$x \in K$ also $tx \in K$ holds.
The \emph{normal cone} of a nonempty, convex, closed set $C$ is defined as 
\begin{eq*}
N_C(x) := \left\{ x^* \in X^* : \langle x^*, x' - x \rangle_{X^*,X} \leq 0 \text{ \msa{for}{} all } x' \in C \right\}.
\end{eq*}
Consider a functional $f\msa{: X}{} \rightarrow \R \cup \{\pm\infty\}$. It is called \lsc{}, if for all sequences $x_n \rightarrow x$ in $X$ also
\begin{eq*}
f(x) \leq \liminf_{n \rightarrow \infty} f(x_n)
\end{eq*}
holds true and it is called \emph{proper}, if $f(x) > -\infty$ for all $x \in X$ holds and there exists $x \in X$ with $f(x) < +\infty$. The \emph{epigraph} of $f$
is defined as the set
\begin{eq*}
\epi{f} := \left\{ (x,\alpha) \in X \times \R : f(x) \leq \alpha \right\}.
\end{eq*}
The \emph{sublevel set} \msa{of $f$}{} with respect to a threshold $\alpha \in \R$ is defined as
\begin{eq*}
\lev_\alpha{f} := \left\{x \in X : f(x) \leq \alpha \right\}.
\end{eq*}
The functional $f$ is called \emph{convex}, if for all $x_0,x_1 \in X$ and $t \in (0,1)$ holds
\begin{eq*}
f(t x_1 + (1-t) x_0) \leq t f(x_1) + (1-t) f(x_0).
\end{eq*}
In fact, many of the above defined analytical properties can be characterized
geometrically. As proven in \cite[Proposition 2.5]{bib:BarbuPrecupanu} or \cite[Proposition 2.5]{bib:BonnansShapiro}, the lower semi-continuity is equivalent to the closedness of the epigraph as well as to the closedness of all sublevel sets. The convexity of $f$ is equivalent with the convexity of $\epi{f}$ and implies the convexity of $\lev_\alpha{f}$ for all $\alpha \in \R$ (see \cite[Proposition 2.3]{bib:BarbuPrecupanu}, \cite[Section 3.4.2]{bib:BoydVandenberghe}).\\
The \emph{closed unit ball} of
$X$ is denoted as
\begin{eq*}
\B_X := \{x \in X : \|x\| \leq 1\}.
\end{eq*}
The \emph{interior} of a set $M \subseteq X$ is defined as
\begin{eq*}
\interior{M} := \left\{ x \in M : \text{ there exists } \varepsilon > 0 :  x + \varepsilon \B_X \subseteq M \right\}.
\end{eq*}
For two sets $M_1,M_2 \subseteq X$ their \emph{Minkowski sum} is defined by
\begin{eq*}
M_1 + M_2 = \left\{ x_1 + x_2 : x_i \in M_i \text{ for } i = 1, 2 \right\}.
\end{eq*}
For another given Banach space $Y$, a function $F: X \rightarrow \mathcal{P}(Y)$ is called a \emph{set-valued operator} or \emph{correspondence} and is denoted by $F: X \rra Y$. Its \emph{graph} is defined by
\begin{eq*}
\gph{F} := \left\{ (x,y) \in X \times Y : y \in F(x) \right\}
\end{eq*}
and its domain by
\begin{eq*}
\mathcal{D}(F) := \left\{ x\in X : F(x) \neq \emptyset \right\}.
\end{eq*}
The graph is called \emph{closed}, if for all sequences $(x_n, y_n)_{n \in \N} \subseteq \gph{F}$ with $(x_n,y_n) \rightarrow (x,y)$ in $X\times Y$ also $(x,y) \in \gph{F}$ follows. It is called \emph{weakly closed}, if for all sequences $(x_n, y_n)_{n \in \N} \subseteq \gph{F}$ with $(x_n,y_n) \rightharpoonup (x,y)$ in $X\times Y$ also $(x,y) \in \gph{F}$ follows.\\
An operator $A : X \rightarrow Y$ is called \emph{bounded}, if bounded subsets of $X$ are mapped to bounded subsets of $Y$. It is called \emph{continuous}, if for a sequence $x_n \rightarrow x$ in $X$ the convergence $A(x_n) \rightarrow A(x)$ in $Y$ holds. The operator is called \emph{weakly continuous}, if for a sequence $x_n \rightharpoonup x$ in $X$ the convergence $A(x_n) \rightharpoonup A(x)$ in $Y$ holds. \ms{Moreover, the operator $A$ is called \emph{Fr\'echet differentiable} in $x \in X$, if there exists an operator $DA(x) \in \calL(X,Y)$ such that 
\begin{eq*}
\lim_{h \rightarrow 0} \frac{\| A(x+h) - A(x) - DA(x)h \|_Y}{\|h\|_X} = 0
\end{eq*}
holds. If the space $X$ admits a product structure $X := \prod_{i = 1}^n X_i$ for some $n \in \N$ and a family of Banach spaces $(X_i)_{i = 1}^N$, then we denote with $\partial_i A(x) \in \calL(X_i,Y)$ the \emph{partial derivative} of $A$ in the $i$-th component defined by \msa{$DA(x) h := \sum_{i = 1}^N \partial A(x)h_i$}{} with $h = (h_i)_{i = 1}^{\msa{N}} \in X$.}\\
Let $d \in \N\backslash\{0\}$ and let $\Omega \subseteq \R^d$ be a bounded, open
domain. For $p \in [1, \infty)$ denote the Lebesgue space as
\begin{eq*}
L^p(\Omega) := \left\{ u: \Omega \rightarrow \R \text{ measurable }: \int_\Omega |u|^p \dx < +\infty \right\}
\end{eq*}
with its elements only identified up to null sets, i.e. sets of Lebesgue measure zero.
This space equipped with the norm $\|u\|_{L^p\msa{(\Omega)}} := \left( \int_\Omega |u|^p \dx\right)^{\frac{1}{p}}$ is a Banach space for all $p \in [1,\infty)$ and a reflexive Banach space for $p \in (1,\infty)$. The \emph{Sobolev spaces} 
$W^{1,p}(\Omega)$ are defined as 
\begin{eq*}
W^{1,p}(\Omega) := \left\{ u \in L^p(\Omega) : \nabla u \in L^p(\Omega; \mathbb{R}^d) \right\},
\end{eq*}
where $\nabla u$ denotes the distributional derivative of $u$. Equipped with the norm $\|u\|_{W^{1,p}} := \left( \|u\|^p_{L^p\msa{(\Omega)}} + \sum_{i = 1}^d \| \partial_i u\|^p_{L^p\msa{(\Omega)}} \right)^{\frac{1}{p}}$, the space $W^{1,p}(\Omega)$ is a Banach space, and for $p \in (1,\infty)$ a reflexive Banach space. In case of $p = 2$ one also denotes $H^1(\Omega) = W^{1,2}(\Omega)$.
\mh{\section{Generalized Nash Equilibrium Problems}}\label{sec:general}
\msa{After the introduction of notions and concepts from functional and convex analysis we devote our attention to Nash games. In the scope of this section we introduce first the notion of (generalized) Nash equilibrium problems as well as of related equilibrium problems. After the discussion of essential properties and relations between them, we turn our attention to the derivation of necessary first order conditions in form of a Lagrange multiplier systems.
\subsection{Definition of (Generalized) Nash Equilibrium Problems and Variational Equilibrium Problems}}{}
Let a family of Banach spaces $U_i$ be given. Define the space $U := U_1 \times \cdots \times U_N$ \ms{as well as the \emph{strategy sets} $U_\ad^i \subseteq U_i$ and the joint strategy set $U_\ad := \prod_{i = 1}^N U_\ad^i$}{} together with a family of real-valued functionals $\calJ_i: \ms{U_\ad} \rightarrow \R$ for $i = 1, \dots, N$. \mh{With the index $-i$ we denote strategies, where the $i$-th component has been omitted. A joint strategy $(u_1, \dots, u_{i - 1}, v_i, u_{i + 1}, \dots, u_N) \in U$ is written as $(v_i,u_{-i})$ \msa{---}{} with no change of the ordering.}{}
\begin{defn}[cf. \cite{bib:NashNPerson}]\label{defn:nep}
A point $u \in U_\ad$ is called a \emph{Nash equilibrium}, if for all $i = 1, \dots, N$ the inclusion
\begin{eq*}
u_i \in \argmin\left\{\calJ_i(v_i,u_{-i}) \text{ subject to } v_i \in U_\ad^i\right\}
\end{eq*}
holds. The problem of finding such a point is called a \emph{Nash equilibrium problem} (abbr.: NEP).
If, moreover, the sets $U_\ad^i$ are convex and the objectives $v_i \mapsto \cJnu(v_i, u_{-i})$ are convex on $U_\ad^i$, the \nep{} is called \emph{convex}.
\end{defn}
\msa{Nash equilibrium problems are an important concept to model a competition between entities. Therefore, we also refer to NEPs as \emph{Nash games}.}{}
A possible generalization of this concept with practical relevance is the case, where the strategies of the players influence the set of feasible choices for each of the other players. This mechanism can be modeled by the introduction of a \emph{\msa{strategy}{} mapping} $C_i : U_\ad^{-i} \rra U_\ad^i$ for each player. The combination of all of them $C: U_\ad \rra U_\ad$ with $C(u) := C_1(u_{-1}) \times \dots \times C_N(u_{-N})$ is called the \emph{\msa{joint}{} strategy mapping}. Subsequently, it is possible to generalize the \msa{concept given in \refer{defn:nep}}. 
\begin{defn}\label{defn:gnep}
A point $u \in U_\ad$ is called a \emph{(generalized) Nash equilibrium}, if for all $i = 1, \dots, N$ holds
\begin{eq*}
u_i \in \argmin\left\{\calJ_i(v_i,u_{-i}) \text{ subject to } v_i \in C_i(u_{-i})\right\}.
\end{eq*}
The problem of finding such a point is called a \emph{generalized Nash equilibrium problem} (abbr.: \gnep).
If, moreover, the sets $C(u)$ are convex for all choices $u \in \dom{C}$ and the objectives $v_i \mapsto \cJnu(v_i, u_{-i})$ are convex on $\range{C_i}$, the \gnep{} is called \emph{convex}. 
\end{defn}
For a better understanding of \msa{the problems}{} under investigation, we propose the \msa{subsequent}{} example:
\begin{ex}\label{ex:intro}
Consider the following \mh{g}eneralized Nash \mh{e}quilibrium \mh{p}roblem governed by a semi-linear elliptic PDE with distributed control\mh{:}{}
\begin{eq*}
&\text{minimize } \calJ_i(u_i,u_{-i}) := \frac{1}{2}\int_{\omega_i} (y-y_d^i)^2 \dx + \frac{\alpha}{2}\int_{\omega_i} u_i^2 \dx \text{ over } u_i \in L^2(\omega_i)\\
&\text{subject to } a_i \leq u_i \leq b_i \text{ a.e. on } \omega_i, \ \underline{\psi} \leq y \leq \overline\psi \text{ a.e. on } \Omega \text{ and}\\
&\hspace*{4.31em}-\Delta y + y^3 = \sum_{i = 1}^4 B_i u_i \text{ in } \Omega,\ y = 0 \text{ on } \partial \Omega.
\end{eq*}
\end{ex}
\noindent Here, $a_i, b_i \in L^2(\Omega)$ and $\underline{\psi}, \overline{\psi} \in H^2(\Omega)$ with $\operatorname{tr}_{\partial \Omega}(\underline\psi) < 0 < \operatorname{tr}_{\partial \Omega}(\overline\psi)$ a.e. on $\partial \Omega$ are given parameters and $B_i$ denotes the extension of the strategy $u_i$ to \msa{the whole domain}{} $\Omega$ by the value zero. In this example, neither the functional, nor the feasible set contain an explicit dependence on the controls of the other players. Instead, the coupling takes place via the state as the solution of the PDE leading to a spatial coupling between the control regions $\omega_i$.\\
Often\ms{, as demonstrated in \refer{ex:intro},}{} the players' decisions are coupled via a single condition that needs to be fulfilled by all players simultaneously. In this case, the strategy mapping can be modeled via the feasible set \msa{of}{} all players.
\begin{defn}\label{defn:gnep:jointly_constrained}
A \gnep{} is called \emph{jointly constrained} or is said to have \emph{shared constraint\mh{s}}\mh{,}{} if there exists a set $\calF$ such that \emph{Rosen's Law} (cf. \cite[p. 484]{bib:AusselCoreaMarechal2011})
\begin{eq*}
v_i \in C_i(u_{-i}) \Leftrightarrow (v_i, u_{-i}) \in \calF
\end{eq*}
holds. If\mh{,}{} moreover, a \gnep{} is convex and jointly constrained and the set $\calF$ is convex as well, then the \gnep{} is called \emph{jointly convex}. 
\end{defn}
In the jointly constrained case it is worthwhile investigating a \ms{modified solution concept}{} based on the previously introduced joint constraint set $\calF$.
\begin{defn}[see also \cite{bib:Rosen}]\label{defn:vep}
A point $u \in U_\ad$ is called a \emph{variational equilibrium}, if for all $i = 1, \dots, N$ \mh{it}{} holds \mh{that}
\begin{eq*}
u \in \argmin\left\{\sum_{i = 1}^N \calJ_i(v_i,u_{-i}) \text{ subject to } v \in \calF \right\}.
\end{eq*}
The problem of finding such a point is called a \emph{variational equilibrium problem} \msa{(abbr.: \vep)}.
\end{defn}

The relationship between these concepts will be investigated next. For this sake, the concept of the \emph{Nikaido-Isoda functional} (cf. \cite{bib:NikaidoIsoda}) is introduced.
\begin{defn}\label{defn:nikaido_isoda}
Let a family of functionals $\calJ_i: U \rightarrow \R$ on a Banach space $U$ be given. The \emph{Nikaido--Isoda} functional $\Psi: U \times U \rightarrow \R$ is defined as
\begin{eq*}
\Psi(u,v) := \sum_{i = 1}^N \left( \cJnu(u_i, u_{-i}) - \cJnu(v_i, u_{-i}) \right).
\end{eq*}
\end{defn}
\begin{thm}[\ms{compare to \cite[Lemma 3.1]{bib:NikaidoIsoda} and \cite[Theorems 2.1.2 and 2.2.3]{bib:vonHeusinger}}]\label{thm:nikaido_isoda:merit}
Consider the following \mh{value}{} functions \mh{associated with}{} the Nikaido--Isoda functional:
\begin{eq*}
V(u) = \sup_{v \in C(u)} \Psi(u,v)\mh{,} \text{ and }\, \Vhat(u) = \sup_{v \in \mathcal{F}} \Psi(u,v)\mh{,}{}
\end{eq*}
corresponding to the \msa{\gnep}{} respectively \msa{\vep}{}. Then $V(u) = 0$ holds, if and only if $u$ is a Nash equilibrium, and $\Vhat(u) = 0$, if and only if $u$ is a variational equilibrium.
\end{thm}
\begin{proof}
For the proof we refer to references given in \refer{thm:nikaido_isoda:merit} to \cite[Theorem 2]{bib:HintStenglEquiGamma} in combination with the discussion given in Subsection 2.2 therein.
\end{proof}

From the definition of the \mh{value}{} functions $V$ and $\Vhat$ one sees that the condition $v \in C(u)$ is substituted by the feasible set and hence is restricted to the fixed point set of $C$. The relationship between these two solution concepts is discussed in the following \msa{theorem}.
\begin{thm}\label{thm:vep->gnep}
In a jointly constrained \gnep{} every variational equilibrium is a Nash equilibrium.
\end{thm}
\begin{proof}
Let $u \in \calF$ be a variational equilibrium. Then $u \in C(u)$, and selecting an arbitrary $v_i \in C_i(u_{-i})$ and setting $\tilde v = (v_i, u_{-i}) \in \calF$ yields
\begin{eq*}
\calJ_i(u_i, u_{-i}) - \calJ_i(v_i, u_{-i}) = \Psi(u,\tilde v) \leq \Vhat(u) = 0
\end{eq*}
and hence $\calJ_i(u_i, u_{-i}) \leq \calJ_i(v_i, u_{-i})$.
\end{proof}
Within the scope of this work we do not address existence of equilibria explicitly. We just refer to \cite{bib:Dutang} for the finite dimensional case and \msa{to \cite{bib:HintStenglEquiGamma}, where this issue is discussed}. It is just \msa{mentioned here}, that the presence of \makebox{(quasi-)}convexity in the objective plays a crucial role for providing existence.\\
In the following, we do not only want to discuss the Nash equilibrium problem, but also want to discuss \msa{the variational equilibrium problem}{} in the \msa{presence}{} of shared constraints. \msa{Hence, we discuss both of them in parallel in the scope of our investigation.}{} As the proofs turn out to be very similar, we proceed \msa{most of the time}{} with the statement formulated for both cases, but only \msa{proven for the}{} Nash equilibrium \msa{problem}.

\subsection{First Order Conditions}
After setting up the basic definitions we proceed with the characterization of solutions. For this sake the optimization based structure is used to derive first order \mh{optimality}{} conditions in the form of a Lagrange-multiplier system for each player. The combined system serves as a first order necessary condition for a point to be a Nash equilibrium. Therefore we want to \ms{specify}{} the framework used in the sections and propose the following \msa{set of}{} assumptions.

\begin{ass}\label{ass:gnep_vep:strategy_map}\hfill%
\begin{enum}
\item\label{enum:part:gnep:strategy_map} Let the strategy mapping $C: U_\ad \rra U_\ad$ admit the form
\begin{eq*}
C_i(u_{-i}) = \left\{ u_i \in U_\ad^i : g_i(u_i,u_{-i}) \in K_i \right\}
\end{eq*}
involving a family of \mh{Fr\'echet differentiable}{} mappings $g_i: U \rightarrow X_i$ together with $X_i$ a Banach space and $K_i \subseteq X_i$ a non-empty, closed, convex cone.
\item\label{enum:part:vep:strategy_map} Let in the above Assumption (i) the mappings $g_i$ coincide with a \mh{Fr\'echet differentiable}{} $g:U \rightarrow X$ and let the constraint sets coincide with a $K \subseteq X$. The \emph{joint constraint set} of the resulting jointly constrained Nash game is denoted as
\begin{eq*}
\calF := \left\{ u \in U_\ad : g(u) \in K \right\}.
\end{eq*}
\end{enum}
\end{ass}
\msa{The part \refer{ass:gnep_vep:strategy_map}\ref{enum:part:gnep:strategy_map} corresponds solely to the case of a jointy constrained GNEP.}{} In view of Nash games involving partial differential equations \msa{in the sense of}{} \refer{ex:intro}, \mh{the operator $g$ can be identified with a composition $g = \calG \circ S$, where $S$ denotes the solution mapping of an operator equation and $\calG$ \msa{is}{} another mapping establishing the state constraint in terms of a cone condition, see also \refer{sec:example}.}{} \mh{Having this consideration in mind,}{} we will in the following refer to the conditions $g_i(u) \in K_i$ respectively $g(u) \in K$ as \emph{state constraint(s)}.
\msa{Next, we derive a first order necessary condition for Nash equilibria}.  

\begin{thm}\label{thm:gnep:firstorder}
\msa{Consider a Nash equilibrium $u \in U_\ad$. Let \refer{ass:gnep_vep:strategy_map}\ref{enum:part:gnep:strategy_map} hold and assume the objectives $v_i \mapsto \cJnu(v_i,u_{-i})$ to be Fr\'echet differentiable.}{} Moreover, \msa{let}{} the following $\mathrm{(RZK)}$-condition
\begin{eq}[\tag{$\mathrm{RZK}_{\text{Nash}}$}]\label{eq:gnep:rzk}
\partial_i g_i(u)U_\ad^i(u_i) - K_i(g_i(u)) = X_i \quad \text{ for all } i = 1, \dots, N
\end{eq}
\msa{hold.}{} Then, there exist Lagrange multipliers $\mu_i \in X_i^*$, such that the complementarity system
\begin{sys}[\label{sys:gnep:firstorder}]
\hspace*{6ex}0 &= q_i - \partial_i g_i(u)^* \mu_i + \lambda_i &&\hspace*{-12ex}\text{ in } U_i^*,\hspace*{6ex}\\
\hspace*{6ex}q_i &= \partial_i \calJ_i(u) &&\hspace*{-12ex}\text{ in } U_i^*,\hspace*{6ex}\\
\hspace*{6ex}K_i^+ \ni \mu_i &\perp g_i(u) \in K_i,  &&\hspace*{-6ex}\\
\hspace*{6ex}\lambda_i &\in N_{U_\ad^i}(u_i) &&\hspace*{-12ex}\text{ in } U_i^*\hspace*{6ex}
\end{sys}
is fulfilled.
\end{thm}
\begin{proof} In essence we utilize \cite{bib:ZoweKurcyusz}.\\
Since $u \in U_\ad$ is a Nash equilibrium every player solves the minimization problem
\begin{eq*}
&\text{minimize} &&\calJ_i(v_i, u_{-i}) \text{ over } \msa{v_i}{} \in U^i_\ad,\\
&\text{subject to} &&g_i(v_i, u_{-i}) \in K_i.
\end{eq*}
The strategy $u_i$ is the minimizer of this problem. The constraint qualification reads as $\partial_i g_i(u) U_\ad^i(u_i) - K_i(g_i(u)) = X_i$ and leads by \cite[Theorem 4.1]{bib:ZoweKurcyusz} to the existence of a Lagrange multiplier $\mu_i \in K_i^+$ such that the system 
\begin{sys*}
0 &= \langle \mu_i, g_i(u) \rangle,\\
0 &\in \partial_i \calJ_i(u_i, u_{-i}) - \partial_i g_i(u)^* \mu_i + N_{U_\ad^i}(u_i)
\end{sys*}
is fulfilled. Combining the collective first order systems of each player together with defining $q_i = \partial_i \calJ_i(u)$ and setting $\lambda_i = \partial_i g_i(u)^* \mu_i - q_i \in N_{U_\ad^i}(u_i)$ leads to \refer{sys:gnep:firstorder}.
\end{proof}

The discussion of the first order system for \msa{variational equilibria}{} follows the same pattern and is addressed in the following theorem.

\begin{thm}\label{thm:vep:firstorder}
\msa{Consider a variational equilibrium $u \in U_\ad$. let \refer{ass:gnep_vep:strategy_map}\ref{enum:part:vep:strategy_map} hold and assume the objectives $v_i \mapsto \cJnu(v_i,u_{-i})$ to be Fr\'echet differentiable. Moreover, let the following $\mathrm{(RZK)}$-condition
\begin{eq}[\tag{$\mathrm{RZK}_{\text{Var}}$}]\label{eq:vep:rzk}
Dg(u)U_\ad(u) - K(g(u)) = X
\end{eq}
hold.}{} Then, there exists a Lagrange multiplier $\mu \in X^*$, such that the complementarity system
\begin{sys}[\label{sys:vep:firstorder}]
\hspace*{6ex}0 &= q_i - \partial_i g_i(u)^* \mu + \lambda_i &&\hspace*{-12ex}\text{ in } U_i^*,\hspace*{6ex}\\
\hspace*{6ex}q_i &= \partial_i \calJ_i(u) &&\hspace*{-12ex}\text{ in } U_i^*,\hspace*{6ex}\\
\hspace*{6ex}K^+ \ni \mu &\perp g(u) \in K,  &&\hspace*{-6ex}\\
\hspace*{6ex}\lambda_i &\in N_{U_\ad^i}(u_i) &&\hspace*{-12ex}\text{ in } U_i^*\hspace*{6ex}
\end{sys}
is fulfilled.
\end{thm}
\begin{proof}
The proof is essentially the same as the one of \refer{thm:gnep:firstorder} \msa{and therefore omitted}{}.
\end{proof}

To deepen our understanding of the interconnection between \msa{Nash}{} and variational equilibria, we compare the two systems derived in \msa{\refer{thm:gnep:firstorder} and \refer{thm:vep:firstorder}:}{} Instead of one multiplier
for each player in \refer{sys:gnep:firstorder}, we only have a single multiplier for all players arising form the joint state constraint.
\begin{rem}\label{rem:gnep_vep:rzk}
For the jointly constrained case \mh{it}{} holds:
\eqref{eq:gnep:rzk} $\Rightarrow$ \eqref{eq:vep:rzk}.\\
Let \msa{$u \in U_\ad$ be a Nash equilibrium and assume}{} \eqref{eq:gnep:rzk} \msa{to}{} be fulfilled. Choose $x \in X$ arbitrarily. By assumption there exist $\alpha_i, \beta_i \geq 0$ together with $v_i \in U_\ad^i$ and $k_i \in K$ such that 
\begin{align*}
\alpha_i \partial_i g(u)(v_i - u_i) - (k_i - \beta_i g(u)) = x \text{ for all } i = 1, \dots, N.
\end{align*}
Without loss of generality we assume $\alpha_i > 0$. Otherwise we can as well write
\begin{align*}
\partial_i g(u)(u_i -  u_i) - (k_i - \msa{\beta_i}{} g(u)) = x.
\end{align*}
Thus the above equation is fulfilled with $v_i = u_i$ and $\alpha_i = 1$.\\
Set $\alpha := \left(\sum_{i = 1}^N \frac{1}{\alpha_n}\right)^{-1}$ and multiply each equation with ${\alpha_i}^{-1}\alpha$. By \msa{subsequent}{} addition we obtain
\begin{align*}
\sum_{i = 1}^N \partial_i g(u)\alpha (v_i - u_i)  - \brackets{2}{\sum_{i = 1}^N k_i - \brackets{2}{\sum_{i = 1}^N \frac{\alpha \beta_i}{\alpha_i}} g(u)} = \brackets{2}{\sum_{i = 1}^N \frac{\alpha}{\alpha_i}}x = x.
\end{align*}
Since $K$ is a convex cone, it holds $k:= \sum_{i = 1}^N k_i \in K$ and with $\beta := \sum_{i = 1}^N \frac{\alpha \beta_i}{\alpha_i} > 0$ we get
\begin{align*}
\alpha Dg(u)(v- u)- (k - \beta g(u)) = x.
\end{align*}
Thus, \eqref{eq:vep:rzk} is fulfilled.\\
This observation tells us, that we have traded a stronger solution concept with a \msa{potentially}{} weaker constraint qualification. A similar observation has been made in \cite{bib:HintSurKaemmler} for the Slater condition proposed therein.
\end{rem}
\msa{\section{Penalization and Path-Following Technique}}{}\label{sec:penalization}
\msa{In the previous section we introduced and discussed the Nash equilibrium problem.
One of the major analytic challenges is the presence of the other players decision in the constraint set.
This has been partially addressed by the introduction of the variational equilibrium problem in the jointly constrained case.
In view of numerical calculations however one may be interested in performing an approximation.
For this sake, we relax the state constraint using a penalty functional and leave the control constraint realized by $U_\ad^i$ unchanged.
Again, we will derive first order necessary conditions. Subsequently, a convergence analysis not only regarding the solution of the penalized problems but also the multipliers occurring in the first order system is provided together with a practically applicable update strategy for the underlying penalty parameter.}{}
\msa{\subsection{Definition of Penalized Equilibrium Problems}}{}
\msa{Returning to the aforementioned dependence of the constraint set, the idea behind a penalty approach is the relaxation of the constraint and shifting}{} the dependence \msa{into the objectives}. For this sake, we use the approach derived in \cite{bib:HintSurKaemmler} \msa{and mentioned in}{} \cite[Sections 3.1 and 3.2]{bib:HintStenglEquiGamma} and introduce a penalty functional. Therefore, we introduce the following two problems:
\begin{defn}\label{defn:gnep_vep:penalized}\makebox{}
\begin{enum}
\item\label{enum:part:gnep:penalized} Let  \refer{ass:gnep_vep:strategy_map}\ref{enum:part:gnep:strategy_map} hold and let convex, \lsc{} penalty functionals $\beta_i: X_i \rightarrow [0,+\infty)$ with $\beta_i(x_i) = 0$, if and only if $x_i \in K_i$ be given together with a \emph{penalty parameter} $\gamma > 0$. Define the \emph{penalized Nash equilibrium problem} \pengnep{} as
\begin{eq}[\tag{\ensuremath{\mathrm{GNEP}_\gamma}}]\label{eq:gnep:penalized}
\hspace*{1ex}u_i \in \argmin\left\{ \cJnu(v_i, u_{-i}) + \gamma \beta_i(g_i(v_i, u_{-i})) \text{ subject to } v_i \in U_\ad^i \right\}
\end{eq}
for all $i = 1, \dots, N$.
\item\label{enum:part:vep:penalized} Let  \refer{ass:gnep_vep:strategy_map}\ref{enum:part:vep:strategy_map} hold and let a convex, \lsc{} penalty functional $\beta: X \rightarrow [0,+\infty)$ with $\beta(x) = 0$, if and only if $x \in K$ be given together with a \emph{penalty parameter} $\gamma > 0$. Define the \emph{penalized variational equilibrium problem} \penvep{} as
\begin{eq}[\tag{\ensuremath{\mathrm{VEP}_\gamma}}]\label{eq:vep:penalized}
u \in \argmin\left\{ \sum_{i = 1}^N \cJnu(v_i, u_{-i}) + \gamma \beta(g(v)) \text{ subject to } v \in U_\ad \right\}.
\end{eq}
\end{enum}
\end{defn}
In the light of \cite{bib:HintStenglEquiGamma}, we defined the penalty functionals 
\begin{eq*}
\pi_C(v,u) := \sum_{i = 1}^N \beta_i(g_i(v_i,u_{-i}))  \text{ and } \pi_\calF(v) := \beta(g(v))
\end{eq*}
\msa{using the notation given therein.}{}
To meet the conditions for the $\Gamma$-convergence discussed therein, it is sufficient to assume the weak continuity of $v_i \mapsto g_i(v_i,u_{-i})$ and $g$, respectively and use the convexity of $\beta_i$, \msa{respectively}{} $\beta$ to obtain the weak lower semi-continuity and subsequent of its composition.\\
\msa{\subsection{First Order Conditions for the Penalized Equilibrium Problems}}{}
\msa{One way to achieve the convergence of a sequence of penalized equilibrium problems has been addressed in \cite{bib:HintStenglEquiGamma}. Therein, $\Gamma$-convergence has been employed to achieve this goal. Within the scope of this work however, we seek to obtain a more detailed view and establish our convergence theory via the first order system. Hence, of}{} particular interest is the approximation of the multipliers in \refer{sys:gnep:firstorder} (respectively \refer{sys:vep:firstorder}). We proceed \msa{by the derivation of}{} the first order systems for the penalized problems. 

\begin{thm}\label{thm:gnep:penalized:firstorder}
Let \refer{ass:gnep_vep:strategy_map}\ref{enum:part:gnep:strategy_map} hold and let the objectives $U_\ad^i \ni v_i \mapsto \cJnu(v_i,u_{-i}) \in \R$ as well as the penalty mappings $\beta_i: X_i \rightarrow [0,+\infty)$ and the operators $U_i \ni u_i \mapsto g_i(u_i,u_{-i}) \in X_i$ be Fr\'echet differentiable in their respective spaces.\\
Then, every Nash equilibrium $u^\gamma\in U_\ad$ of \pengnep{} fulfils the following first order system 
\begin{sys}[\label{sys:gnep:penalized:firstorder}]
\hspace*{7.5ex}0 &= q_i^\gamma - \partial_i g_i(u^\gamma)^* \mu_i^\gamma + \lambda_i^\gamma &&\hspace*{-15ex}\text{ in } U_i^*,\hspace*{7.5ex}\\
\hspace*{7.5ex}q_i^\gamma &= \partial_i \calJ_i(u^\gamma) &&\hspace*{-15ex}\text{ in } U_i^*,\hspace*{7.5ex}\\
\hspace*{7.5ex}\mu_i^\gamma &= -\gamma D\beta_i(g_i(u^\gamma)) &&\hspace*{-15ex}\text{ in } X_i^*,\hspace*{7.5ex}\\
\hspace*{7.5ex}\lambda_i^\gamma &\in N_{U_\ad^i}(u_i^\gamma) &&\hspace*{-15ex}\text{ in } U_i^*\hspace*{7.5ex}
\end{sys}
for all $i = 1, \dots, N$.
\end{thm}
\begin{proof}
Let $u^\gamma \in U_\ad$ be a Nash equilibrium of \pengnep. Then, by \refer{defn:gnep_vep:penalized}\ref{enum:part:gnep:penalized} $u_i^\gamma$ is a solution of
\begin{eq*}
u_i^\gamma \in \argmin_{v_i \in U_\ad^i}\left\{ \cJnu(v_i,u_{-i}^\gamma) + \gamma \beta_i(g_i(v_i,u_{-i}^\gamma)) \right\}.
\end{eq*}
Using the Fr\'echet differentiability \msa{this}{} yields the following first order condition
\begin{eq*}
0 \in \partial_i \cJnu(u_i^\gamma,u_{-i}^\gamma) + \gamma \partial_i g_i(u^\gamma)^* D\beta_i(g_i(u^\gamma)) + N_{U_\ad^i}(u_i^\gamma) \text{ in } U_i^*.
\end{eq*}
Setting $q_i^\gamma = \partial_i \cJnu(u_i^\gamma,u_{-i}^\gamma)$ and $\mu_i^\gamma = -\gamma D\beta_i(g_i(u^\gamma))$ there exists a $\lambda_i^\gamma \in N_{U_\ad^i}(u_i^\gamma)$ with
\begin{eq*}
0 = q_i^\gamma - \partial_i g_i(u^\gamma)^* \mu_i^\gamma + \lambda_i^\gamma \text{ in } U_i^*
\end{eq*}
proving the assertion.
\end{proof}

Comparing the system \msa{in \refer{thm:gnep:penalized:firstorder}}{} with the one \msa{given in \refer{thm:gnep:firstorder}}{}, we observe that the role of the multiplier for the state constraint is (up to a sign) fulfilled by the derivative of the penalty function scaled by the penalty parameter. \msa{Next}, we discuss the first order system for the penalized variational equilibrium problem. 
\begin{thm}\label{thm:vep:penalized:firstorder}
Let \refer{ass:gnep_vep:strategy_map}\ref{enum:part:vep:strategy_map} hold and let the objectives $U_\ad^i \ni v_i \mapsto \cJnu(v_i,u_{-i}) \in \R$ as well as the penalty mapping $\beta: X \rightarrow [0,+\infty)$ and the operator $U \ni u \mapsto g(u_i,u_{-i}) \in X$ be Fr\'echet differentiable.\\
Then, every equilibrium $u^\gamma\in U_\ad$ of \penvep{} fulfils the following first order system 
\begin{sys}[\label{sys:vep:penalized:firstorder}]
\hspace*{7.5ex}0 &= q_i^\gamma - \partial_i g(u^\gamma)^* \mu^\gamma + \lambda_i^\gamma &&\hspace*{-15ex}\text{ in } U_i^*,\hspace*{7.5ex}\\
\hspace*{7.5ex}q_i^\gamma &= \partial_i \calJ_i(u^\gamma) &&\hspace*{-15ex}\text{ in } U_i^*,\hspace*{7.5ex}\\
\hspace*{7.5ex}\mu^\gamma &= -\gamma D\beta(g(u^\gamma)) &&\hspace*{-15ex}\text{ in } X^*,\hspace*{7.5ex}\\
\hspace*{7.5ex}\lambda_i^\gamma &\in N_{U_\ad^i}(u_i^\gamma) &&\hspace*{-15ex}\text{ in } U_i^*\hspace*{7.5ex}
\end{sys}
for all $i = 1, \dots, N$.
\end{thm}
\begin{proof}
Let $u^\gamma \in U_\ad$ be an equilibrium of \penvep{} (see \refer{defn:gnep_vep:penalized}\ref{enum:part:vep:penalized}), which means   $u^\gamma$ fulfills
\begin{eq*}
u^\gamma \in \argmin_{v \in U_\ad^i}\left\{ \sum_{i = 1}^N \cJnu(v_i,u_{-i}^\gamma) + \gamma \beta(g(v)) \right\}.
\end{eq*}
Using the Fr\'echet differentiability in combination with the product structure of the space $U$ \msa{and}{} the set $U_\ad$ yields the following first order condition
\begin{eq*}
0 &\in D\left( \sum_{i = 1}^N \cJnu((\ccdot)_i ,u_{-i}^\gamma) + \gamma\beta(g(\ccdot))\right)(u^\gamma) + N_{U_\ad}(u^\gamma)\\
&= (\partial_1 \calJ_1(u^\gamma), \dots, \partial_N \calJ_N(u^\gamma)) + \gamma Dg(u^\gamma)^* D\beta(g(u^\gamma)) + N_{U_\ad}(u^\gamma)\\
&= \prod_{i = 1}^N \left(\partial_i \cJnu(u^\gamma) + \gamma \partial_i g(u^\gamma)^* D\beta(g(u^\gamma)) + N_{U_\ad^i}(u_i^\gamma) \right) \text{ in } U^*.
\end{eq*}
Setting $q_i^\gamma:= \partial \cJnu(u_i,u_{-i})$ with $\mu^\gamma := -\gamma D\beta(g(u^\gamma))$ there exist $\lambda_i^\gamma \in N_{U_\ad^i}(u_i^\gamma)$ with
\begin{eq*}
0 &= q_i^\gamma - \partial_i g(u^\gamma)^* \mu^\gamma + \lambda_i^\gamma \text{ in } U_i^*
\end{eq*}
proving the assertion.
\end{proof}

We observe that both---the penalized, \msa{jointly constrained}{} Nash game and the penalized variational
equilibrium problem---yield the same first order system. Hence, one might ask as
well for the relation between the two concepts. To gain insight into this question we provide in analogy to \refer{thm:vep->gnep} the following \msa{result}.
\begin{thm}
Let \refer{ass:gnep_vep:strategy_map}\ref{enum:part:vep:strategy_map} hold. Then every equilibrium of \penvep{} is an equilibrium of \pengnep.
\end{thm}
\begin{proof}
Let $u^\gamma \in U_\ad$ be a variational equilibrium of \pengnep. Take $\msa{j} \in \{1, \dots, N\}$ and $v_j \in U^j_\ad$ arbitrarily and set $\tilde v = (v_j,u^\gamma_{-j})$. \msa{Then}, we deduce
\begin{eq*}
\sum_{i = 1}^N \cJnu(u^\gamma_i, u^\gamma_{-i}) &+ \gamma \beta(g(u^\gamma)) \leq \sum_{i = 1}^N \cJnu(v_i, u^\gamma_{-i}) + \gamma \beta(g(\tilde v))\\
&= \sum_{i \neq j} \cJnu(u^\gamma_i,u^\gamma_{-i}) + \cJnu(v_j,u^\gamma_j) + \gamma \beta(g(v_j,u^\gamma_{-j}))
\end{eq*}
and hence
\begin{eq*}
\calJ_j(u^\gamma_j,u^\gamma_{-j}) + \gamma \beta(g(u^\gamma_j,u^\gamma_{-j})) \leq \calJ_j(v_j,u^\gamma_{-j}) + \gamma \beta(g(v_j,u^\gamma_{-j}))
\end{eq*}
proving $u^\gamma$ to be a Nash equilibrium for \pengnep{} since the choice of \msa{$j$}{} was arbitrary.
\end{proof}
In fact, in the case of a convex game with $g$ being a concave vector-valued operator as in \cite{bib:HintStenglKConvex} and $\beta$ being isotone, we obtain the convexity of the optimization problem and hence the first order system being a sufficient condition.

\msa{\subsection{Asymptotic Behaviour of the Primal-Dual Path}}
\msa{After the derivation of the first order systems of the penalized equilibrium problems we investigate the behaviour of equilibria and  multipliers as the penalty parameter $\gamma$ goes to infinity. The naive expectation is the convergence of the multipliers in \refer{sys:gnep:penalized:firstorder} and \refer{sys:vep:penalized:firstorder} towards the corresponding ones in \refer{sys:gnep:firstorder} and \refer{sys:vep:firstorder} respectively. In fact, $\Gamma$-convergence of the equilibrium problems has been confirmed in \cite{bib:HintStenglEquiGamma} but this does not show the behavior of the multipliers yet. For the investigation of the latter we propose the following definition first.}
\begin{defn}
Define the set-valued solution operator
\begin{eq*}
\calS : (0,\infty) \rra  U \times U^* \times U^* \times \textstyle\prod_{i = 1}^N X_i,\ \gamma \mapsto (u^\gamma, q^\gamma, \lambda^\gamma, \mu^\gamma)
\end{eq*}
by mapping $\gamma$ to the set of solutions of \refer{sys:gnep:penalized:firstorder}. The set of \emph{primal-dual paths} $\mathscr{P}$ is defined as
\begin{eq*}
\mathscr{P} := \left\{\rule{0pt}{1em} ((u^\gamma, q^\gamma, \lambda^\gamma,  \mu^\gamma))_{\gamma > 0} : (u^\gamma, q^\gamma, \lambda^\gamma,  \mu^\gamma) \in \calS(\gamma) \right\}
\end{eq*}
and every of its elements is called a \emph{(primal-dual) path} for \refer{sys:gnep:penalized:firstorder}.\\
Define the set-valued solution operator
\begin{eq*}
\widehat\calS : (0,\infty) \rra  U \times U^* \times U^* \times X,\ \gamma \mapsto (u^\gamma, q^\gamma, \lambda^\gamma, \mu^\gamma)
\end{eq*}
by mapping $\gamma$ to the set of solutions of  \refer{sys:vep:penalized:firstorder}. Its set of \emph{primal-dual paths} $\widehat{\mathscr{P}}$ is defined as
\begin{eq*}
\widehat{\mathscr{P}} := \left\{ ((u^\gamma, q^\gamma, \lambda^\gamma,  \mu^\gamma))_{\gamma > 0} : (u^\gamma, q^\gamma, \lambda^\gamma,  \mu^\gamma) \in \widehat \calS(\gamma) \right\}.
\end{eq*}
\end{defn}
Next, we are investigating the asymptotic behaviour of the multipliers.
\begin{lem}\label{lem:gnep_vep:path:bdd}
Under the assumptions of \refer{thm:gnep:penalized:firstorder} (respectively \refer{thm:vep:penalized:firstorder}) assume additionally that for all $i = 1, \dots, N$ the mappings $u \mapsto \partial_i \calJ_i(u_i,u_{-i})$ and $u \mapsto \partial_i g_i(u)$ (respectively $u \mapsto D g(u)$) are bounded (i.e.: images of bounded sets are bounded).\\
Moreover, let the following \emph{uniform Robinson-type condition} hold:\\
There exists $\varepsilon > 0$ such that for all $u^\gamma$ in the path \mh{it}{} holds \mh{that}{}
\begin{eq}\label{eq:gnep:path:uniformrobinson}
\varepsilon \B_{X_i} &\subseteq \partial_i g_i(u^\gamma)(U_\ad^i - u_i^\gamma) - (K_i - g_i(u^\gamma)) \text{ for all } i = 1, \dots, N
\end{eq}
\begin{eq}\label{eq:vep:path:uniformrobinson}
\left(\text{respectively}\hspace*{1em} \varepsilon \B_X \right. &\left. \subseteq D g(u^\gamma)(U_\ad - u^\gamma) - (K - g(u^\gamma)) \right).
\end{eq}
Then, the path $\mathscr{P}$ (respectively $\widehat{\mathscr{P}}$) is bounded.
\end{lem}
\begin{proof}(similar to \cite{bib:HintSurKaemmler}, inspired by \cite{bib:HintKunischPath})\\
The boundedness of the sets $U_\ad^i$ for all $i = 1, \dots, N$ yields $\|u^\gamma\|_U \leq C_U$
for some constant $C_U > 0$. \msa{The}{} assumed boundedness of the operators $\partial_i\cJnu$ \msa{yields}{} the boundedness of $q_i^\gamma$ for all $i = 1, \dots, N$. \msa{By}{} the convexity of the penalty function $\beta_i$ and the set $K_i$ we obtain for all $z_i \in K_i$, that
\begin{eq*}
0 = \gamma \beta_i(z_i) \geq \gamma \beta_i(g_i(u^\gamma)) + \langle \gamma D\beta_i(g_i(u^\gamma)) , z_i - g_i(u^\gamma) \rangle \geq -\langle \mu^\gamma_i , z_i - g_i(u^\gamma) \rangle.
\end{eq*}
To prove the boundedness of $\mu^\gamma_i$ we utilize condition \eqref{eq:gnep:path:uniformrobinson}. Hence for all $\gamma > 0$ and $x_i \in X_i$ with $\|x_i\|_{X_i} \leq \varepsilon$ there exist $v_i^\gamma \in U_\ad^i$ as well as $k_i^\gamma \in K_i$ fulfilling
\begin{eq*}
x_i = \partial_i g_i(u^\gamma)(v_i^\gamma - u_i^\gamma) - (k_i^\gamma - g_i(u^\gamma)).
\end{eq*}
Applying $\mu_i^\gamma$ yields
\begin{eq*}
\langle \mu_i^\gamma, x_i \rangle_{\msa{X_i^*,X_i}} &= \langle \mu_i^\gamma, \partial_i g_i(u^\gamma)(v_i^\gamma - u^\gamma_i) - (k_i^\gamma - g_i(u^\gamma)) \rangle_{\msa{X_i^*,}X_i}\\
&= \langle \partial_i g_i(u^\gamma)^* \mu_i^\gamma, v_i^\gamma - u_i^\gamma \rangle_{\msa{U_i^*,}U_i} - \langle \mu_i^\gamma, k_i^\gamma - g_i(u^\gamma) \rangle_{\msa{X_i^*,}X_i}\\
&\leq \langle q_i^\gamma + \lambda_i^\gamma, v_i^\gamma - u_i^\gamma \rangle_{\msa{U_i^*,}U_i} \leq \langle q_i^\gamma, v_i^\gamma - u_i^\gamma \rangle_{\msa{U_i^*,U_i}}.
\end{eq*}
The last expression is bounded by some constant $C$ due to the boundedness of $U_\ad$ and the boundedness of $q^\gamma_i$.
Hence, we obtain
\begin{eq*}
\|\mu_i^\gamma\|_{X_i^*} = \frac{1}{\varepsilon}\sup_{\|x_i\|_{X_i} \leq \varepsilon} \langle \mu_i^\gamma , x_i \rangle_{\msa{X_i^*,X_i}} \leq  \frac{C}{\varepsilon}
\end{eq*}
and conclusively also $\lambda_i^\gamma = -\left( q_i^\gamma - \partial_i g_i (u^\gamma)^* \mu_i^\gamma \right)$ is bounded.
\end{proof}
Establishing a uniform Robinson condition implies the RZK condition for each $\gamma$, but it is not clear, whether the opposite holds true, because \cite[Theorem 2.1]{bib:ZoweKurcyusz} just guarantees the existence of a ball, but does not specify its size. Hence, the radius might vanish \msa{in the limit}{} as the penalty parameter goes to infinity.
\begin{rem}\label{rem:gnep_vep:uniformrobinson}
In analogy to \refer{rem:gnep_vep:rzk} we observe that also for the uniform Robinson-type conditions in \refer{lem:gnep_vep:path:bdd} in the jointly constrained case \mh{it}{} holds that  \refer{eq:gnep:path:uniformrobinson} implies \refer{eq:vep:path:uniformrobinson}:\\
Since $U^i_\ad$ is assumed to be convex, for all $\gamma>0$ and $\lambda \geq 1$
\begin{eq*}
\partial_i g(u^\gamma) \left( U_\ad^i - u_i^\gamma \right) \subseteq \lambda \partial_i g(u^\gamma) (U_\ad^i - u_i^\gamma).
\end{eq*}
Hence, we deduce with $\lambda = N$ the inclusion
\begin{eq*}
\eps \B_X \subseteq \partial_i g(u^\gamma) (U_\ad^i - u_i^\gamma) - (K - g(u^\gamma)) \subseteq N \partial_i g(u^\gamma) (U_\ad^i - u_i^\gamma) - (K - g(u^\gamma))
\end{eq*}
for all $i = 1,\dots,N$. Taking an arbitrary $x \in \B_X$ yields hence the existence of $v_i^\gamma \in U_\ad^i$ and $k_i^\gamma \in K$ such that
\begin{eq*}
\eps x = N \partial_i g(u^\gamma) (v_i^\gamma - u_i^\gamma) - (k_i^\gamma - g(u^\gamma)) \text{ for all } i = 1,\dots,N.
\end{eq*}
Summing over $i=1, \dots, N$ and dividing by $N$ yields
\begin{eq*}
\eps x &= \sum_{i = 1}^N \partial_i g(u^\gamma)(v_i^\gamma - u_i^\gamma) - \left( \frac{1}{N}\sum_{i = 1}^N k_i^{\msa{\gamma}} - g(u^\gamma) \right)\\
&= Dg(u^\gamma) (v - u^\gamma) - (k^{\msa{\gamma}} - g(u^\gamma))
\end{eq*}
with $k^{\msa{\gamma}} := \frac{1}{N}\textstyle\sum_{i = 1}^N k_i^{\msa{\gamma}} \in K$. Hence, we obtain the uniform Robinson-type condition.
\end{rem}
%
\msa{As}{} the non-linearity of the operator $g$ and its derivative play a role \msa{in \refer{sys:gnep:penalized:firstorder} and \refer{sys:vep:penalized:firstorder}}{} \msa{for the behaviour of the multipliers}{} we need to discuss the influence of weak-to-weak as well as complete continuity on the first derivative. Therefore, we introduce the concept of uniform Fr\'echet differentiability.
%
%
\begin{defn}\label{defn:uniformfrechet}
(see \cite[Definition 5.1.2]{bib:Llavona}.) Let $X, Y$ be Banach spaces. A Fr\'echet differentiable operator $T: X \rightarrow Y$ is called \emph{uniformly Fr\'echet differentiable} on a set $M \subseteq X$, if 
\begin{eq*}
\lim_{\|h\|_X \rightarrow 0}\sup_{x \in M} \frac{\|T(x + h) - T(x) - DT(x)h\|_Y}{\|h\|_X} = 0
\end{eq*}
holds.
\end{defn}
As preparation for the upcoming convergence result in \refer{thm:gnep_vep:path:conv} we provide the following Lemma\ms{ta containing properties of uniformly Fr\'echet differentiable. The proofs are given in the appendix.}
\begin{lem}\label{lem:uniformfrechet:properties}
Let $X,Y$ be Banach spaces and let furthermore $X$ be reflexive.
\begin{enum}
\item\label{enum:state:uniformfrechet:completelycont} Let $T: X \rightarrow Y$ be completely continuous and uniformly Fr\'echet differentiable on every bounded set.
Then, for weakly convergent sequences $x_n \rightharpoonup x$, $h_n \rightharpoonup h$ in $X$ and \msa{weakly$^*$-convergent sequence}{} $y_n \rightharpoonup^* y^*$ in $Y^*$ \mh{it}{} holds \mh{that}{}
\begin{eq*}
DT(x_n)h_n \rightarrow DT(x)h \text{ in } Y \text{ and } DT(x_n)^* y_n^* \rightarrow DT(x)^* y^* \text{ in } X^*.
\end{eq*}
\item\label{enum:state:uniformfrechet:weaklycont} Let $T: X \rightarrow Y$ be weakly continuous and uniformly Fr\'echet differentiable on every bounded set.
Then, for weakly convergent sequences $x_n \rightharpoonup x$, $h_n \rightharpoonup h$\linebreak[2] in $X$ and \msa{weakly$^*$-convergent sequence}{} $y_n \rightharpoonup^* y^*$ in $Y^*$ \mh{it}{} holds \mh{that}{} 
\begin{eq*}
DT(x_n)h_n \rightharpoonup DT(x)h \text{ in } Y \text{ and } DT(x_n)^* y_n^* \rightharpoonup^* DT(x)^* y^* \text{ in } X^*.
\end{eq*}
\end{enum}
\end{lem}

\mh{In view of the previous result let us point out}{} that in the case of uniform Fr\'echet differentiability over every bound\-ed set, continuity properties of the operator itself are inherited by its derivative as well as the pointwise dual derivative.
\begin{lem}\label{lem:bddderivative}
Let $X,Y$ be Banach spaces \mh{with}{} $X$ reflexive. Consider a weakly continuous and Fr\'echet differentiable operator $T: X \rightarrow Y$. Then, the first derivative $DT: X \rightarrow \calL(X,Y)$ is a bounded operator.
\end{lem}

Next, further properties of uniform Fr\'echet differentiable operators are discussed. Of interest is the discussion of properties  of \emph{inverse operators} as well as \emph{compositions}.
\begin{lem}\label{lem:uniformfrechet:inverse}
Let $T:X \rightarrow Y$ be uniformly Fr\'echet differentiable on every bounded subset as well as bijective\mh{,}{} and \mh{let}{} its inverse $S:Y \rightarrow X$ be continuously Fr\'echet differentiable. Moreover assume, that $S$ and the mapping $x \mapsto DT(x)^{-1}$ are bounded operators.\newline
Then, $S$ is uniformly Fr\'echet differentiable on every bounded subset of $Y$.
\end{lem}
The corresponding result for compositions is addressed in the following \mh{l}emma.
\begin{lem}\label{lem:uniformfrechet:composition}
Let $T_1:X \rightarrow Y$ and $T_2: Y \rightarrow Z$ be uniformly Fr\'echet differentiable on every bounded subset and let $T_1$ and $DT_1: X \rightarrow \calL(X,Y)$ as well as $DT_2: Y \rightarrow \calL(Y,Z)$ be bounded operators as well. Then, the composition $T_2\circ T_1$ is uniformly Fr\'echet differentiable on every bounded subset of $Y$.
\end{lem}
As an application of uniform Fr\'echet differentiability we discuss the relationship between the uniform $\mathrm{RZK}$-condition \refer{eq:gnep:path:uniformrobinson} and \refer{eq:gnep:rzk}.
\begin{prop}
Let a weakly convergent sequence $u^\gamma \rightharpoonup u$ in $U$ for $\gamma \rightarrow +\infty$ with $u^\gamma \in U_\ad$ be given and assume $g_i:U_\ad \to X_i$ to be weakly continuous and uniformly Fr\'echet differentiable on $U_\ad$. Moreover assume, that the uniform Robinson-type condition \refer{eq:gnep:path:uniformrobinson} is fulfilled. If $g_i(u) \in K_i$, then \refer{eq:gnep:rzk} is fulfilled.
\end{prop}
\begin{proof}
Since $U_\ad$ is a closed, convex subset it holds $u \in U_\ad$. Taking an arbitrary $\msa{x_i} \in X_i$ with $\|x_i\|_{X_i} \leq \varepsilon$. By \refer{eq:gnep:path:uniformrobinson} there exist for all $\gamma > 0$ elements $v_i^\gamma \in U^i_\ad$ and $k^\gamma_i \in K_i$ with
\begin{eq*}
\msa{x_i} = \partial_i g_i(u^\gamma)(v_i^\gamma - u_i^\gamma) - (k_i^\gamma - g_i(u^\gamma)).
\end{eq*}
Since $(v_i^\gamma)_{\gamma > 0} \subseteq U_\ad^i$ is bounded and $U$ is reflexive, there exists a subsequence (not relabelled) with $v_i^\gamma \rightharpoonup v_i$ in $U_i$ \msa{and}{} $v \in U_\ad^i$. Then, we obtain
\begin{eq*}
k_i^\gamma = g_i(u^\gamma) - \msa{x_i} - \partial_i g_i(u^\gamma) (v_i^\gamma - u_i^\gamma) \rightarrow g_i(u) - \msa{x_i} - \partial_i g_i(u)(v_i - u_i) =: k_i \in K_i
\end{eq*}
by the use of \refer{lem:uniformfrechet:properties}\ref{enum:state:uniformfrechet:weaklycont} as well as the convexity and closedness of $K_i$. Hence, we have found $v \in U_\ad^i$ and $k_i \in K_i$ with 
\begin{eq*}
\msa{x_i} = \partial_i g_i(u)(v_i-u_i) - (k_i-g_i(u)),
\end{eq*}
\msa{which implies}{} $0 \in \interior{\partial_i g_i(u)(v_i-u_i) - (K_i-g_i(u))}$. Since by assumption $g_i(u)\in K_i$ \msa{holds true}{} we conclude by \cite[see Equations (3.2)--(3.4)]{bib:ZoweKurcyusz} the relation \refer{eq:gnep:rzk}.
\end{proof}
\mh{W}e are now ready \msa{to}{} \mh{return to Nash games and}{} prove \mh{a}{} convergence theorem for the first order condition.
\begin{thm}\label{thm:gnep_vep:path:conv}
Let the conditions of \refer{lem:gnep_vep:path:bdd} be fulfilled. Moreover, assume:
\begin{enum}
\item\label{enum:cond:gnep_vep:path:conv:objectivederivative} \mh{For $i = 1, \dots, N,$}{} \mh{t}he operator $u \mapsto \partial_i \calJ_i(u)$ has a weakly-weakly* closed graph, (i.e. for all sequences with $u^n \rightharpoonup u$ and $q_i^n := \partial_i \cJnu(u^n) \rightharpoonup^* q_i $ \mh{it}{} holds \mh{that}{} $q_i = \partial_i \cJnu(u)$). Let furthermore
\begin{eq*}
\langle \partial_i \cJnu(u), u_i \rangle_{\msa{U_i^*,}U_i} \leq \limsup_{n \rightarrow \infty}\, \langle \partial_i \cJnu(u^n), u_i^n \rangle_{\msa{U_i^*,}U_i}
\end{eq*}
be fulfilled.
\item\label{enum:cond:gnep_vep:path:conv:weakcont} The mappings $g_i: U \rightarrow X_i$ are weakly continuous \mh{for}{} all $i = 1,\dots,N\mh{,}$ and for all sequences $u^n \rightharpoonup u$ in $U$ and $\mu_i^n \rightharpoonup^* \mu_i$ in $X_i^*$ \mh{it}{} holds \mh{that}{}
\begin{eq*}
\partial_i g_i(u^n)^* \mu_i^n \rightarrow \partial_i g_i(u)^* \mu_i \text{ in } U_i^*
\end{eq*}
(respectively $g:U \rightarrow X$ is weakly continuous and for all sequences $u^n \rightharpoonup u$ in $U$ and $\mu^n \rightharpoonup^* \mu$ in $X^*$ \mh{it}{} holds \mh{that}{} $Dg(u^n)^*\mu^n \rightarrow Dg(u)^* \mu$ in $U^*$).
\end{enum}
Then, every path has a limit point $(u,q,\lambda,\mu)$ and every such limiting point fulfils the first order system \eqref{sys:gnep:firstorder} in \refer{thm:gnep:firstorder} \msa{for a}{} \gnep{} (resp. the first order system \eqref{sys:vep:firstorder} in \refer{thm:vep:firstorder} \msa{for a}{} \vep{}).
\end{thm}
\begin{proof}
Since the path \mh{induced by $\gamma_n \nearrow +\infty$}{} is bounded, there exist weakly convergent subsequences with
\begin{eq*}
u_i^{\gamma_n} \rightharpoonup u_i \text{ in } U_i,\quad q_i^{\gamma_n} \rightharpoonup q_i \text{ in } U_i^*,\quad \mu_i^{\gamma_n} \rightharpoonup^* \mu_i \text{ in } X_i^*,\quad \lambda_i^{\gamma_n}\rightharpoonup \lambda_i \text{ in } U_i^*
\end{eq*}
due to the reflexivity of $U_i$ and the Banach-Alaoglu theorem for \msa{$X_i^*$}.
Due to the weak-weak* closedness of the graph of $\partial_i \calJ_i$ we conclude by $(u^{\gamma_n},q_i^{\gamma_n}) \in \gph{\partial_i \calJ_i}$ that $q_i \in \partial_i \calJ_i(u)$.
Applying the assumption of this theorem we deduce 
\begin{eq*}
\partial_i g_i(u^{\gamma_n})^* \mu_i^{\gamma_n} \rightarrow \partial_i g_i(u)^* \mu_i \text{ in } U_i^*.
\end{eq*}
For arbitrary $v_i \in U_\ad^i$ we obtain
\begin{eq*}
\langle -\lambda_i, u_i \rangle_{\msa{U_i^*,}U_i} &= \langle q_i - \partial_i g_i(u)^* \mu_i, \bar u_i \rangle_{\msa{U_i^*,}U_i} \leq \limsup_{n \rightarrow \infty}\, \langle q_i^{\gamma_n} - \partial_i g_i(u^{\gamma_n})^* \mu_i^{\gamma_n}, u_i^{\gamma_n}\rangle_{\msa{U_i^*,U_i}} \\
&= \limsup_{n \rightarrow \infty}\, \langle -\lambda_i^{\gamma_n}, u_i^{\gamma_n} \rangle_{\msa{U_i^*,U_i}} \leq \limsup_{n \rightarrow \infty}\, \langle -\lambda_i^{\gamma_n}, v_i \rangle_{\msa{U_i^*,U_i}} \\
&= \langle q_i - \partial_i g_i(u)^* \mu_i, v_i \rangle_{\msa{U_i^*,U_i}} = \langle -\lambda_i , v_i\rangle_{\msa{U_i^*,U_i}}\mh{,}
\end{eq*}
and hence $\lambda_i \in N_{U_\ad^i}(u_i)$. Finally\msa{,}{} we show $\mu_i \in N_{K_i}(g_i(u))$. \mh{For this purpose}{} take an arbitrary $k_i \in K_i$. By the convexity of $\beta_i$ \mh{it}{} holds \mh{that}{}
\begin{eq*}
0 &\leq \gamma \beta_i(g_i(u^{\gamma_n})) \leq \gamma_n \beta_i(k_i) - \gamma_n \langle D\beta_i(g_i(u^{\gamma_n})), k_i - g_i(u^{\gamma_n}) \rangle_{\msa{X_i^*,X_i}} \\
&= - \langle \mu_i^{\gamma_n}, k_i - g_i(u^{\gamma_n}) \rangle_{\msa{X_i^*,X_i}}.
\end{eq*}
Due to the boundedness of $\mu_i^{\gamma_n}$ and $g_i(u^{\gamma_n})$ also $\gamma_n \beta_i(g_i(u^{\gamma_n}))$ is bounded\mh{,}{} and therefore $\beta_i(g_i(u^{\gamma_n})) \rightarrow 0$. Since $g_i$ is weakly continuous, it holds \mh{that}{} $0 \leq \beta_i(g_i(u)) \leq \liminf_{n \rightarrow \infty} \beta_i(g_i(u^n)) = 0$ and therefore $g_i(u) \in K_i$. Since $k_i \in K_i$ is a minimizer of $\beta_i$, \mh{one has}{} $D \beta_i(k_i) = 0$. Using the monotonicity of $D\beta_i$ \mh{we}{} observe \mh{that}{}
\begin{eq*}
\langle \mu_i, \msa{k_i}{} - g_i(u) \rangle_{\msa{X_i^*,X_i}} &= \lim_{n \rightarrow \infty} \langle \mu_i^{\gamma_n}, k_i - g_i(u^{\gamma_n}) \rangle_{\msa{X_i^*,X_i}}\\
&\geq \lim_{n \rightarrow \infty} \langle -\gamma_n D \beta_i(k_i), k_i - g_i(u^{\gamma_n}) \rangle_{\msa{X_i^*,X_i}} = 0.
\end{eq*}
This implies eventually $\mu_i \in K_i^+$ and $\langle \mu_i, g_i(u) \rangle_{\msa{X_i^*,X_i}} = 0$.
\end{proof}

\def\bg{\bar\gamma}

With these results at hand, we laid the theoretical foundation for the path-following strategy. The next step towards an applicable numerical method is the careful choice of an update strategy for the penalty parameter.
\msa{\subsection{Path-Following Strategy}}
\msa{So far,}{} we \mh{approximated}{} \msa{\gnep{}s}{} respectively \vep{}s \mh{by}{} a sequence of \msa{equilibrium}{} problems without moving sets. \msa{In}{} practice however, the careful selection of the penalty parameter \mh{$\gamma$}{} is of crucial importance. This has
several reasons: \mh{On}{} the one hand\mh{,}{} all results are inherently asymptotic and hence one is interested in a \mh{rapidly}{} increasing sequence of penalty parameters. On the other hand, every \mh{selection}{} of $\gamma$ \mh{leads to}{} the solution of an equilibrium
problem which---depending on the particular form of the problem---might be a
challenging task \mh{in}{} its own right \msa{and is  addressed by iterative methods.}{} The latter decisively depend on a careful selection of starting points. \mh{Often, a good}{} choice is \mh{given by}{} the solution of the previous iterate \mh{associated with \msa{its respective}{} $\gamma$-value. Naturally,}{} one is \mh{then}{} interested that the new
solution is not too far away form the previous one \mh{in order}{} to \mh{enable}{} fast convergence.
The \mh{main}{} aim of \mh{our}{} \emph{path-following strategy} is to balance these \mh{often conflicting}{} interests by
using a \mh{value}{} function as an indication of a rapid change in the solution and
hence a \mh{more moderate}{} update of $\gamma$.\\
First, we prove the following \mh{rather general l}emma on the differen\ms{tiability}{} of a certain type of merit function. For more on this compare to \cite[Section 4.3]{bib:BonnansShapiro}.

\begin{lem}\label{lem:merit}
Let $U$ be a reflexive Banach space and $f: U \rightarrow \overline{\R}$ be a proper functional. Consider a penalty functional $\pi: U \rightarrow [0,+\infty)$ and introduce the sequence of \mh{merit}{} functionals
\begin{eq*}
f_\gamma(u) := f(u) + \gamma \pi(u)\mh{,}{}
\end{eq*}
together with the \mh{(}optimal\mh{)}{} value function
\begin{eq*}
W:(0,+\infty) \rightarrow \R,\quad W(\gamma) := \inf_{u \in U} f_\gamma(u).
\end{eq*}
Then, the following assertions hold true:
\begin{enum}
\item\label{enum:state:merit:monotoneconcave} The functional $W :(0,+\infty) \rightarrow \R$ is non-decreasing and concave.
\item\label{enum:state:merit:diffestimate} Let $f, \pi$ be \msa{lower semi-continuous}{} and the set of minimizers $\argmin f_\gamma$ be non-empty for all $\gamma > 0$. Additionally define for $\gamma > 0$ the set 
\begin{eq*}
R(\gamma) := \left\{ u_\gamma \in \argmin f_\gamma : \forall \eta \rightarrow 0 \ \exists u_{\gamma + \eta} \in \argmin f_{\gamma + \eta} \text{ with } u_{\gamma + \eta} \rightharpoonup u_\gamma\right\}.
\end{eq*}
Then
\begin{eq*}
\sup_{u_\gamma \in R(\gamma)} \pi(u_\gamma) &\leq \liminf_{\eta \searrow 0} \frac{W(\gamma + \eta) - W(\gamma)}{\eta}\\
&\leq \limsup_{\eta \searrow 0} \frac{W(\gamma + \eta) - W(\gamma)}{\eta} \leq \inf_{u_\gamma \in \argmin f_\gamma} \pi(u_\gamma)
\end{eq*}
holds. If $\pi$ is continuous, then also
\begin{eq*}
\sup_{u_\gamma \in \argmin f_\gamma} \pi(u_\gamma) &\leq \liminf_{\eta \nearrow 0} \frac{W(\gamma + \eta) - W(\gamma)}{\eta}\\
&\leq \limsup_{\eta \nearrow 0} \frac{W(\gamma + \eta) - W(\gamma)}{\eta} \leq \inf_{u_\gamma \in R(\gamma)} \pi(u_\gamma)
\end{eq*}
holds.
\item\label{enum:state:merit:suffconddiff} If $f, \pi$ are \msa{lower semi-continuous}{} and the solution
mapping $\gamma \mapsto \argmin f_\gamma$ is singleton and continuous, then $W$ is 
differentiable with
\begin{eq*}
W'(\gamma) = \pi(u_\gamma).
\end{eq*}
\item\label{enum:state:merit:derivativewhenexists} If $\argmin f_\gamma$ is non-empty for all $\gamma > 0$ and the functional $W$ is differentiable in $\bar\gamma$, then $W'(\bar \gamma) = \pi(x_{\bar\gamma})$ holds for any $x_{\bar\gamma} \in \argmin f_{\bar\gamma}$.
\end{enum}
\end{lem}
\begin{proof}
(compare also to \cite[Proposition 4.12]{bib:BonnansShapiro})\\
ad \ref{enum:state:merit:monotoneconcave}: The functionals $\gamma \mapsto f(u) + \gamma \pi(u)$ are non-decreasing and concave (affine) for every fixed $u \in U$. Hence the functional $\gamma \mapsto W(\gamma)$ is non-decreasing and concave as it is the pointwise infimum.\\
ad \ref{enum:state:merit:diffestimate}: Choose $\gamma,\eta > 0$ as well as $u_\gamma \in \argmin f_\gamma$ arbitrarily. By definition holds
\begin{align*}
W(\gamma + \eta) - W(\gamma) \leq f_{\gamma + \eta}(u_\gamma) - f_\gamma(u_\gamma) = \eta \pi(u_\gamma)
\end{align*}
and hence
\begin{align*}
\limsup_{\eta \searrow 0} \frac{W(\gamma + \eta) - W(\gamma)}{\eta} \leq \inf_{u_\gamma \in \argmin f_\gamma } \pi(u_\gamma).
\end{align*}
In the same fashion we get for $\eta < 0$, that
\begin{align*}
\liminf_{\eta \nearrow 0} \frac{W(\gamma + \eta) - W(\gamma)}{\eta} \geq \sup_{u_\gamma \in \argmin f_\gamma} \pi(u_\gamma).
\end{align*}
Let now $u_\gamma \in R(\gamma)$ and $\eta_k \searrow 0$ be arbitrary. There exists a sequence $(u_k)_{k \in \N}$ with $u_k \in \argmin f_{\gamma + \eta_k}$ and $u_k \rightharpoonup u_\gamma$. Then, we get
\begin{eq*}
W(\gamma + \eta_k) - W(\gamma) \geq f_{\gamma + \eta_k}(u_k) - f_\gamma (u_k) = \eta_k \pi(u_k).
\end{eq*}
By the \lscity{} we obtain
\begin{eq*}
\pi(u_\gamma) \leq \liminf_{k \rightarrow \infty} \pi(u_k) \leq \liminf_{k \rightarrow \infty} \frac{W(\gamma + \eta_k) - W(\gamma)}{\eta_k}.
\end{eq*}
Hence, we see
\begin{eq*}
\sup_{u_\gamma \in R(\gamma)} \pi(u_\gamma) \leq \liminf_{\eta \searrow 0} \frac{W(\gamma + \eta) - W(\gamma)}{\eta}.
\end{eq*}
Using the continuity of $\pi$ one obtains analogously
\begin{eq*}
\limsup_{\eta \nearrow 0} \frac{W(\gamma + \eta) - W(\gamma)}{\eta} \leq \inf_{u_\gamma \in R(\gamma)} \pi(u_\gamma).
\end{eq*}
ad \ref{enum:state:merit:suffconddiff}: If now the mapping $\gamma \mapsto \argmin f_\gamma$ is singleton and continuous, we obtain $R(\gamma) = \argmin f_\gamma$ and the assertion follows by the use of Statement \ref{enum:state:merit:diffestimate}.\\
ad \ref{enum:state:merit:derivativewhenexists}: Assuming the differentiability of $W$ in $\bar\gamma$ choose an arbitrary $u_{\bar\gamma} \in \argmin f_{\bar\gamma}$ and define the functional
\begin{eq*}
V : (0,+\infty) \rightarrow \R,\quad V(\gamma) := W(\bar\gamma) - W(\gamma) + (\gamma - \bar\gamma) \pi(u_{\bar\gamma}).
\end{eq*}
Then $V$ is differentiable in $\bar\gamma$ and attains a global minimum there:
\begin{align*}
V(\gamma) &= f_{\bar\gamma}(u_{\bar\gamma}) - W(\gamma) + (\gamma - \bar\gamma) \pi(u_{\bar \gamma})\\
&\geq f(u_{\bar\gamma}) + \bar\gamma \pi(u_{\bar\gamma}) - f(u_{\bar\gamma}) - \gamma\pi(u_{\bar\gamma}) + (\gamma - \bar\gamma)\pi(u_{\bar\gamma}) = 0 = V(\bar\gamma).
\end{align*}
\msa{Thus}{} $0 = V'(\bar\gamma) = \pi(u_{\bar\gamma}) - W'(\bar\gamma)$ \msa{holds true}, yielding the assertion.
\end{proof}
\begin{rem}\label{rem:merit:references}
The above theorem can directly be applied to cover the corresponding results \msa{in}{} \cite[Theorems 4.3 -- 4.5]{bib:HintRasch}, \cite[Theorem 2.1]{bib:AdamHintSur} and partly in \cite[Proposition 4.2]{bib:HintKunischPath} for the case of $\bar \lambda = 0$ therein.
\end{rem}

The update principle of $\gamma$ is based on the expected influence after a change
of this parameter. This is achieved via its directional
derivative. If the increase is \mh{large}, the $\gamma$ update is chosen more
conservatively, whereas with a small derivative a more aggressive increase is
\mh{admitted}. In contrast to the corresponding techniques for optimization problems (see
\cite{bib:HintItoKunisch} and the references in \refer{rem:merit:references}) the \mh{value}{} functions do not only depend explicitly on the parameter $\gamma$
itself but also on the calculated equilibrium. Therefore a sequence of \mh{value}{}
functions is considered \msa{---}{} one for each $\gamma$-update.
\begin{defn}\label{defn:gnep_vep:penalized:merit}\makebox{}
\begin{enum}
\item\label{enum:part:gnep:penalized:merit} Consider \pengnep{}. Let for $\bar\gamma > 0$ a strategy $u^{\bar\gamma} \in U_\ad$ be given.
Define for \pengnep{} \msa{its}{} \emph{\mh{value}{} function} $\calW(\ccdot,u^{\bar\gamma}): (0,+\infty) \rightarrow \R$ by
\begin{eq*}
\calW(\gamma,u^{\bar\gamma}):= \inf_{v \in U_\ad} \left(\sum_{i = 1}^N \left(\rule{0pt}{11pt} J_i(v_i, u_{-i}^{\bar\gamma}) + \gamma\beta(g(v_i,u_{-i}^{\bar\gamma})) \right) \right).
\end{eq*}
\item\label{enum:part:vep:penalized:merit} Consider \penvep{}. Let for $\bar\gamma > 0$ an equilibrium $u^{\bar \gamma} \in U_\ad$ be given.
Define for \penvep{} \msa{its}{} \emph{\mh{value}{} function} $\calW(\ccdot, u^{\bar \gamma}): (0,+\infty) \rightarrow \R$ by
\begin{eq*}
\widehat{\calW}(\gamma,u^{\bar\gamma}):= \inf_{v \in U_\ad} \left(\sum_{i = 1}^N J_i(v_i,u_{-i}^{\bar\gamma}) + \gamma\beta(g(v))\right).
\end{eq*}
\end{enum}
\end{defn}
\msa{Here, the value functions $\calW$ and $\widehat\calW$ depend on the penalty parameter and the associated approximation of the equilibrium. In fact, a rapid change in these value functionals indicates in first place a change in the respective solution of the encoded minimization problem. The latter is also called the \emph{best response}, see eg. \cite{bib:HintStenglEquiGamma} and the references therein. This however can as well be interpreted as an indication of a change for the current iterate $u^{\bar\gamma}$, as the condition to be an equilibrium translates into the iterate and the best-response need to coincide.}{}\\
Utilizing the abstract differentiability results in \refer{lem:merit}, we obtain for the setup introduced in \refer{defn:gnep_vep:penalized:merit} the following \mh{c}orollary.
\begin{cor}
Let $\bar\gamma > 0$ be fixed and $u^{\bar\gamma} \in U_\ad$ be chosen. Let the functionals $v_i \mapsto \cJnu(v_i,u_{-i})$ be weakly \lsc{}. Moreover, assume the operators $g_i:U\to X_i$ (respectively the operator $g:U\to X$) to be weakly continuous.\\
For the \mh{value}{} function in \refer{defn:gnep_vep:penalized:merit}\mh{\ref{enum:part:gnep:penalized:merit} there}{} holds
\begin{eq*}
\limsup_{\eta \searrow 0} \frac{\calW(\bar\gamma + \eta, u^{\bar\gamma}) - \calW(\bar\gamma,u^{\bar\gamma})}{\eta} \leq \sum_{i = 1}^N \beta_i(g_i(u^{\bar\gamma})).
\end{eq*}
If for every fixed $u \in U_\ad$ the best response mapping $\calB^\gamma$ defined by
\begin{eq*}
\calB_i^\gamma(u_{-i}) := \argmin_{v_i \in U_\ad^i} \left( \calJ_i(v_i, u_{-i}) + \gamma \beta_i(g_i(v_i, u_{-i})) \right)
\end{eq*}
and $\calB^\gamma(u) := \prod_{i = 1}^N \calB_i^\gamma(u_{-i})$ is \mh{a}{} singleton and continuous in $\gamma$, then $\calW$ is directionally differentiable with 
$\calW'(\bar\gamma; 1) = \sum_{i = 1}^N \beta_i(g_i(u^{\bar\gamma}))$.\\
For the \mh{value}{} function in \refer{defn:gnep_vep:penalized:merit}\mh{\ref{enum:part:vep:penalized:merit} there}{} holds
\begin{eq*}
\limsup_{\eta \searrow 0} \frac{\widehat{\calW}(\bar\gamma + \eta, u^{\bar\gamma}) - \widehat{\calW}(\bar\gamma, u^{\bar\gamma})}{\eta} \leq \beta(g(u^{\bar\gamma})).
\end{eq*}
If for every fixed $u \in U_\ad$ the response mapping $\widehat{\calB}^\gamma$ defined by
\begin{eq*}
\widehat{\calB}^\gamma(u) := \argmin_{v \in U_\ad} \left( \sum_{i = 1}^N \calJ_i(v_i, u_{-i}) + \gamma \beta(g(v)) \right)
\end{eq*}
is \mh{a}{} singleton and continuous in $\gamma$, then $\widehat{\calW}$ is differentiable with $\widehat{\calW}'(\gamma;1) = \beta(g(u^\gamma))$.
\end{cor}
\begin{proof}
We check the conditions of \refer{lem:merit}. \mh{For this purpose,}{} define for chosen $u^\gamma \in U_\ad$ the functionals
$f(u) := \sum_{i = 1}^N \cJnu(u_i,u^\gamma_{-i}) + I_{U_\ad}(u)$ and $\pi(u) = \sum_{i = 1}^N \beta_i(g_i(u_i,u_{-i}^\gamma))$.
As $\cJnu(\ccdot,u^\gamma_{-i})$ \mh{is}{} weakly lsc. \mh{for all $i \ms{= 1, \dots, N$}{}}{} and $U_\ad$ is closed, convex and hence weakly closed, $f$ is weakly lsc.
Due to the assumed weak continuity of $g_i$ and the weak \lscity{} of $\beta_i$ (see \refer{defn:gnep_vep:penalized}), also the functional $\pi$ is weakly \lsc{}.
By \refer{lem:merit}\ref{enum:state:merit:monotoneconcave} the functional $\calW(\ccdot,u^\gamma)$ is non-decreasing and concave.
Using the inequality for the right-sided difference quotient in \refer{lem:merit}\ref{enum:state:merit:diffestimate} proves the asserted inequality.\\
For the response mapping holds $\calB^\gamma(u^{\bar\gamma}) = \argmin f^\gamma$. Hence, we can apply \refer{lem:merit}\ref{enum:state:merit:suffconddiff} and obtain the directional differentiability in positive direction.
\end{proof}
As update strategy \ms{one approach is to establish a model function aiming at
resembling the behavior of the value function as successfully performed in
\cite{bib:HintKunischPath}, \cite{bib:HintKunischPathLowMultiplier},
\cite{bib:HintRasch} and \cite{bib:AdamHintSur}). Due to the presence of the
computed equilibrium we use the simplified strategy in \cite{bib:HintSurKaemmler}.
The intention is to perform a more conservative $\gamma$-update for larger changes in the solution and more aggressive ones for small changes in the solution. As a way to indicate these changes the value functions $\calW, \hat\calW$ respectively the estimates of their directional derivatives are used. Is the estimate large, a smaller $\gamma$-update is performed. Is the estimate small, a larger update is added. From this consideration}{} we utilize the technique in
\cite{bib:HintSurKaemmler} and recover the \ms{subsequent algorithms}{}:\newline
\begin{algorithm}[H]
	\SetAlgoLined
	\DontPrintSemicolon
	\caption{\textsc{Nikaido--Isoda-based path-following---GNEPs}}
	\KwData{Choose $\gamma_0 > 0$, $C_\mathrm{Path} > 0$, $\varepsilon > 0$.}
	\KwResult{Approximation of a Nash Equilibrium.}
	\For{$k = 1, 2, \dots $}{%
		Solve \pengnep{} for $\gamma = \gamma_k$ with result $u^{\msa{\gamma_k}}$.\;
		\eIf{$\beta_{\ms{i}}(g_i(u^{\gamma_k})) = 0$ for all $i = 1, \dots, N$}%
		{\Return $u^{\msa{\gamma_k}}$.\;}%
		{Set $\gamma_{k+1} := \gamma_k + \max\left( \frac{C_{\mathrm{Path}}}{\sum_{i = 1}^N \beta_i(g_i(u^{\gamma_k}))}, \varepsilon \right)$.\;}%
	}
\end{algorithm}
\vspace*{0.666em}\noindent%
\ms{Here, the speed of the method is controlled by the two parameters $\eps$ as a lower bound imposed on the $\gamma$-update and the parameter $C_{\mathrm{Path}}>0$ controlling the growth in the first order approximation
\begin{eq*}
\calW(\gamma + \eta, u^\gamma) \approx \calW(\gamma,u^\gamma) + \eta \sum_{i = 1}^N \beta_i(g_i(u^\gamma)).
\end{eq*}
With the condition $\calW(\gamma + \eta, u^\gamma) - \calW(\gamma,u^\gamma) \approx C_{\mathrm{Path}}$ this leads to the estimate $\eta \approx \frac{C_{\mathrm{Path}}}{\sum_{i = 1}^N \beta_i(g_i(u^\gamma))}$.}{} In the same fashion variational equilibria can be treated as well.\\[0.666em]
\begin{algorithm}[H]
	\SetAlgoLined
	\DontPrintSemicolon
	\caption{\textsc{Nikaido--Isoda-based path-following---VEPs}}
	\KwData{Choose $\gamma_0 > 0$, $C_\mathrm{Path} > 0$, $\varepsilon > 0$.}
	\KwResult{Approximation of a Variational Equilibrium.}
	\For{$k = 1, 2, \dots $}{%
		Solve \penvep{} for $\gamma = \gamma_k$ with result $u^{\msa{\gamma_k}}$.\;
		\eIf{$\beta(g(u^{\gamma_k})) = 0$}%
		{\Return $u^k$.\;}%
		{Set $\gamma_{k+1} := \gamma_k + \max\left( \frac{C_{\mathrm{Path}}}{\beta(g(u^{\msa{\gamma_k}}))}, \varepsilon \right)$.\;}%
	}
\end{algorithm}
With these results at hand we established a penalization technique suitable for the treatment of GNEPs and VEPs, \msa{investigated}{} the interconnections between these two concepts, \msa{developed}{} the convergence analysis as well as a practically usable update strategy for the penalty parameter. \msa{Returning to the}{} treatment of equilibrium problems involving PDE constraints, we transfer our results to a specialized framework \msa{which is tailored}{} to that particular \mh{setting}.
\section{Specialized Framework}\label{sec:specialized}
In proximity to tracking-type objectives in optimal control we focus our attention on the following class of Nash equilibrium problems, similarly to the one discussed in \cite{bib:HintSurKaemmler}:
\begin{eq}\label{eq:gnep:specialized}
&\minimize &&J_i^1(y) + J_i^2(u_i) \text{ over } (u_i, y) \in U_i\times Y\\
&\subjectto &&u_i \in U_\ad^i,~\calG(y) \in K,\\
&&&A(y) = f + B (u_i,u_{-i}) \text{ in } W.
\end{eq}%
The role of the operator $g$ is then fulfilled by the composition of the solution operator and the operator $\calG$ \msa{realizing}\mh{, together with the set $K$,}{} the state constraint. Rewriting \refer{eq:gnep:specialized} using \mh{a}{} reduced formulation yields
\begin{eq*}
&\minimize &&\calJ_i(u_i, u_{-i}) := J_i^1(u_i) + J_i^2(S(u_i,u_{-i})) \text{ over } u_i \in U_i\\
&\subjectto &&u_i \in U_\ad^i \text{ and } g(u) \in K\mh{,}{}
\end{eq*}
\mh{where $S$ corresponds to the solution map of the operator equation in the constraints of \refer{eq:gnep:specialized}.}{}
The reduced form fits directly into the framework of jointly constraint GNEPs (see \refer{defn:gnep:jointly_constrained}) as proposed in \refer{ass:gnep_vep:strategy_map}\ref{enum:part:vep:strategy_map}. In this setting the operator $g \mh{= \calG \circ S}{}$, which is responsible for the state constraint, is effectively determined by the solution operator of the underlying PDE and hence only implicitely defined. To guarantee the properties required by the results in the previous section we propose the following group of standing assumptions, which are in proximity \mh{of}{} the ones \msa{given}{} in \cite[Assumption 2.1]{bib:HintSurKaemmler}:
\begin{ass}\label{ass:specializedframework}\makebox{}
\begin{enum}
\item The spaces $U_i$ are reflexive, separable Banach spaces, $Y$, $W$ are reflexive Banach spaces, $X$ is a Banach space.
\item The embedding $Y \hookrightarrow X$ is continuous.
\item The space $U := \prod_{i = 1}^N U_i$ is equipped with the \mh{product}{} topology.
\item The operator $A: Y \rightarrow W$ is bijective and continuously Fr\'echet differentiable.
\item The mapping $B \in \calL(U,W)$ is a bounded, linear operator and $B_i \in \calL(U_i, W)$ is defined by $B_i u_i := B(u_i, 0_{-i})$.
\item The constraint set $K \subseteq X$ is a non-empty, convex, closed cone.
\item The strategy $U_\ad^i \subset U_i$ is non-empty, bounded, closed, convex.
\item The solution operator $S := A^{-1}(f + B \,\cdot\,): U \rightarrow Y$ of the equation
\begin{eq*}
A(y) = f + Bu
\end{eq*}
is assumed to be completely continuous.
\item The mapping $\calG: Y \rightarrow X$ is continuously Fr\'echet differentiable and we define $g:U \rightarrow X$ by $g = \calG \circ S$.
\item The first derivative $DA(y) \in \calL(Y,W)$ is invertible for all $y \in Y$ and the map $y \mapsto DA(y)^{-1}$ is bounded.
\item The set $\calF := \{ u \in U_\ad : g(u) \in K\}$ is non-empty.
\item The functionals $J_i^1: Y \rightarrow \R$ are completely continuous and $J_i^2: U_i \rightarrow \R$ are continuous and weakly \lsc. Moreover, all of them are assumed to be bounded from below and continuously Fr\'echet differentiable.
\item The operator $DJ_i^2 : U_i \rightarrow U_i^*$ has a weak-weak* closed graph.
\end{enum}
\end{ass}
Following closely the steps of the previous section, we discuss the path-following scheme in the specialized framework. Hence, we proceed with the derivation of the first order system.
\begin{cor}\label{cor:gnep_vep:firstorder:special}
Let $u \in U_\ad$ be a Nash equilibrium together with $y = S(u)$ and let the following (RZK) constraint qualification 
\begin{eq}[\label{eq:gnep:specialized:rzk}]
\left( D\calG(y) DA(y)^{-1} B_i \right) U_\ad^i(u_i) - K(\calG(y)) = X
\end{eq}
hold. Then there exist $p_i \in W^*$, $\mu_i \in X^*$ and $\lambda_i \in U_i^*$ for $i = 1,\dots,N$ fulfilling the following necessary first order system
\begin{sys}[\label{sys:gnep:firstorder:special}]
0 &= DJ_i^2(u_i) + B_i^* p_i + \lambda_i &&\hspace*{-14ex}\text{ in } U_i^*,\\
A(y) &= f + B u &&\hspace*{-14ex}\text{ in } W,\\
DA(y)^* p_i &= DJ_i^1(y) - D\calG(y)^* \mu_i &&\hspace*{-14ex}\text{ in } Y^*,\\
\lambda_i &\in N_{U_\ad^i}(u_i) &&\hspace*{-14ex}\text{ in } U_i^*,\\
K^+ \ni \mu_i &\perp \calG(y) \in K.
\end{sys}
Let $u \in U_\ad$ be a variational equilibrium together with $y = S(u)$ and let the following (RZK) constraint qualification
\begin{eq}[\label{eq:vep:specialized:rzk}]
\left( D\calG(y) DA(y)^{-1} B \right) U_\ad(u) - K(\calG(y)) = X
\end{eq}
hold. Then there exist $p_i \in W^*, \lambda_i \in U_i^*$ for $i = 1,\dots,N$ and $\mu \in X^*$ fulfilling the following necessary first order system
\begin{sys}[\label{sys:vep:firstorder:special}]
0 &= DJ_i^2(u_i) + B_i^* p_i + \lambda_i &&\hspace*{-14ex} \text{ in } U_i^*,\\
A(y) &= f + B u &&\hspace*{-14ex}\text{ in } W,\\
DA(y)^* p_i &= DJ_i^1(y) - D\calG(y)^* \mu &&\hspace*{-14ex}\text{ in } Y^*,\\
\lambda_i &\in N_{U_\ad^i}(u_i) &&\hspace*{-14ex}\text{ in } U_i^*,\\
K^+ \ni \mu &\perp \calG(y) \in K.\hspace*{7ex}
\end{sys}
\end{cor}
\begin{proof}
Let $u \in U_\ad$ be a Nash equilibrium. Using \refer{ass:specializedframework} one calculates the following derivatives
\begin{eq*}
\partial_i S(u) &= \partial_i \left( A^{-1}(f + B \ccdot) \right)(u) = (DA^{-1}(f + Bu)) B_i = DA(y)^{-1} B_i,\\
\partial_i g(u) &= \partial_i(\calG \circ S)(u) = D\calG(y) \partial_i S(u) = D\calG(y) DA(y)^{-1} B_i,\\
\partial_i \cJnu(u) &= \partial_i \left( J_i^1(S(\ccdot,u_{-i})) + J_i^2(\ccdot) \right)(u_i) = \partial_i S(u)^* DJ_i^1(y) + D J_i^2(u_i)\\
&= B_i^* DA(y)^{-*}DJ_i^1(y) + DJ_i^2(u_i).
\end{eq*}
\msa{Then}{} we obtain for \eqref{eq:gnep:rzk} that
\begin{eq*}
X = \partial_i g(u)U_\ad^i(u_{\msa{i}}) - K(g(u)) = \left(D\calG(y) DA(y)^{-1} B_i\right) U_\ad^i(u_i) - K(\calG(y))
\end{eq*}
and for \eqref{eq:vep:rzk} that
\begin{eq*}
X = Dg(u)U_\ad(u) - K(g(u)) = \left(D\calG(y) DA(y)^{-1} B\right) U_\ad(u) - K(\calG(y)).
\end{eq*}
Applying \refer{thm:gnep:firstorder} and \refer{thm:vep:firstorder} yields the existence of $\mu_i \in \msa{K^+}$ and $\lambda_i \in N_{U_\ad^i}(u_i)$ with
\begin{eq*}
0 = \langle \mu_i, g(u) \rangle_{\msa{X^*, X}} = \langle \mu_i, \calG(y) \rangle_{\msa{X^*, X}}
\end{eq*}
and
\begin{eq*}
0 &= \partial_i \cJnu(u) - \partial_i g(u)^* \mu_i + \lambda_i\\
&= B_i^* DA(y)^{-*} DJ_i^1(y) + D J_i^2(u_i) - B_i^* DA(y)^{-*} D\calG(y)^*\mu_i\\
&= DJ_i^2(u_i) + B_i^*DA(y)^{-*}\left( \msa{DJ_i^1(y)}{} - D\calG(y)^*\mu_i \right).
\end{eq*}
Setting $p_i = DA(y)^{-*}\left(\msa{DJ_i^1(y)} - D\calG(y)^* \mu_i \right)$ and using $y = S(u)$ yields the system
\begin{eq*}
0 &= DJ_i^2(u_i) + B_i^* p_i + \lambda_i,\\
A(y) &= f + Bu.\\
DA(y)^* p_i &= D J_i^{\msa{1}}(y) - D\calG(y)^* \mu_i
\end{eq*}
and thus the assertion. The proof for the variational equilibrium reads essentially the same and is hence omitted.
\end{proof}
\mh{N}ext, the penalty technique is applied. This leads to the following penalized NEP
\begin{eq}\label{eq:gnep:penalized:specialized}
&\minimize &&J_i^1(y') + J_i^2(v_i) + \gamma \beta(\calG(y)) \text{ over } v_i \in U_i,\\
&\subjectto &&\msa{v_i} \in U_\ad^i,\ A(y') = f + B (v_i,u_{-i}) \text{ in } W.
\end{eq}
Analogously one obtains for \penvep{} \mh{the}{} problem
\begin{eq}\label{eq:vep:penalized:specialized}
&\minimize \sum_{i = 1}^N \left( J_i^1(y'_i) + J_i^2(v_i) \right) + \gamma \beta(\calG(y')) \text{ over } \msa{v} \in U,\\
&\subjectto \hspace*{1.2ex} \hphantom{\sum_{i = 1}^N} \msa{v} \in U_\ad,\ A(y'_i) = f + B (v_i,u_{-i}) \text{ in } W \text{ for } i = 1, \dots, N,\\
&\hspace*{20.65ex}\text{and } A(y') = f + Bu' \text{ in } W.
\end{eq}
It is worth noting that the latter equilibrium problem requires another solution of the PDE involving \mh{feasible \ms{test}{} strategies $v_i$}{} only.\\
The first order system for both problems is formulated in the following \mh{c}orollary.
\begin{cor}\label{cor:gnep_vep:penalized:firstorder:special}
Let $u^\gamma \in U_\ad$ be a solution of \refer{eq:gnep:penalized:specialized} or \refer{eq:vep:penalized:specialized} together with states $y^\gamma = S(u^\gamma)$. Then there exist $p_i^\gamma \in W^*$, $\lambda_i^\gamma \in U_i^*$ for $i = 1, \dots N$ and $\mu^\gamma \in X^*$ fulfilling the following \msa{necessary}{} first order system
\begin{sys}[\label{sys:gnep_vep:penalized:firstorder:special}]
\hspace*{7ex}0 &= D J_i^2(u_i^\gamma) + B_i^* p_i^\gamma + \lambda_i^\gamma &&\hspace*{-14ex} \text{ in } U_i^*,\hspace*{7ex}\label{eq:gnep_vep:penalized:firstorder:specialized:opt}\\
A(y^\gamma) &= f + B u^\gamma &&\hspace*{-14ex}\text{ in } W,\\
DA(y^\gamma)^* p_i^\gamma &= DJ_i^1(y^\gamma) - D\calG(y^\gamma)^* \mu^\gamma &&\hspace*{-14ex}\text{ in }Y^*,\\
\lambda_i^\gamma &\in N_{U_\ad^i}(u_i^\gamma) &&\hspace*{-14ex}\text{ in } U_i^*,\label{eq:gnep_vep:penalized:firstorder:specialized:normalcone}\\
\mu &= -\gamma D\beta(\calG(y^\gamma)) &&\hspace*{-14ex}\text{ in } X^*.
\end{sys}
\end{cor}
\begin{proof}
Using the calculations \msa{in}{} the \mh{p}roof\mh{s}{} of \refer{cor:gnep_vep:firstorder:special}, \refer{thm:gnep:penalized:firstorder} and \refer{thm:vep:penalized:firstorder} one obtains \refer{sys:gnep_vep:penalized:firstorder:special}.
\end{proof}
\section{Examples}\label{sec:example}

\begin{figure}[H]
\includestandalone{domain}
\includegraphics[trim = {850 200 500 200}, clip, width=0.45\textwidth]{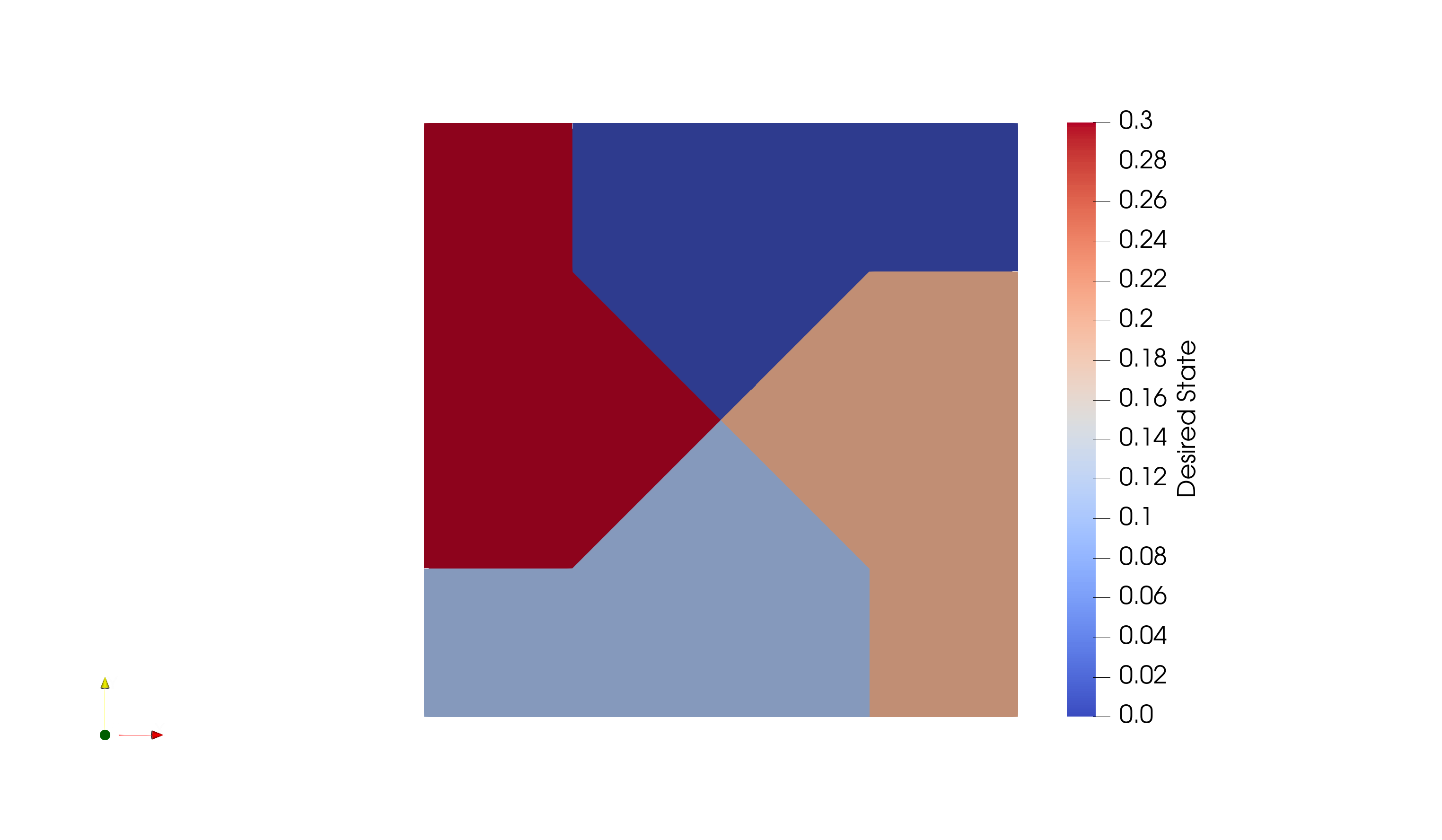}
\hspace*{0.05\textwidth}
\caption{Left: Decomposition of $\Omega = (0,1)^2$ into subdomains, Right: Plot of the desired state.}\label{fig:subdomains}
\end{figure}

\begin{figure}[H]
\includegraphics[trim = {850 200 500 200}, clip, width=0.45\textwidth]{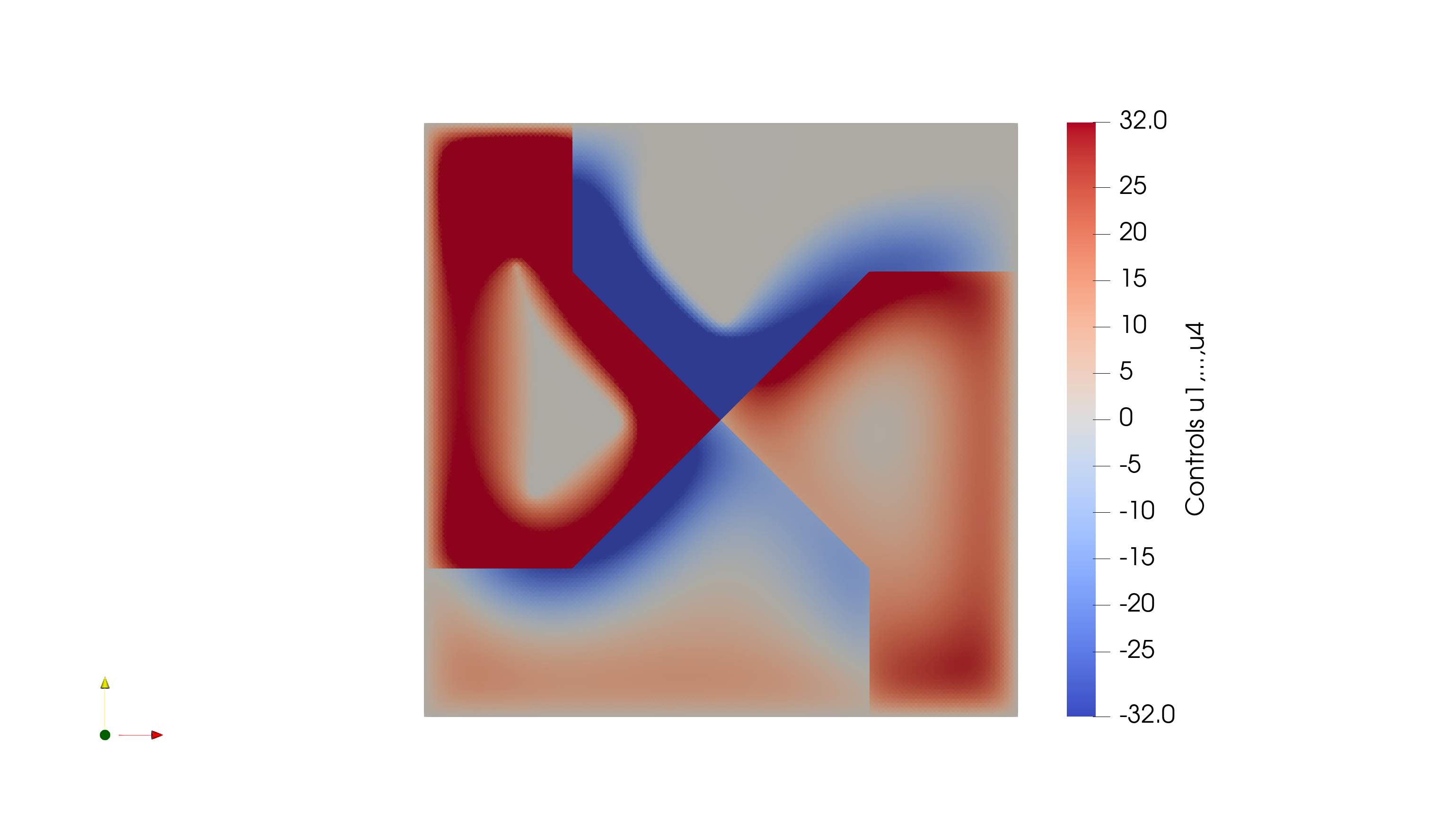}
\hspace*{0.05\textwidth}
\includegraphics[trim = {850 200 500 200}, clip, width=0.45\textwidth]{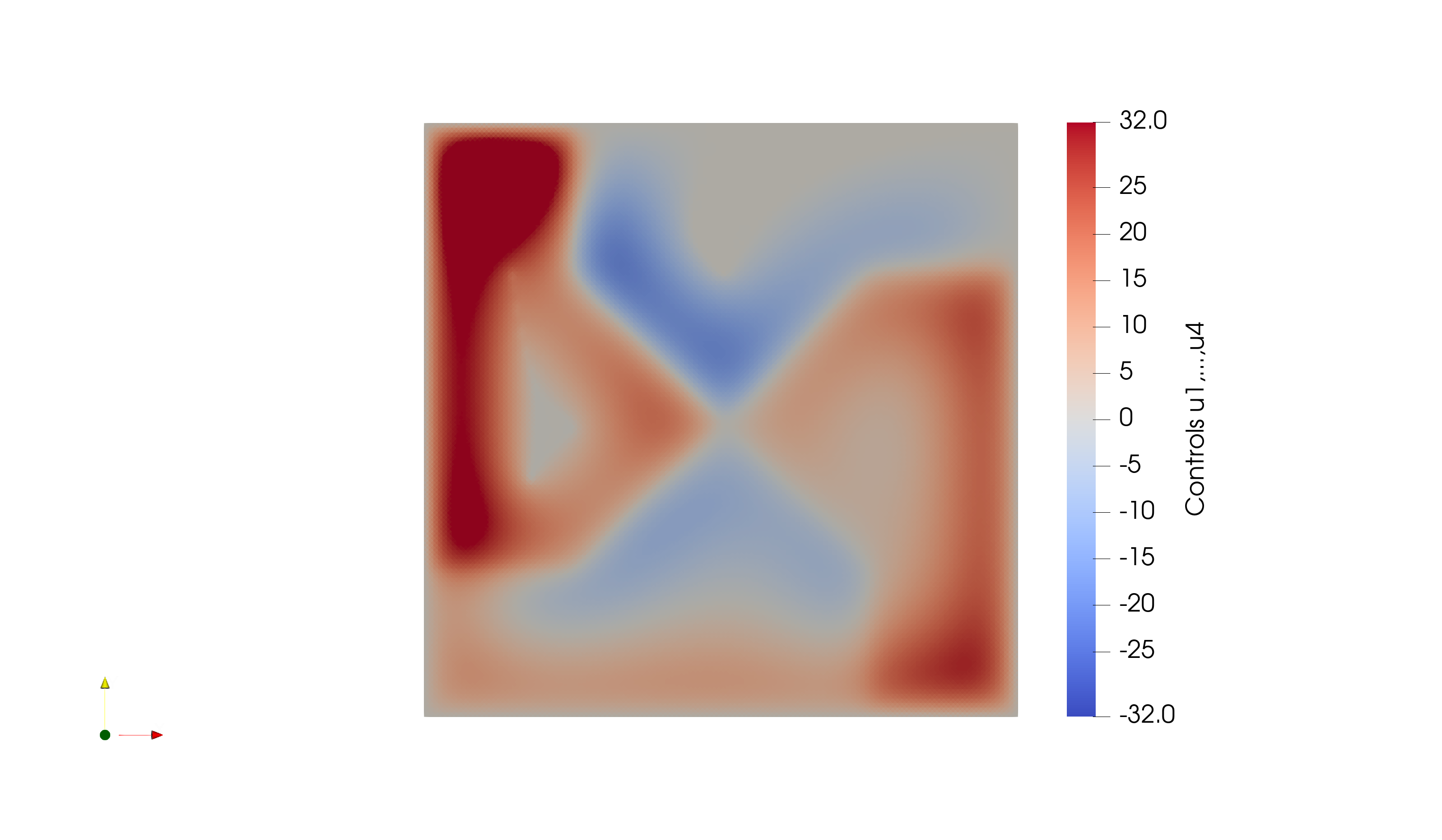}
\caption{Plot of the controls. Left: Results for Example \ref{ex:semielliptic:distr}, Right: Results for the cooperative version in \eqref{eq:ex:semielliptic:distr:coop}.}\label{fig:exp:elliptic_semi_distr:controls}
\end{figure}

\begin{figure}[H]
\includegraphics[trim = {850 200 500 200}, clip, width=0.45\textwidth]{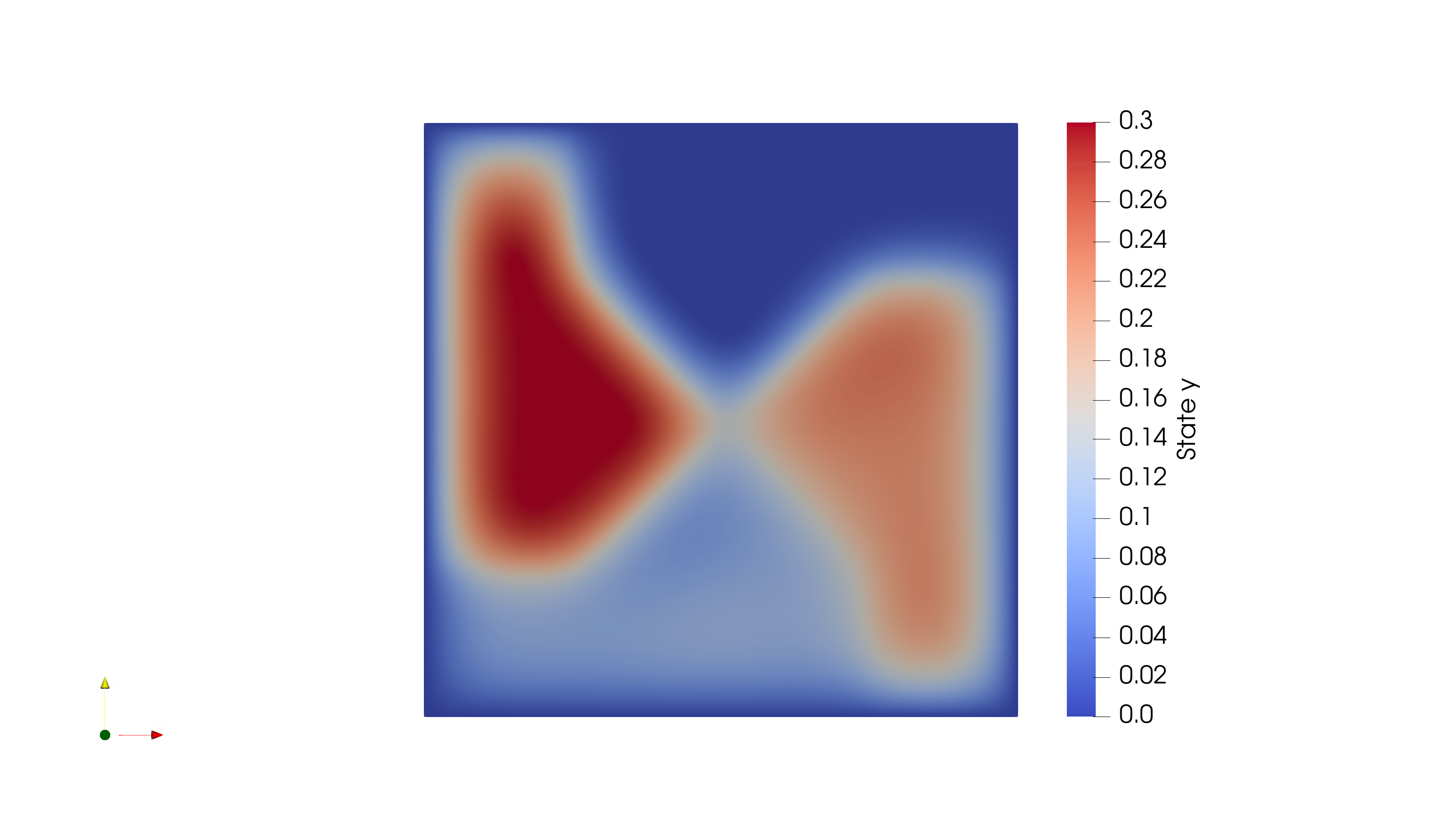}
\hspace*{0.05\textwidth}
\includegraphics[trim = {850 200 500 200}, clip, width=0.45\textwidth]{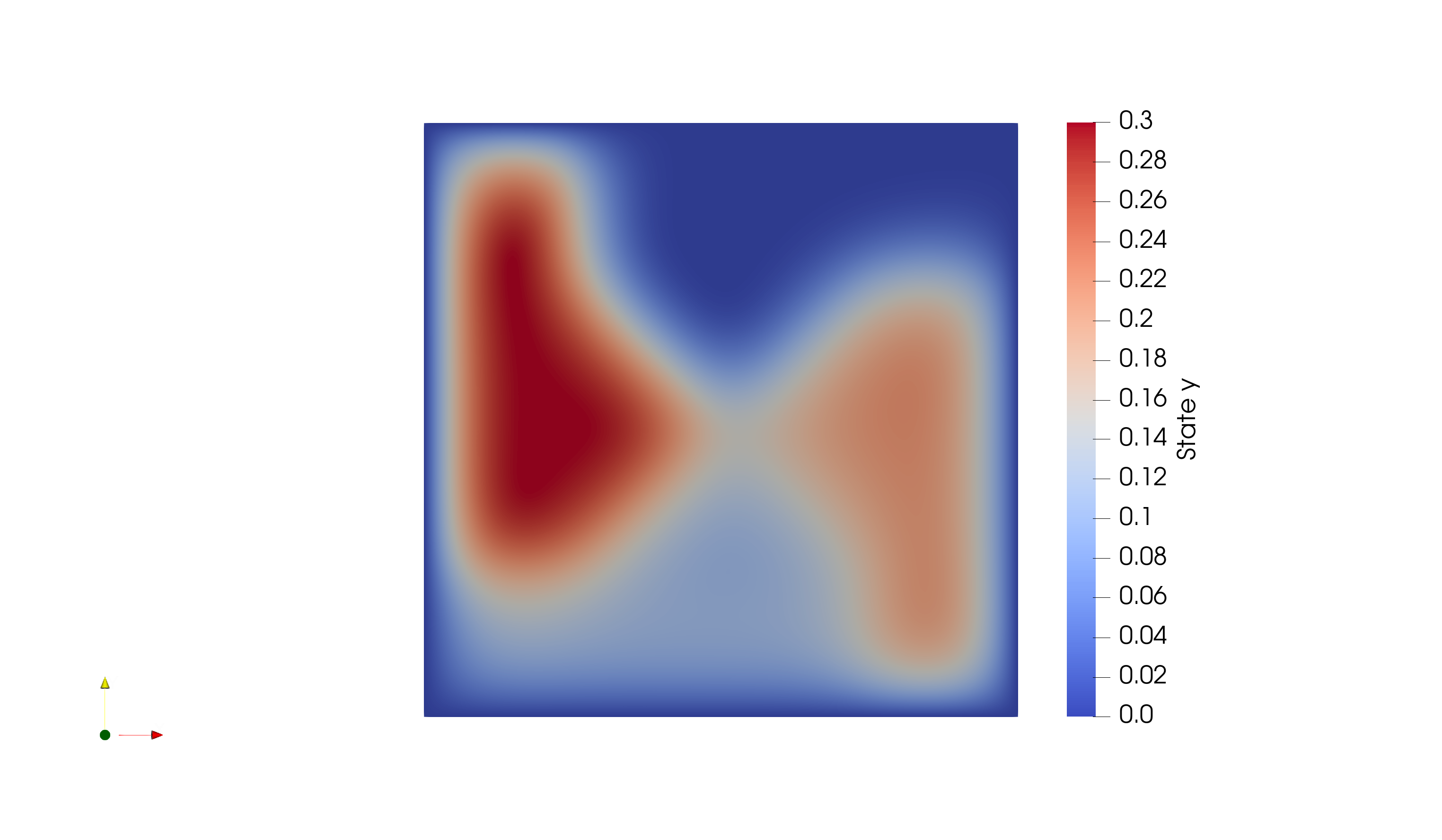}
\caption{Plot of the states. Left: Results for Example \ref{ex:semielliptic:distr}, Right: Results for the cooperative version \eqref{eq:ex:semielliptic:distr:coop}.}\label{fig:exp:elliptic_semi_distr:states}
\end{figure}

\begin{figure}[H]
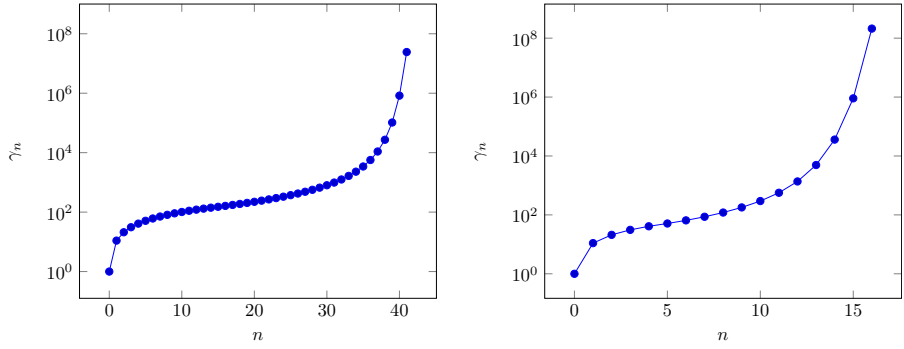

\begin{tabular}{cc}
\resizebox{0.47\textwidth}{!}{%
\includestandalone{nash_elliptic_semi_distr/gamma}} &
\resizebox{0.47\textwidth}{!}{%
\includestandalone{nash_elliptic_semi_distr_coop/gamma}}
\end{tabular}
\caption{Plot of the $\gamma$-updates. Left: Results for Example \ref{ex:semielliptic:distr}, Right: Results for the cooperative version in \eqref{eq:ex:semielliptic:distr:coop}.}\label{fig:exp:elliptic_semi_distr:gamma}
\end{figure}

In the remainder of this work we apply the \msa{developed}{} tools and methods to a selection of examples involving semi-linear elliptic PDEs.
\mh{In this context}{}, we establish two cases of Nash games with players \msa{having}{} tracking-type objectives. In the first example we will see an instance of distributed control and in the second one an example o\mh{f}{} boundary control.
In addition to the illustration \msa{of}{} the numerical methods developed within this work, this is taken as an opportunity to study the influence of the non-cooperative aspect of Nash games in comparison to \msa{an associated}{} cooperative optimization problem.\\
As a common setup we take as domain the unit square $\Omega = (0,1)^2 \subseteq \R^2$ and establish a partition into four subdomains as depicted in \refer{fig:subdomains}. Each of these domains serve\mh{s}{} as area of interest for exactly one player. The tracking type functional \msa{as well as}{} each player's control \msa{depend}{} on this players' domain only. However, their strategies have an effect on a state being the solution of semi-linear elliptic PDE. By this, a spatial coupling is established \msa{leading}{} to an interaction between these players. Additionally\ms{,}{} a box constraint for the state is established\ms{, which can be rewritten as a cone constraint in a product space}{}.\\
In the cooperative case however, the sum over all objectives is taken and minimized with respect to the combined strategies of all players simultaneously. As desired state we take a piecewise constant function with values $y_d = y_d^i$ a.e. on $\omega_i$ with constants $y_d^1 = 0.1$, $y_d^2 = 0.2$, $y_d^3 = 0$ and $y_d^4 = 0.3$. The desired state is plotted in \refer{fig:subdomains} as well. As constraints for the controls as well as the state we use box \msa{constraints}. 
For the solution of the corresponding first order systems we use in both cases a \emph{semi-smooth Newton method} (cf. \mh{\cite{bib:ChenNashedZuhair, bib:HintItoKunisch}}{}). Therefore, the semi-smoothness of the operator $\max(\ccdot, 0) : L^{r_1}(\Omega) \rightarrow L^{r_2}(\Omega)$ for $r_1 > r_2$ is utilized.
Since Newton-type methods only guarantee convergence, when the starting point is in a sufficiently small neighbourhood from a solution. This issue is addressed in two \msa{ways:}{}
On the one hand by the selection of parameters for the update strategy to avoid big changes in the solutions between the iterates. Then the solution with respect to the previous parameter can be used as a starting value.\\
On the other hand, one way to extend the basin of attraction is the use of a very simple damping technique. In our case, this is established by trying at first the iteration with\mh{out damping}{} and $25$ as the maximal number of iterat\mh{ions}.
Whenever the method fails to achieve both---the relative tolerance of $10^{-10}$ as well as the absolute tolerance of $10^{-10}$, the damping factor is halved, the maximal iteration number is doubled and the Newton method is attempted with the same
starting value again. In case of \msa{success,}{} the new settings are kept for all
subsequent $\gamma$-iterates. \mh{Clearly}{}, this procedure is heuristic \mh{only, but}{} \msa{it}{} performed well in our experiments. For more sophisticated damping strategies we just refer to \cite[Chapter 3]{bib:DeuflhardNewton} and for other techniques to \cite{bib:DirkseFerris}.\\
As discretization method we use in all examples \msa{the}{} finite element method (FEM) and for its implementation the package \fenics{} (cf. \cite{bib:fenics}). \msa{We discretize the domain using a}{} crossed mesh with $2^7$ segments along every side of the unit square. For the discretization of the state as well as each adjoint state we take the finite element \msa{spaces}{} 
\begin{eq*}
\calS^1(\calT) &:= \{z \in H^1(\Omega) : v|_{T} \in P_1(T) \text{ for all } T \in \calT\}
\end{eq*}
\msa{and}{}
\begin{eq*}
\calS_0^1(\calT) &:= \calS^1(\calT) \cap H_0^1(\Omega).
\end{eq*}
We are making use of the variational discretization technique as presented in \cite[Section 3.2.5]{bib:HinzePinnauMUlbrichSUlbrich}, see also the references therein. By that, the control remains undiscretized and the state equation is discretized by a Galerkin approach with the above introduced spaces. The derivation of the first order system leads to the corresponding discretization of the adjoint equation. The elimination of the control via the adjoint state leads to a non-linear coupled system on the discrete level and is hence numerically accessible.
\begin{ex}\label{ex:semielliptic:distr}
\ms{First, we return to the introductory \refer{ex:intro}, which is restated for convenience, and c}onsider the following \mh{g}eneralized Nash \mh{e}quilibrium \mh{p}roblem governed by a semi-linear elliptic PDE with distributed control\mh{:}{}
\begin{eq*}
&\text{minimize }  \calJ_i(u_i,u_{-i}) := \frac{1}{2}\int_{\omega_i} (y-y_d^i)^2 \dx + \frac{\alpha}{2}\int_{\omega_i} u_i^2 \dx \text{ over } u_i \in L^2(\omega_i)\\
&\text{subject to } a_i \leq u_i \leq b_i \text{ a.e. on } \omega_i, \ \underline{\psi} \leq y \leq \overline\psi \text{ a.e. on } \Omega \text{ and}\\
&\hspace*{4.31em}-\Delta y + y^3 = \sum_{i = 1}^4 B_i u_i \text{ in } \Omega,\ y = 0 \text{ on } \partial \Omega.
\end{eq*}
\end{ex}
Let $a_i = -32$ and $b_i = 32$ and choose the regularization parameter $\alpha = 10^{-5}$. We use the subdomain partition depicted in \refer{fig:subdomains}. As upper and lower obstacle for the state we take $\underline\psi = 0$ together with $\overline\psi = 0.3$ \msa{being}{} the upper and lower pointwise bounds of $y_d$ on \mh{all of}{} $\Omega$.
To use the results of this work we formulate this example in a way suitable to the framework established in \refer{ass:specializedframework}.
The operators $B_i \in \mathcal{L}(L^2(\omega_i), H^{-1}(\Omega))$ are just the extension of a function on $\omega_i$ by zero on the whole domain combined with \msa{the}{} embedding of $L^2(\Omega)$ in $H^{-1}(\Omega)$.
The operator $A$ will be defined via the relation
\begin{eq*}
\langle A(y), z \rangle_{H_0^1\msa{(\Omega)}, H^{-1}\msa{(\Omega)}} = (\nabla y, \nabla z)_{L^2\msa{(\Omega;\R^d)}} + (y^3, z)_{L^2\msa{(\Omega)}}
\end{eq*}
for sufficiently regular arguments and test functions. \msa{It is straightforward to show existence and uniqueness of a solution in $H_0^1(\Omega)$ with test functions in $H^{-1}(\Omega)$.}{} However, to discuss the choice of suitable function spaces, we use regularity theory for PDEs. \msa{Therefore, set}{} $Y := W^{1,r}(\Omega)$ and $X := C(\bar \Omega)$. \msa{Then, we}{} obtain \msa{the}{} continuous embedding $X \hookrightarrow Y$. To establish the complete continuity of the solution operator $S: U \rightarrow W^{1,r}_0(\Omega)$ we observe by \cite[Theorem 3.2.1.2]{bib:Grisvard} that there exists a constant $C > 0$ with 
\begin{eq*}
\|z\|_{H^2(\Omega)} \leq C\left( \|f\|_{L^2(\Omega)} + \|z\|_{L^2(\Omega)}\right),
\end{eq*}
where $f \in L^2(\Omega)$ and $z \in H^2(\Omega) \cap H_0^1(\Omega)$ is the solution of the Poisson problem 
\begin{eq*}
-\Delta z  &= f \text{ in } \Omega,\ z = 0 \text{ on } \partial \Omega.
\end{eq*}
In the following we will adapt the notation $a \lesssim b$ for $a,b > 0$, if there exists a constant $c > 0$ such that $a \leq c\cdot b$ holds.\\
Let a sequence $(u^n)_{n \in \N} \subseteq U = \prod_{i = 1}^4 L^2(\omega_i)$ with $u^n \rightharpoonup u$ be given and let \msa{$(y^n)_{n \in \N}$}{} denote the corresponding sequence of solutions \msa{of the state equations in the sense of $H_0^1(\Omega)$}. Then we obtain the following estimate
\begin{eq*}
\|y^n\|_{H^2\msa{(\Omega)}} &\lesssim \left\|\sum_{i = 1}^4 B_i u_i^n - (y^n)^3\right\|_{L^2(\Omega)} \leq \sum_{i = 1}^4 \|u_i^{\msa{n}}\|_{L^2(\omega_i)} + \|y^n\|^3_{L^6\msa{(\Omega)}}\\
&\lesssim \sum_{i = 1}^4 \|u_i^{\msa{n}}\|_{L^2(\omega_i)} + \|y^n\|^3_{H_0^1\msa{(\Omega)}} \lesssim \sum_{i = 1}^4 \|u_i^{\msa{n}}\|_{L^2(\omega_i)} + \left( \sum_{i = 1}^4 \|u_i^{\msa{n}}\|_{L^2(\omega_i)} \right)^3,
\end{eq*}
where we used the continuous embedding $H^1(\Omega) \hookrightarrow L^6(\Omega)$.
Therefore the boundedness of $(y^n)_{\msa{n \in \N}}$ in $H^2(\Omega)$ holds and for every
subsequence we can extract another subsequence with \msa{weak}{} limit point $y^* \in
H^2(\Omega)$.
Since $H^2(\Omega) \overset{c}{\hookrightarrow} H_0^1(\Omega)$ we obtain the strong
convergence in $H_0^1(\Omega)$ and by the continuity of $A$ that $A(y_n)
\rightarrow A(y^*)$ in $H^{-1}(\Omega)$.
By $A(y^n) = \sum_{i=1}^4 B_i u_i^n \rightarrow \sum_{i = 1}^4 B_i u_i^*$ in
$H^{-1}(\Omega)$ one observes $A(y^*) = \sum_{i = 1}^4 B_i u_i^*$. 
Consequently, since the choice of the first subsequence was arbitrary, we obtain
the weak convergence of the whole sequence in $H^2(\Omega)$.
Since the embedding in $H^2(\Omega) \overset{c}{\hookrightarrow} W^{1,r}(\Omega)$
is compact, we obtain eventually the complete continuity of the solution operator
$S:\prod_{i = 1}^4 L^2(\omega_i) \rightarrow W^{1,r}_0(\Omega)$. As space $W$ we
use $W = W^{-1,r}(\Omega)$ \msa{being}{} \mh{the dual of $W^{1,s}(\Omega)$ with $\frac{1}{r} + \frac{1}{s} = 1$}{}. \mh{We further define the}{} mapping $\calG: W_0^{1,r}(\Omega) \rightarrow
W^{1,r}(\Omega; \R^2)$ \msa{by}{} $\calG(y):= (y - \underline\psi, \overline\psi - y)$
\msa{with}{} the cone \msa{%
\begin{eq*}
K := \{ (z_1,z_2) \in W^{1,r}(\Omega;\R^2): z_1,z_2 \geq 0 \text{ a.e. on } \Omega\}.
\end{eq*}}
\msa{and take}{} $J_i^1(y) := \frac{1}{2}\|y - y_d^i\|_{L^2(\omega_i)}^2$ and $J_i^2(u_i) := \frac{\alpha}{2}\|u_i\|_{L^2(\omega_i)}^2$. \msa{Then}, $DJ_i^2(u_i) = \alpha u_i$ is even a weakly continuous mapping.\\
The fulfillment of the constraint qualification is a delicate question and is in the following assumed within the scope of this work. \ms{In the light of the embedding $W^{1,r}(\Omega) \hookrightarrow C(\bar\Omega)$ for $r > d$ a formulation in the space of continuous functions might be of advantage.}{}\\
Hence, we establish the penalization technique discussed in \refer{sec:penalization}
and \refer{sec:specialized}. For this \msa{sake}{} we introduce as penalty function \msa{the}{} Moreau--Yosida regularization \msa{with respect
to the $L^2$-norm}{} of the \msa{set}{} 
$\{z \in L^2(\Omega) : \underline\psi \leq z \leq \overline\psi  \}$ \msa{corresponding to the state constraint}{}\mh{. This}{} lead\mh{s}{} to
\begin{eq*}
\beta(y) := \frac{1}{2}\int_\Omega \left(y - \overline\psi\right)^{2+} \dx + \frac{1}{2}\int_\Omega \left(\underline\psi - y\right)^{2+} \dx.
\end{eq*}
As first order system, we derive by the use of \refer{thm:gnep_vep:path:conv} the following system:
\begin{sys}[\label{sys:exp:elliptic_semi_distr:firstorder}]
u_i &= \operatorname{Proj}_{U_\ad}\left(-\frac{1}{\alpha}\mathbbm{1}_{\omega_i} p_i\right) &&\text{in } \Omega,\\
-\Delta y + y^3 &= \sum_{i = 1}^4 B_i u_i &&\text{in } \Omega,\\
y &= 0  &&\text{on } \partial \Omega,\\
-\Delta p_i + 3y^2 p_i &= \mathbbm{1}_{\omega_i} (y-y_d^i) + \gamma (y - \overline\psi)^+ - \gamma (\underline\psi - y)^+  &&\text{in } \Omega,\\
p_i &= 0 &&\text{on } \partial \Omega.
\end{sys}
In this setup the cooperative comparison \msa{problem}{} reads as
\begin{eq}\label{eq:ex:semielliptic:distr:coop}
&\text{minimize} &&\sum_{i = 1}^4 \left(J_i^1(y) + J_i^2(u_i)\right) = \frac{1}{2}\|y-y_d\|^2_{L^2(\Omega)} + \sum_{i = 1}^4 \frac{\alpha}{2}\|u_i\|_{L^2(\omega_i)}^2\\
&\text{subject to} &&u_i \in L^2(\omega_i), a_i \leq u_i \leq b_i, \underline\psi \leq y \leq \overline\psi \text{ a.e. on } \Omega \text{ and}\\
&\phantom{\text{subject to}}&&-\Delta y + y^3 = \sum_{i = 1}^4 B_i u_i \text{ in } \Omega,\ y = 0 \text{ on } \partial \Omega.
\end{eq}
The first order system can be derived by standard techniques in optimization or by the results in this paper for the case $N=1$ with exactly one player. This leads to
\begin{sys}[\label{sys:exp:elliptic_semi_distr:coop:firstorder}]
u_i &= \operatorname{Proj}_{U_\ad}\left(-\frac{1}{\alpha}\mathbbm{1}_{\omega_i} p_i\right) &&\text{in } \Omega,\\
-\Delta y + y^3 &= \sum_{i = 1}^4 B_i u_i &&\text{in } \Omega,\\
y &= 0 &&\text{on } \partial \Omega,\\
-\Delta p + 3y^2 p &= (y-y_d) + \gamma (y - \overline\psi)^+ - \gamma (\underline\psi - y)^+  &&\text{in } \Omega,\\
p &= 0 &&\text{on } \partial \Omega.
\end{sys}
\msa{Comparing \refer{sys:exp:elliptic_semi_distr:coop:firstorder} and \refer{sys:exp:elliptic_semi_distr:firstorder} it is worth noting that $(y-y_d) = \sum_{i = 1}^4 (y-y_d^i)\mathbbm{1}_{\omega_i}$ by our partition of the domain $\Omega$.}{}
In the following numerical experiments, we choose as parameters $C_{\mathrm{Path}} = 10^{-5}$ and $\eps = 10$ for the path-following strategy for both \msa{---}{} the cooperative as well as the non-cooperative case.\\
The plots of the results are depicted below. Since the controls are basically truncations of the adjoint states (see \refer{sys:exp:elliptic_semi_distr:coop:firstorder})
we restrict ourselves to \msa{the depiction of}{} the controls (\msa{in}{} \refer{fig:exp:elliptic_semi_distr:controls}) and states (\msa{in}{} \refer{fig:exp:elliptic_semi_distr:states}).
To illustrate the controls we depicted the combined control $u = \sum_{i = 1}^4 \mathbbm{1}_{\omega_i} u_i $. In both cases the different control regions $\omega_i$ can be clearly \msa{recognized}. The cooperative case admits a solution that is in $H_0^1(\Omega)$, whereas in the non-cooperative case only piecewise $H^1$-regularity for the combined controls is guaranteed. In the latter case the boundaries between the control regions are clearly visible due to \mh{activity}{} of the box constraint.\\
In \refer{fig:exp:elliptic_semi_distr:gamma} the update history for both cases is
depicted. Both iterations are terminated as soon as the penalty function \mh{drops below}{} a threshold of $10^{-15}$ or the penalty parameter $\gamma$ exceeds the
value $10^8$. The direct comparison of the two graphs indicate that the update mechanism is less aggressive in the non-cooperative case. \msa{On the one hand, this might be related}{} to the explicit presence of the computed equilibrium in the \mh{value}{}
function (see also \refer{defn:gnep_vep:penalized:merit}).
\msa{On}{} the other hand, the \mh{game}{} leads the player to \mh{compete}{} in their decisions near boundary regions \mh{with}{} another player. Hence\msa{,}{} the benefit from the competition
rewarded via \msa{the}{} tracking term might be more valuable than the punishment implied by
the penalty functional. This phenomenon also translates to the states\msa{:}{} As it can be
seen from \refer{fig:exp:elliptic_semi_distr:states}\mh{,}{} the interfaces between
the different regions are more distinct in the non-cooperative case. However, in the regions $\omega_1$ and $\omega_2$ the state falls \mh{below}{} respectively
exceed\mh{s}{} the desired values near the regions $\omega_4$ and $\omega_3$.
\begin{figure}[H]
\begin{tabular}{cc}
\resizebox{0.47\textwidth}{!}{%
\includestandalone{nash_elliptic_semi_distr/objectives}} &
\resizebox{0.47\textwidth}{!}{%
\includestandalone{nash_elliptic_semi_distr_coop/objectives}}
\end{tabular}
\caption{Plot of the summed objectives. Left: Results for Example \ref{ex:semielliptic:distr}, Right: Results for the cooperative version in \eqref{eq:ex:semielliptic:distr:coop}.}\label{fig:exp:elliptic_semi_distr:obj}
\end{figure}

\begin{figure}[H]
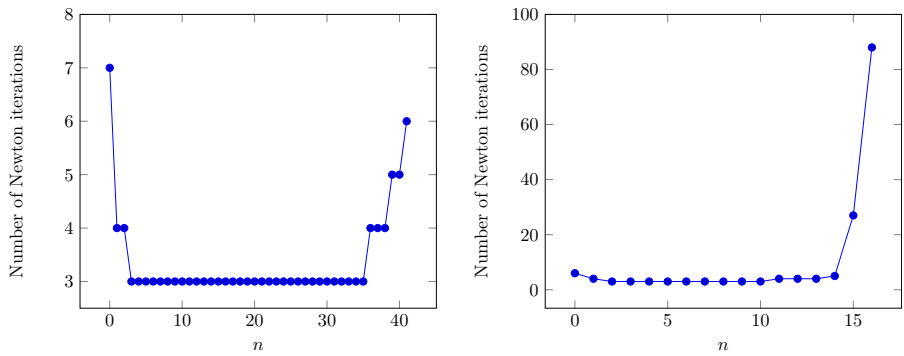

\begin{tabular}{cc}
\resizebox{0.47\textwidth}{!}{%
\includestandalone{nash_elliptic_semi_distr/newton_iter}} &
\resizebox{0.47\textwidth}{!}{%
\includestandalone{nash_elliptic_semi_distr_coop/newton_iter}}
\end{tabular}
\ms{\caption{Plot of the \textcolor{msa}{number} of Newton iterations. Left: Results for Example \ref{ex:semielliptic:distr}, Right: Results for the cooperative version in \eqref{eq:ex:semielliptic:distr:coop}.}}\label{fig:exp:elliptic_semi_distr:newton}
\end{figure}

\noindent The behaviour of the sum of objective values (without the penalty
function) is depicted in \refer{fig:exp:elliptic_semi_distr:obj}. They both follow
the very same pattern that is well known from optimization problems.
A \msa{frequently}{} used concept in game theory related to discrete mathematics, e.g.\mh{,}{} in the context of selfish routing, is called the \emph{price of anarchy} (see e.g. \cite[Definition 20.4]{bib:NisanRoughgardenAlgorithmicGameTheory}) defined by
\begin{eq*}
\operatorname{PoA} = \frac{\sup_{u \in \mathrm{Equil}} \sum_{i = 1}^N \calJ_i(u)}{\min_{u \in \calF} \sum_{i = 1}^N \calJ_i(u)}.
\end{eq*}
Here\msa{,}{} $\calF$ denotes again the joint constraint set for the problem, and `$\mathrm{Equil}$' denotes the set of all Nash equilibria. A calculation with respect to the last iterate yields $\operatorname{PoA} \geq 1.1641$ as a lower approximation. 
\ms{For the Newton iterates we note, that the application of the damping strategy has only been applied in the cooperative case and only for the last two $\gamma$-iterations. The reason for this behaviour lies in the more conservative update behaviour of the penalty parameter for the Nash game.}{}

\begin{figure}[H]
\includegraphics[trim = {750 325 1050 450}, clip, width = 0.45\textwidth]{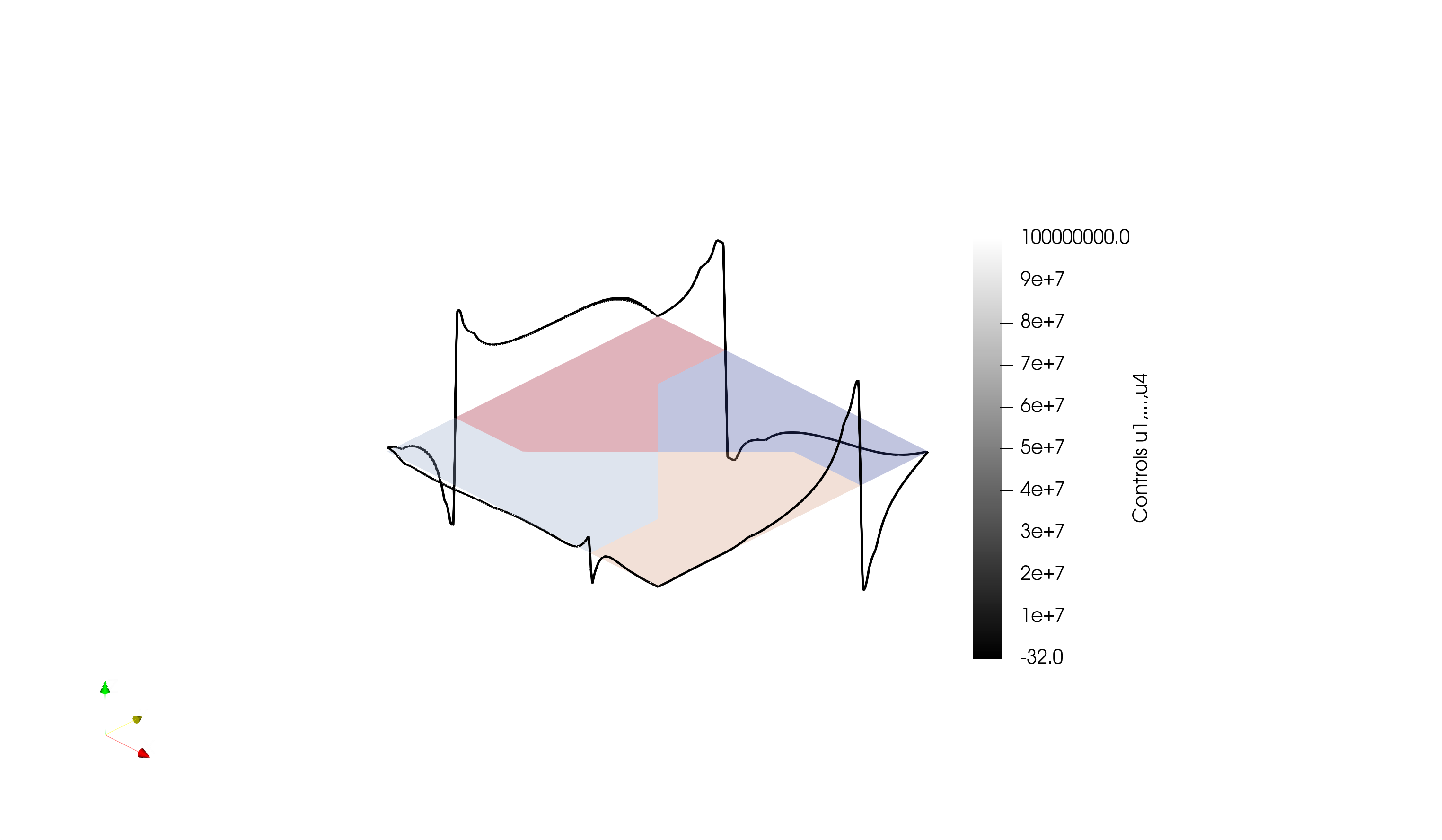}
\hspace*{0.05\textwidth}
\includegraphics[trim = {750 325 1050 450}, clip, width = 0.45\textwidth]{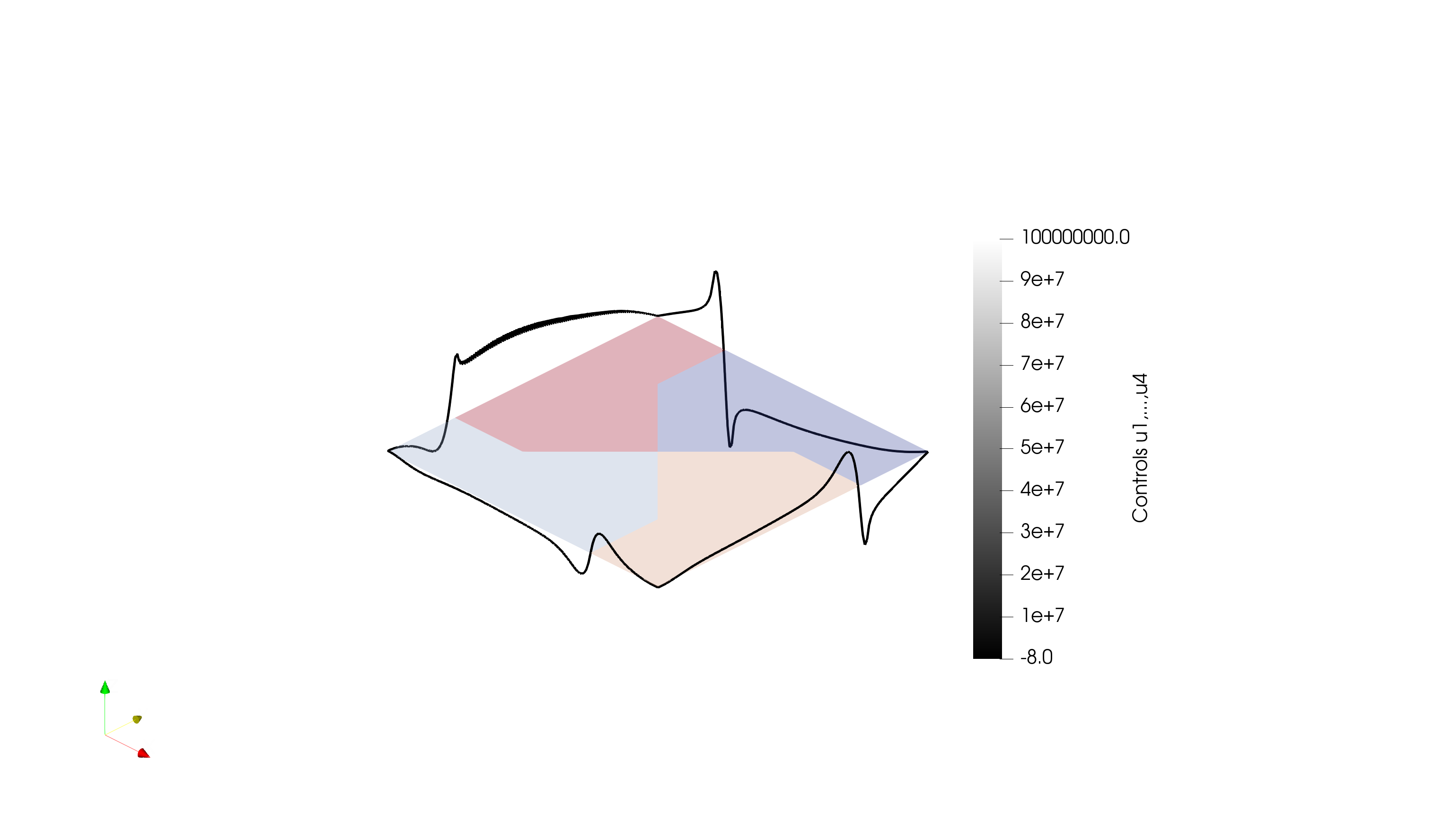}
\ms{\caption{Plot of the controls. Left: Results for Example \ref{ex:semielliptic:bdry} with values between $-32$ and $32$. Right: Results for the cooperative version in \eqref{eq:ex:semielliptic:bdry:coop} with values between $-6.5$ and $5.5$.}}\label{fig:exp:elliptic_semi_bdry:controls}
\end{figure}

\begin{figure}[H]
\includegraphics[trim = {850 200 500 200}, clip, width=0.45\textwidth]{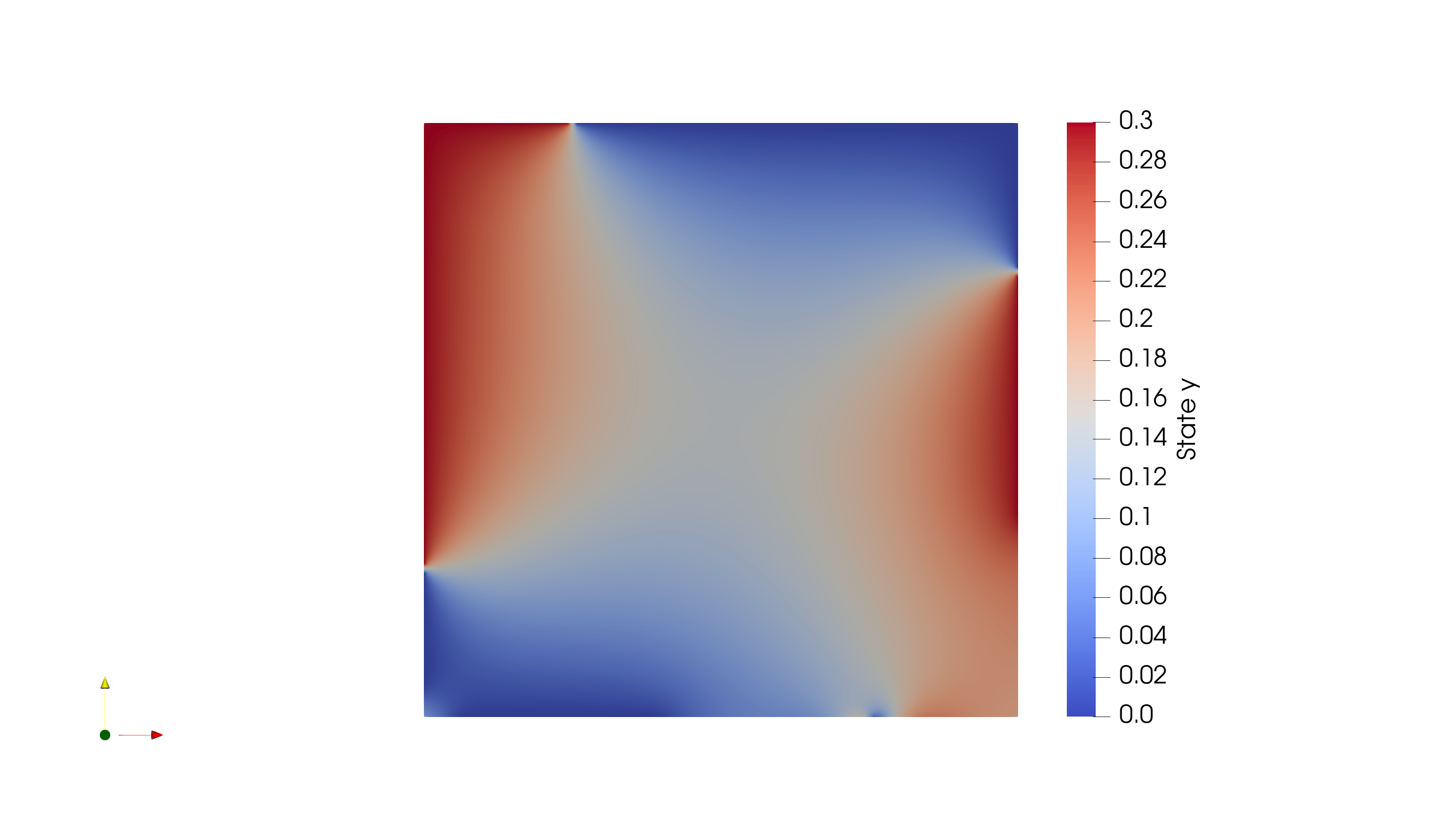}
\hspace*{0.05\textwidth}
\includegraphics[trim = {850 200 500 200}, clip, width=0.45\textwidth]{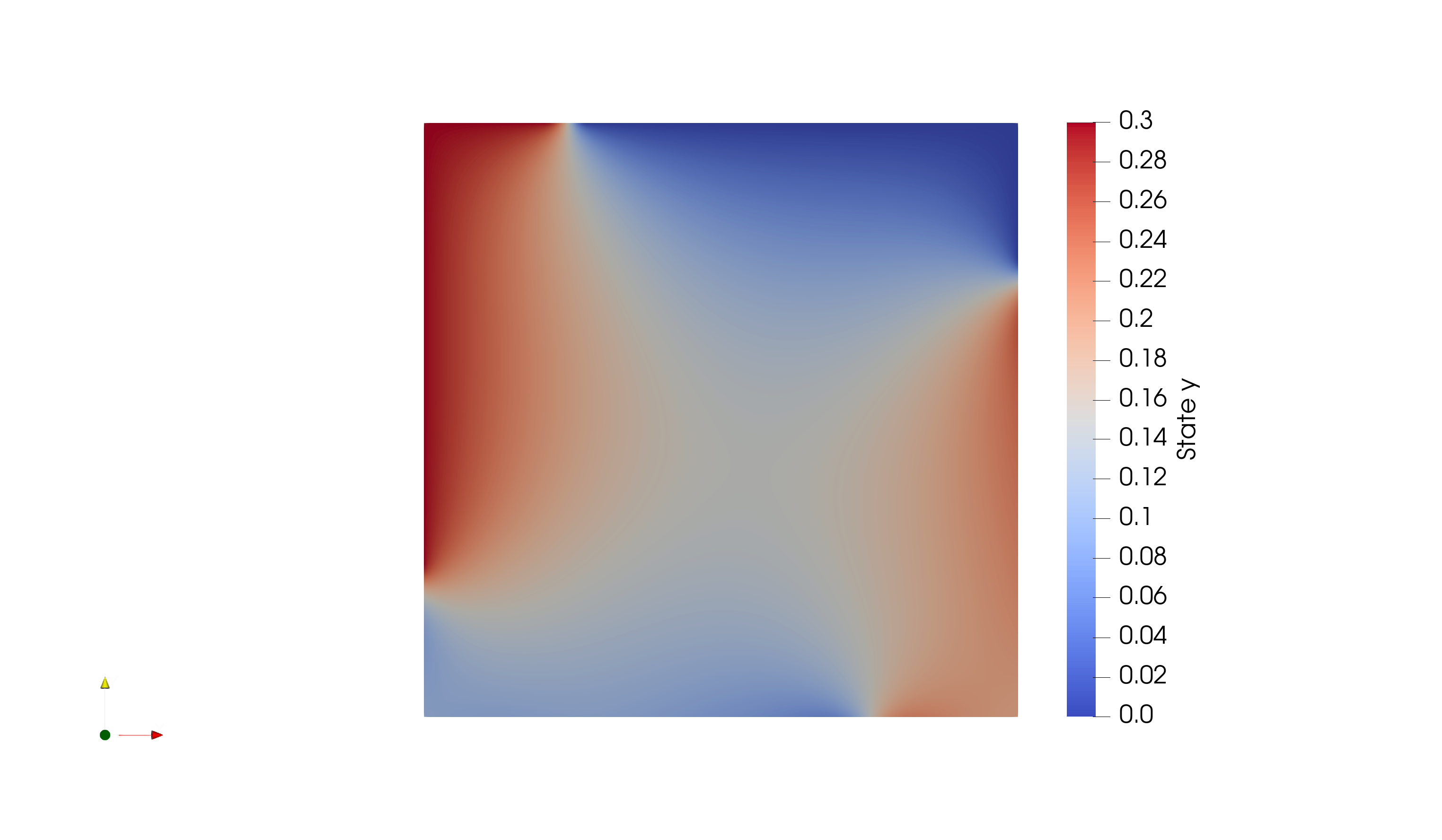}
\caption{Plot of the states. Left: Results for Example \ref{ex:semielliptic:bdry}, Right: Results for the cooperative version \eqref{eq:ex:semielliptic:bdry:coop}.}\label{fig:exp:elliptic_semi_bdry:states}
\end{figure}

\begin{ex}\label{ex:semielliptic:bdry}
Consider the following \mh{g}eneralized Nash \mh{e}quilibrium \mh{p}roblem governed by a semi-linear elliptic PDE with boundary control
\begin{eq*}
&\text{minimize} &&J_i(u_i,u_{-i}) := \frac{1}{2}\int_{\omega_i} (y-y_d^i)^2 \dx + \frac{\alpha}{2}\int_{\Gamma_i} u_i^2 \dS \text{ over } u_i \in L^2(\Gamma_i)\\
&\text{subject to} &&a_i \leq u_i \leq b_i \text{ a.e. on } \Gamma_i, \underline{\psi} \leq y \leq \overline\psi \text{ a.e. on } \Omega \text{ and}\\
&\hspace*{4.775em} &&-\Delta y + y = 0, \text{ in } \Omega,\ \frac{\partial y}{\partial \nu} + y^3 = \sum_{i = 1}^N B_i u_i \text{ on } \partial \Omega.
\end{eq*}
\end{ex}
\msa{We choose again}{} $a_i = -32$ and $b_i = 32$ and as regularization parameter $\alpha = 10^{-5}$. Together with the subdomain partition depicted in \refer{fig:subdomains} we define as control region for player $i$ the set $\Gamma_i := \partial \Omega \cap \partial \omega_i$. As upper and lower obstacle for the state again $\underline\psi = 0$ and $\overline\psi = 0.3$ are taken.\\
To utilize our framework we check the conditions in \refer{ass:specializedframework}. \msa{For this sake,}{} we define again $X= C(\bar \Omega)$ as well as $Y:= W^{1,r}(\Omega)$ with $r \in (2,3)$. Using the regularity result \cite[Theorem 4]{bib:Savare} one can deduce from its proof for $s \in (-\frac{1}{2},\frac{1}{2})$ that for $f \in L^2(\Omega)$ and $g \in H^{-\frac{1}{2}+s}(\Omega)$ \msa{the}{} solution of the PDE
\begin{eq*}
-\Delta z + z &= f \text{ in } \Omega,\ \frac{\partial z}{\partial \nu} = g \text{ in } \partial \Omega\mh{,}{}
\end{eq*}
\mh{satisfies}{} the estimate 
\begin{eq*}
\|z\|_{H^{1+s}\msa{(\Omega)}} \lesssim \|f\|_{L^2\msa{(\Omega)}} + \|g\|_{H^{-\frac{1}{2}+s}(\partial \Omega)}\mh{,}{}
\end{eq*}
with a constant depending on $s$ as well as the domain $\Omega$.\\
Taking a weakly convergent sequence $(u^n)_{n\in \N} \subseteq \prod_{i = 1}^4 L^2(\Gamma_i)$ with limit $u^* \in U$ leads to a sequence of states $(y^n)_{n \in \N} \subseteq H^1(\Omega)$ with $y^n = S(u^n)$. 
By using the embedding $H^{\frac{1}{2}}(\partial \Omega) \hookrightarrow L^p(\partial \Omega)$ for all $p \in (1, +\infty)$ we infer that $\left(\operatorname{tr}_{\partial \Omega} y\right)^3 \in L^2(\partial \Omega) \hookrightarrow H^{-\frac{1}{6}}(\partial \Omega)$. With $s = \frac{1}{3}$ we obtain
\begin{eq*}
\|y^n\|_{H^{\frac{4}{3}}(\Omega)} &\lesssim \left\|\sum_{i = 1}^4 B_i u_i^n - \left(\operatorname{tr}_{\partial \Omega} y^n\right)^3\right\|_{H^{-\frac{1}{6}}(\partial \Omega)}\\
&\lesssim \sum_{i = 1}^4 \|u_i^n\|_{L^2(\Gamma_i)} + \|\operatorname{tr}_{\partial \Omega} y^n\|_{L^6(\partial \Omega)}^3\\
&\lesssim \sum_{i = 1}^4 \|u_i^n\|_{L^2(\Gamma_i)} + \|y^n\|_{H^1(\Omega)}^3\\
&\lesssim \sum_{i = 1}^4 \|u_i^n\|_{L^2(\Gamma_i)} + \left( \sum_{i = 1}^4 \|u_i^n\|_{L^2(\Gamma_i)}\right)^3.
\end{eq*}
\mh{From}{} the boundedness of the sequence $(u^n)_{n \in \N}$ we \mh{also}{} \msa{deduce}{} 
the boundedness of $(y^n)_{n \in \N}$ in $H^{\frac{4}{3}}(\Omega)$. By reflexivity we \mh{further}{} infer along every subsequence the existence of weakly convergent subsequence in $H^{\frac{4}{3}}(\Omega)$. Using the embedding $H^{\frac{1}{3}}(\Omega) \overset{c}{\hookrightarrow} L^r(\Omega)$ for all $r \in (2,3)$. Hence we obtain $H^{\frac{4}{3}}(\Omega) \overset{c}{\hookrightarrow} W^{1,r}(\Omega)$ as well and deduce the strong convergence in the latter space as well as in $H^{1}(\Omega)$. By the continuity of the operator $A: H^1(\Omega) \rightarrow H^{-1}(\Omega)$ we infer $A(y^n) \rightarrow A(y)$ in $H^{-1}(\Omega)$ and by the weak convergence of $u^n$ as well $\left\langle \sum_{i = 1}^4 u^n_i, \operatorname{tr}_{\partial \Omega}(\ccdot) \right\rangle_{H^{-\frac{1}{2}}\msa{(\Omega)},H^{\frac{1}{2}}\msa{(\Omega)}} \rightarrow \left\langle \sum_{i = 1}^4 u_i, \operatorname{tr}_{\partial \Omega}(\ccdot) \right\rangle_{H^{-\frac{1}{2}\msa{(\Omega)}},H^{\frac{1}{2}}\msa{(\Omega)}}$ in $H^{-1}(\Omega)$. Hence, \msa{$y$}{} is the unique solution \msa{of}{} the state equation with respect to $u$. Thus, the whole sequence $(y^n)_{n \in \N}$ converges weakly in $H^{\frac{4}{3}}(\Omega)$ and strongly in $W^{1,r}(\Omega)$ for $r \in (2,3)$. Eventually, we deduce the complete continuity of the solution operator.\\
As in the previous \refer{ex:semielliptic:distr} we take $\calG: W_0^{1,r}(\Omega) \rightarrow W^{1,r}(\Omega; \R^2)$ defined by $\calG(y) := (y - \underline\psi, \overline\psi - y)$ along with \msa{%
\begin{eq*}
K := \{(z_1,z_2) \in W^{1,r}(\Omega;\R^2): z_1, z_2 \geq 0 \text{ a.e. on } \Omega\}.
\end{eq*}}
And again, the constraint qualification \refer{eq:gnep:path:uniformrobinson} is assumed to hold true. For the penalization of the state constraint we introduce
\begin{eq*}
\beta(y) := \frac{1}{2}\int_\Omega \left(y - \overline\psi\right)^{2+} \dx + \frac{1}{2}\int_\Omega \left(\underline\psi - y\right)^{2+} \dx
\end{eq*}
as well. The first order system of the penalized problem can be derived with \mh{the help of}{} \refer{cor:gnep_vep:penalized:firstorder:special} as in the previous example\mh{:}{}
\begin{sys}
u_i &= \operatorname{Proj}_{U_\ad}\left(- \frac{1}{\alpha} \mathbbm{1}_{\Gamma_i} \operatorname{tr}_{\partial \Omega}p_i\right) &&\text{in }\Omega,\\
-\Delta y + y &= 0 &&\text{in } \Omega,\\
\frac{\partial y}{\partial \nu} + y^3 &= \sum_{i = 1}^N B_i u_i &&\text{on } \partial \Omega,\\
-\Delta p_i + p_i &= \mathbbm{1}_{\omega_i}(y - y_d^i) + \gamma (y - \overline\psi)^+ - \gamma (\underline\psi - y)^+ &&\text{in } \Omega,\\
\frac{\partial p_i}{\partial \nu} + 3y^2 p_i &= 0 &&\text{on } \partial \Omega. & 
\end{sys}
Hence, the cooperative comparison problem reads as
\begin{eq}\label{eq:ex:semielliptic:bdry:coop}
&\text{minimize} &&\sum_{i = 1}^4 \left(J_i^1(y) + J_i^2(u_i)\right) = \frac{1}{2} \|y-y_d\|_{L^2(\Omega)}^2 + \sum_{i = 1}^4 \frac{\alpha}{2}\|u_i\|_{L^2(\Gamma_i)}^2 \\
&\text{subject to} &&u_i \in L^2(\Gamma_i), a_i \leq u_i \leq b_i \text{ a.e. on } \Gamma_i, \underline{\psi} \leq y \leq \overline\psi \text{ a.e. on } \Omega \text{ and}\\
&{}&&-\Delta y + y = 0 \text{ in } \Omega,\ \frac{\partial y}{\partial \nu} + y^3 = \sum_{i = 1}^N B_i u_i \text{ on } \partial \Omega.
\end{eq}
The first order system can as in the previous example. This leads to
\begin{sys}
u_i &= \operatorname{Proj}_{U_\ad}\left(- \frac{1}{\alpha} \mathbbm{1}_{\Gamma_i} \operatorname{tr}_{\partial \Omega}p\right) &&\text{in } \Omega,\\
-\Delta y + y &= 0 &&\text{in } \Omega,\\
\frac{\partial y}{\partial \nu} + y^3 &= \sum_{i = 1}^4 B_i u_i &&\text{on } \partial \Omega,\\
-\Delta p + p &= (y - y_d) + \gamma (y - \overline\psi)^+ - \gamma (\underline\psi - y)^+ &&\text{in } \Omega,\\
\frac{\partial p}{\partial \nu} + 3 y^2 p &= 0 &&\text{on } \partial \Omega.
\end{sys}
Below, the results of the experiments are depicted. For better orientation in \refer{fig:exp:elliptic_semi_bdry:controls} the distinction of the domain $\Omega$ has been added in the colouring of the desired state of \refer{fig:subdomains}. Similarly to the previous example the competition between the players \msa{induces}{} in \mh{activity}{} of the control constraints near the boundary points between the control boundaries $\Gamma_i$. In contrast, the cooperative case \mh{does}{} not reach the \mh{bounds}{}. (Note the different scaling.) This drastic difference in behaviour affects the states as well. As one can see \msa{in \refer{fig:exp:elliptic_semi_bdry:states}}{}, the desired state is not even nearly met on $\omega_1$ and $\omega_2$ in the non-cooperative case, but it is observably \msa{better}{} met in the cooperative one.
The influence of this competition can \mh{also}{} be seen in the update behaviour.
As before\mh{,}{} the non-cooperative case has \msa{less}{} aggressive update\mh{s}, but
\mh{when}{} compared to distributed control the difference is even more profound.
The behaviour of the summed objectives is in the light of the previous example
counterintuitive. However, in the optimization case it is clear that due to
imposing the state constraint the objective values are a non-decreasing sequence.
For Nash equilibria such a result is not available and is in fact not a mandatory
behaviour as it can be seen in this example.
This observation can therefore be interpreted as an instance of a \emph{Braess paradoxon} (cf. \cite{bib:BraessParadox}). Using the final iterates of the algorithm we obtain
for the price of anarchy $\operatorname{PoA} \geq 1.29890$ as lower estimate. So
in comparison to \refer{ex:semielliptic:distr} the price of anarchy is considerably
bigger, which indicates a more competitive environment and might hence be an
explanation for the behaviour observed in \refer{fig:exp:elliptic_semi_bdry:obj}.
\ms{For the Newton iterations the established damping rule has only been used towards the end for large $\gamma$-updates.}{}

\begin{figure}[H]
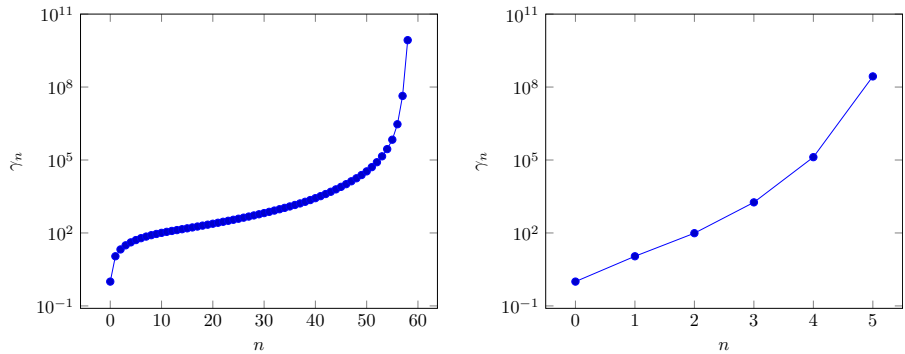

\begin{tabular}{cc}
\resizebox{0.47\textwidth}{!}{%
\includestandalone{nash_elliptic_semi_bdry/gamma}} &
\resizebox{0.47\textwidth}{!}{%
\includestandalone{nash_elliptic_semi_bdry_coop/gamma}}
\end{tabular}
\caption{Plot of the $\gamma$-updates. Left: Results for Example \ref{ex:semielliptic:bdry}, Right: Results for the cooperative version \eqref{eq:ex:semielliptic:bdry:coop}.}\label{fig:exp:elliptic_semi_bdry:gamma}
\end{figure}

\begin{figure}[H]
\begin{tabular}{cc}
\resizebox{0.47\textwidth}{!}{%
\includestandalone{nash_elliptic_semi_bdry/objectives}} &
\resizebox{0.47\textwidth}{!}{%
\includestandalone{nash_elliptic_semi_bdry_coop/objectives}}
\end{tabular}
\caption{Plot of the summed objectives. Left: Results for Example \ref{ex:semielliptic:bdry}, Right: Results for the cooperative version \eqref{eq:ex:semielliptic:bdry:coop}.}\label{fig:exp:elliptic_semi_bdry:obj}
\end{figure}

\begin{figure}[H]
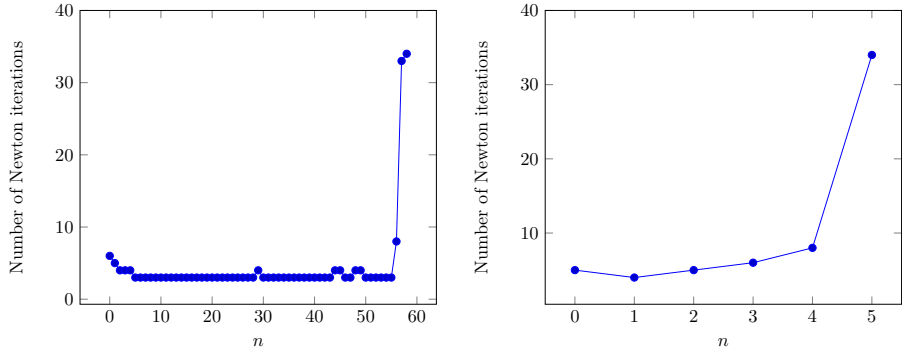

\begin{tabular}{cc}
\resizebox{0.47\textwidth}{!}{%
\includestandalone{nash_elliptic_semi_bdry/newton_iter}} &
\resizebox{0.47\textwidth}{!}{%
\includestandalone{nash_elliptic_semi_bdry_coop/newton_iter}}
\end{tabular}
\ms{\caption{Plot of the \textcolor{msa}{number} of Newton iterations. Left: Results for Example \ref{ex:semielliptic:bdry}, Right: Results for the cooperative version in \eqref{eq:ex:semielliptic:bdry:coop}.}}\label{fig:exp:elliptic_semi_bdry:newton}
\end{figure}

\ms{%
\section{Conclusion}
In the present work a framework for the analytical and numerical treatment of generalized Nash equilibrium problems has been developed for different levels of abstraction. Special emphasis has been put to the application of our results to games derived from optimal control problems involving non-linear partial differential equations. \msa{As}{} approximation technique a penalization concept for equilibria has been discussed along with the treatment of the multipliers in the first order systems. The results have been applied to two GNEPs with optimal control problems with distributed respectively boundary controls. The comparison of these games with their cooperative counterparts has been used to inspect and visualize the influence of competition as well as the enforcement of joint constraints to the choice of strategies of the players.
}{}

\section*{Acknowledgements}
The author thanks Prof. Dr. Michael Hinterm\"uller for valuable remarks and comments that helped to improve this article.

\appendix
\section{Appendix}

\begin{proof}[Proof of \refer{lem:uniformfrechet:properties}]
Consider for $t \in \R\backslash\{0\}$ the sequence of mappings $T_t: X \times X \rightarrow Y$ by
\begin{eq*}
T_t(x,h):= \frac{T(x + th) - T(x)}{t}.
\end{eq*}
Using the uniform Fr\'echet differentiability one obtains that for all bounded subsets $M \subseteq X$ and every 
$\varepsilon > 0$ there exists a $\delta > 0$ such that
\begin{eq*}
\frac{\| T(x' + h') - T(x') - DT(x')h' \|_Y}{\|h'\|_X} \leq \varepsilon
\end{eq*}
for all $x' \in M$ and  $h' \in X$ with $\|h'\|_X \leq \mh{\delta}$. Hence, taking
\mh{$R>0$ and}{} $t\mh{>0}{}$ sufficiently small yields
\begin{eq*}
\sup_{\substack{x' \in M, h \in X,\\ \|h\|_X \leq R}} \|T_t(x,h) - DT(x)h\|_{Y} &= \sup_{\substack{x' \in M, h \in X,\\ \|h\|_X \leq R}} \frac{\|T(x + th) - T(x) - tDT(x)h\|_Y}{|t|}\\
&\leq \varepsilon R,
\end{eq*}
which implies the uniform convergence of $(x,h) \mapsto T_t(x,h)$ towards $(x,h) \mapsto DT(x)h$ on every
 bounded subset of $X \times X$ for $t \rightarrow 0$ (\mh{cf.}{} also \cite[Chapter 10, footnote 5]{bib:Ehrmann}).\\
Let first $T$ be weakly continuous. Take sequences $x_n \rightharpoonup x$ and $h_n \rightharpoonup h$ in $X$ as well as an arbitrary $y^* \in Y^*$. By the boundedness of \msa{the sequences $(x_n)_{n \in \N}, (h_n)_{n \in \N}$}{} and the uniform Fr\'echet differentiability we obtain
\begin{eq*}
|\langle y^*, DT(x_n)h_n &- DT(x)h \rangle_{\msa{Y^*,Y}}| \msa{\leq} |\langle y^*, DT(x_n)h_n - T_t(x_n,h_n) \rangle_{\msa{Y^*,Y}}|\\
&+ |\langle y^*, T_t(x_n,h_n) - T_t(x,h) \rangle_{\msa{Y^*,Y}}| + |\langle y^*, DT(x,h) - T_t(x,h) \rangle_{\msa{Y^*,Y}}|\\
&\leq 2\varepsilon + |\langle y^*, T_t(x_n,h_n) - T_t(x.h) \rangle_{\msa{Y^*,Y}}|.
\end{eq*}
Using the weak continuity of $T$ one obtains for fixed $t \in \R\backslash \{0\}$ also the weak continuity of $T_t$ by its definition. \mh{This yields}{}
\begin{eq*}
0 \leq \limsup_{n \rightarrow \infty} |\langle y^*, DT(x_n)h_n - DT(x)h \rangle_{\msa{Y^*,Y}}| \leq 2\varepsilon
\end{eq*}
for all $\varepsilon > 0$ and\mh{,}{} thus\mh{,}{} the weak convergence \mh{of $DT(x_n)h_n$}{}.\\
Taking arbitrary sequences $x_n \rightharpoonup x$ in $X$ and $y_n^* \rightarrow y^*$ in $Y^*$ we see for all $\mh{h}{} \in X$ that
\begin{eq*}
\langle DT(x_n)^*y_n^*,\mh{h}{} \rangle_{\msa{X^*,X}} = \langle y_n^*, DT(x_n)\mh{h}{} \rangle_{\msa{Y^*,Y}} \rightarrow \langle y^*,DT(x)\mh{h}{} \rangle_{\msa{Y^*,Y}}
\end{eq*}
using the previously proven weak convergence $DT(x_n)\mh{h}{} \rightharpoonup DT(x)\mh{h}{}$ in $Y$.\\
In the case of $T$ being completely continuous take again sequences $x_n \rightharpoonup x$ and $h_n \rightharpoonup h$ in $X$ as well as an arbitrary $y^* \in Y^*$. Using the uniform Fr\'echet differentiability we see
\begin{eq*}
\| DT(x_n)h_n &- DT(x)h\|_Y \leq \| DT(x_n)h_n - T_t(x_n,h_n)\|_Y\\
&+ \| T_t(x_n,h_n) - T_t(x,h)\|_Y + \| DT(x)h- T_t(x,h)\|_Y\\
&\leq 2\varepsilon + \|T_t(x_n,h_n) - T_t(x,h)\|_Y.
\end{eq*}
The complete continuity of $T$ yields as well for fixed $t \in \R\backslash\{0\}$ the complete continuity of $T_t$  and hence we obtain 
\begin{eq*}
0 \leq \limsup_{n \rightarrow \infty} \| DT(x_n)h_n - DT(x)h \|_Y \leq 2\varepsilon
\end{eq*}
for $\varepsilon > 0$, thus yielding the strong convergence.\\
For the remaining assertion we use an indirect proof. Assume, there exist $x_n \rightharpoonup x$ and $y_n^* \rightharpoonup^* y^*$ and $\eps > 0$ such that $\|DT(x_n)^* y_n^* - DT(x)^*y^*\|_{X^*} \geq 2\eps$ holds for all $n \in \mathbb{N}$. Then, for every $n$ there exists a $v_n \in X$ with $\|v_n\|_X \leq 1$ such that
\begin{eq*}
\langle DT(x_n)^* y_n^* - DT(x)^*y^*, v_n \rangle_{\msa{X^*,X}} \geq \frac{1}{2}\|DT(x_n)^* y_n^* - DT(x)^*y^*\|_{X^*} \geq \eps\msa{.}
\end{eq*}
By \mh{the}{} Banach-Alaoglu \mh{theorem}{} we extract a weakly convergent subsequence from \msa{$(v_n)_{n \in \N}$}{} with limit $v$ (not relabeled). Using the previously shown complete continuity of $(x,h) \mapsto DT(x)h$ yields
\begin{eq*}
\eps \leq \langle DT(x_n)^* y_n^* &- DT(x)^*y^*, v_n \rangle_{\msa{X^*,X}} = \langle y_n^*, DT(x_n) v_n \rangle_{\msa{Y^*,Y}} - \langle y^* , DT(x) v_n \rangle_{\msa{Y^*,Y}}\\ &\rightarrow \langle y^*, DT(x) v \rangle_{\msa{Y^*,Y}} - \langle y^*, DT(x) v \rangle_{\msa{Y^*,Y}} = 0
\end{eq*}
and thus the requested contradiction.
\end{proof}

\begin{proof}[Proof of \refer{lem:bddderivative}]
Suppose the contrary. Then, there exists a positive bound $R > 0$ together with a sequence $(x_n)_{n \in \N}$,
$\|x_n\|_X \leq R$ for all $n \in \N$ and $\|DT(x_n)\|_{\calL(X,Y)} \geq 2n$.
Hence, there exists $(h_n)_{n \in \N}$ with $\|h_n\|_X \leq 1$ and $\|DT(x_n) h_n\|_Y \geq n$. By the reflexivity of $X$ one extracts weakly convergent subsequences $x_n \rightharpoonup x$ and $h_n \rightharpoonup h$ in $X$.
\refer{lem:uniformfrechet:properties}\ref{enum:state:uniformfrechet:weaklycont} yields the weak convergence of $DT(x_n)h_n 
\rightharpoonup DT(x)h$ in $Y$ and thus the boundedness of $\|DT(x_n)h_n\|_Y$---a contradiction.
\end{proof}

\begin{proof}[Proof of \refer{lem:uniformfrechet:inverse}]
Let a bounded subset $L \subset Y$ be given. Since $S$ is assumed to be a bounded operator\mh{,}{} the set $M:= S(L+\B_Y) \subset X$ is bounded as well. Using the assumption on the inverses of the first derivatives, \mh{let}{}
\begin{eq*}
0 < B := \sup_{x \in M} \|DT(x)^{-1}\|_{\calL(Y,X)} < +\infty.
\end{eq*}
By the uniform Fr\'echet differentiability there exists for all $\varepsilon > 0$ a positive number $\delta \in (0,1)$ such that for all $x \in M$ and $h \in X$ with $\|h\|_X \leq \delta$ \mh{it}{} holds \mh{that}{}
\begin{eq*}
\frac{1}{\|h\|_X}\| T(x+h) - T(x) - DT(x)h \|_Y \leq \frac{\varepsilon}{B+1}.
\end{eq*}
Let now $y \in L$ and $d \in Y$ with $\|d\|_Y \leq \frac{\delta}{B + 1}$ be chosen arbitrarily. Set ${x\mh{:}=S(y)} \in M$ and $h \mh{:}= S(y+d) - S(y)$. Then, using the relation $DS(y') = DT(S(y'))^{-1}$ for all $y' \in Y$ \mh{as well as $S(y+td) \in M$ for $t \in [0,1]$}{} yields the following estimate\mh{:}
\begin{eq*}
\|h\|_X &= \| S(y+d) - S(y) \|_X = \left\|\int_0^1 DS(y + td)d\, \mathrm{d}t \right\|_X\\
&\leq \int_0^1 \|DT(S(y+td))^{-1}\|_{\calL(Y,X)} \|d\|_Y \mathrm{d}t \leq B \|d\|_Y \leq \frac{B}{B+1}\delta  < \delta\mh{.}{}
\end{eq*}
Hence\msa{,}{} we deduce
\begin{eq*}
\| S(y+d) &- S(y) - DS(y)d \|_Y = \|DT^{-1}(x) \left( DT(x)h - d \right)\|_Y\\
&\leq \| DT(x)^{-1} \|_{\calL(Y,X)} \|DT(x)h - (y+d) + y\|_Y\\
&= \| DT(x)^{-1} \|_{\calL(Y,X)} \|DT(x)h - T(x + h) + T(x)\|_Y\\
&\leq\frac{B}{B+1}\varepsilon \|h\|_X \leq \left(\frac{B}{B+1}\right)^2 \varepsilon \|d\|_X \leq \varepsilon \|d\|_X.
\end{eq*}
This implies the uniform Fr\'echet differentiability of the map $S$ on every bounded subset.
\end{proof}

\begin{proof}[Proof of \refer{lem:uniformfrechet:composition}]
Let a bounded subset $M \subseteq X$ be given and take $\eps > 0$ arbitrarily. By the boundedness of $T_1$ the set $L := \conv{T_1(M + \B_X)}$ is bounded. Using the boundedness of $DT_1$ respectively $DT_2$ define
\begin{eq*}
B_1 &:= \sup_{x \in M + \B_X} \| DT_1(x) \|_{\calL(X,Y)} < +\infty,\\
\text{and}\rule{0.3\textwidth}{0pt}&\rule{0.65\textwidth}{0pt}\\
B_2 &:= \sup_{y \in L + \B_Y} \| DT_2(y) \|_{\calL(Y,Z)} < +\infty.
\end{eq*}
By the uniform Fr\'echet differentiability of $T_1$ there exist $\eps > 0$ and $\delta \in (0,1)$ with 
\begin{eq*}
\frac{1}{\|h\|_X}\| T_1(x+h) - T_1(x) - DT_1(x)h \|_Y \leq \frac{\eps}{B_1 + B_2 + 1}
\end{eq*}
for all $h \in X$ with $\|h\|_X \leq \frac{\delta}{B_1 + 1}$ and all $x \in M$ and
\begin{eq*}
\frac{1}{\|k\|_X}\| T_2(y+k) - T_2(y) - DT_2(y)k \|_Y \leq \frac{\eps}{B_1 + B_2 + 1}
\end{eq*}
for all $k \in Y$ with $\|k\|_Y \leq \delta$ and all $y \in L$. Then\msa{,}{} we obtain for all $y', y \in L + \B_Y$
\begin{eq}[\label{eq:tmp}]
\|T_2(y') - T_2(y)\|_Z &= \left\| \int_0^1 DT_2(y+t(y'-y))(y'-y) \mathrm{d}t \right\|_Z\\
&\leq \int_0^1 \|DT_2(y + t(y'-y))\|_{\calL(Y,Z)} \|y'-y\|_{\msa{Y}} \mathrm{d}t \leq B_2 \|y'-y\|_Y.
\end{eq}
For $x \in M$ and $h \in X$ with $\|h\|_X \leq \frac{\delta}{B_1+1}$ we \msa{get}{} $\|DT_1(x)h\|_Y\leq \delta$ and hence using the uniform Fr\'echet differentiablity \msa{the}{} estimate \refer{eq:tmp} yields
\begin{eq*}
\|T_2(T_1(x+h)) &- T_2(T_1(x)) - DT_2(T_1(x)) DT_1(x) h\|_Z\\
&\leq \| T_2(T_1(x+h)) - T_2(T_1(x) + DT_1(x)h)\|_Z\\
&+ \| T_2(T_1(x) + DT_1(x)h) - T_2(T_1(x)) - DT_2(T_1(x)) DT_1(x) h\|_Z\\
&\leq B_2 \|T_1(x+h) - T_1(x) - DT_1(x)h\|_{\msa{Y}} + \eps \|DT_1(x)h\|_Y\\
&\leq \frac{\eps}{B_1 + B_2 + 1} B_2 \|h\|_X + \frac{\eps}{B_1 + B_2 + 1} B_1 \|h\|_X \leq \eps \|h\|_X\mh{,}{}
\end{eq*}
\mh{which ends the proof.}{}
\end{proof}

\bibliographystyle{alpha}
\bibliography{../../../mystyle/mylibrary}
\end{document}

%% file: domain.tex
\begin{tikzpicture}[scale = 4]
\draw[color = gray, ->] (0,0) -- (1.10,0);
\draw[color = gray, ->] (0,0) -- (0,1.10);

\draw (0,0) -- (1,0) -- (1,1) -- (0,1) -- cycle;

\draw[color = gray] (-0.025,1) -- (0.025,1);
\node[anchor = east, inner sep = 4pt] at (0,1) {$1$};

\draw[color = gray] (-0.025,0.25) -- (0.025,0.25);
\node[anchor = east, inner sep = 4pt] at (0,0.25) {$0.25$};

\draw[color = gray] (0.75,-0.025) -- (0.75,0.025);
\node[anchor = north, inner sep = 4pt] at (0.75,0) {$0.75$};

\draw[color = gray] (1,-0.025) -- (1,0.025);
\node[anchor = north, inner sep = 4pt] at (1,0) {$1$};

\draw[color = gray] (0.25,-0.025) -- (0.25,0.025);
\draw[color = gray, dotted] (0.25,0) -- (0.25,0.25);
\node[anchor = north, inner sep = 4pt] at (0.25,0) {$0.25$};

\draw[color = gray] (-0.025,0.5) -- (0.025,0.5);
\draw[color = gray, dotted] (0,0.5) -- (0.13,0.5);
\draw[color = gray, dotted] (0.26,0.5) -- (0.5,0.5);
\node[anchor = east, inner sep = 4pt] at (0,0.5) {$0.5$};

\node[anchor = north east, inner sep = 2pt] at (0,0) {$0$};
\node[anchor = south west, inner sep = 2pt] at (1,1) {$(1,1)$};

\draw (0,0) -- (0,0.25) -- (0.25,0.25) -- (0.5,0.5) -- (0.75,0.25) -- (0.75,0) -- cycle;
\node[anchor = center] at (0.5,0.2) {$\omega_1$};

\draw (1,0) -- (0.75,0) -- (0.75,0.25) -- (0.5,0.5) -- (0.75,0.75) -- (1,0.75) -- cycle;
\node[anchor = center] at (0.8,0.5) {$\omega_2$};

\draw (1,1) -- (1,0.75) -- (0.75,0.75) -- (0.5,0.5) -- (0.25,0.75) -- (0.25,1) -- cycle;
\node[anchor = center] at (0.5,0.8) {$\omega_3$};

\draw (0,1) -- (0.25,1) -- (0.25,0.75) -- (0.5,0.5) -- (0.25,0.25) -- (0,0.25) -- cycle;
\node[anchor = center] at (0.2,0.5) {$\omega_4$};

\end{tikzpicture}

%% file: nash_elliptic_semi_distr/gamma.tex
\begin{tikzpicture}
\begin{axis}[
  xlabel=$n$,
  ylabel=$\gamma_n$,
  ymode=log,
  ymax = 1e9]
\addplot table {nash_elliptic_semi_distr/gamma.dat};
\end{axis}
\end{tikzpicture}

%% file: nash_elliptic_semi_distr_coop/gamma.tex
\begin{tikzpicture}
\begin{axis}[
  xlabel=$n$,
  ylabel=$\gamma_n$,
  ymode = log]
\addplot table {nash_elliptic_semi_distr_coop/gamma.dat};
\end{axis}
\end{tikzpicture}

%% file: nash_elliptic_semi_distr/objectives.tex
\begin{tikzpicture}
\begin{axis}[
  xlabel=$n$,
  ylabel=$\sum_{i = 1}^N \calJ_i\left(u^{\gamma^n}\right)$]
\addplot table[x index=0,y index=5] {nash_elliptic_semi_distr/objectives.dat};
\end{axis}
\end{tikzpicture}

%% file: nash_elliptic_semi_distr_coop/objectives.tex
\begin{tikzpicture}
\begin{axis}[
  xlabel=$n$,
  ylabel=$\sum_{i = 1}^N \calJ_i\left(u^{\gamma^n}\right)$]
\addplot table[x index=0,y index=5] {nash_elliptic_semi_distr_coop/objectives.dat};
\end{axis}
\end{tikzpicture}

%% file: nash_elliptic_semi_distr/newton_iter.tex
\begin{tikzpicture}
\begin{axis}[
  xlabel=$n$,
  ylabel=Number of Newton iterations,
  ymax = 8]
\addplot table {nash_elliptic_semi_distr/newton_iter.dat};
\end{axis}
\end{tikzpicture}

%% file: nash_elliptic_semi_distr_coop/newton_iter.tex
\begin{tikzpicture}
\begin{axis}[
  xlabel=$n$,
  ylabel=Number of Newton iterations,
  ymax = 100]
\addplot table {nash_elliptic_semi_distr_coop/newton_iter.dat};
\end{axis}
\end{tikzpicture}

%% file: nash_elliptic_semi_bdry/gamma.tex
\begin{tikzpicture}
\begin{axis}[
  xlabel=$n$,
  ylabel=$\gamma_n$,ymode=log,ymax = 1e11]
\addplot table {nash_elliptic_semi_bdry/gamma.dat};
\end{axis}
\end{tikzpicture}

%% file: nash_elliptic_semi_bdry_coop/gamma.tex
\begin{tikzpicture}
\begin{axis}[
  xlabel=$n$,
  ylabel=$\gamma_n$,ymode=log,ymax = 1e11]
\addplot table {nash_elliptic_semi_bdry_coop/gamma.dat};
\end{axis}
\end{tikzpicture}

%% file: nash_elliptic_semi_bdry/objectives.tex
\begin{tikzpicture}
\begin{axis}[
  xlabel=$n$,
  ylabel=$\sum_{i = 1}^N \calJ_i\left(u^{\gamma^n}\right)$]
\addplot table[x index=0,y index=5] {nash_elliptic_semi_bdry/objectives.dat};
\end{axis}
\end{tikzpicture}

%% file: nash_elliptic_semi_bdry_coop/objectives.tex
\begin{tikzpicture}
\begin{axis}[
  xlabel=$n$,
  ylabel=$\sum_{i = 1}^N \calJ_i\left(u^{\gamma^n}\right)$]
\addplot table[x index=0,y index=5] {nash_elliptic_semi_bdry_coop/objectives.dat};
\end{axis}
\end{tikzpicture}

%% file: nash_elliptic_semi_bdry/newton_iter.tex
\begin{tikzpicture}
\begin{axis}[
  xlabel=$n$,
  ylabel=Number of Newton iterations,
  ymax = 40]
\addplot table {nash_elliptic_semi_bdry/newton_iter.dat};
\end{axis}
\end{tikzpicture}

%% file: nash_elliptic_semi_bdry_coop/newton_iter.tex
\begin{tikzpicture}
\begin{axis}[
  xlabel=$n$,
  ylabel=Number of Newton iterations,
  ymax = 40]
\addplot table {nash_elliptic_semi_bdry_coop/newton_iter.dat};
\end{axis}
\end{tikzpicture}

%% file: preprint_shortened.bbl
\newcommand{\etalchar}[1]{$^{#1}$}
\begin{thebibliography}{KKSW19}

\bibitem[ABH{\etalchar{+}}15]{bib:fenics}
M.~S. Aln{\ae}s, J.~Blechta, J.~Hake, A.~Johansson, B.~Kehlet, A.~Logg,
  C.~Richardson, J.~Ring, M.~E. Rognes, and G.~N. Wells.
\newblock The {FEniCS} {P}roject {V}ersion 1.5.
\newblock {\em Archive of Numerical Software}, 3(100), 2015.

\bibitem[ACM11]{bib:AusselCoreaMarechal2011}
D.~Aussel, R.~Correa, and M.~Marechal.
\newblock Gap {F}unctions for {Q}uasivariational {I}nequalities and
  {G}eneralized {N}ash {E}quilibrium {P}roblems.
\newblock {\em Journal of {O}ptimization {T}heory and {A}pplications},
  151(3):474, 2011.

\bibitem[AHS18]{bib:AdamHintSur}
L.~Adam, M.~Hinterm\"uller, and T.~M. Surowiec.
\newblock A semismooth {N}ewton method with analytical path-following for the
  {$H^1$}-projection onto the {G}ibbs simplex.
\newblock {\em IMA Journal of Numerical Analysis}, 39(3):1276--1295, 06 2018.

\bibitem[BK13]{bib:BorziKanzow}
A.~Borzi and C.~Kanzow.
\newblock {Formulation and numerical solution of Nash equilibrium
  multiobjective elliptic control problems}.
\newblock {\em SIAM Journal on Control and Optimization}, 51(1):718--744, 2013.

\bibitem[BP12]{bib:BarbuPrecupanu}
V.~Barbu and T.~Precupanu.
\newblock {\em Convexity and {O}ptimization in {B}anach Spaces}.
\newblock Springer Monographs in Mathematics. Springer Netherlands, 2012.

\bibitem[Bra68]{bib:BraessParadox}
D.~Braess.
\newblock {\"U}ber ein {P}aradoxon aus der {V}erkehrsplanung.
\newblock {\em Unternehmensforschung}, 12(1):258--268, 1968.

\bibitem[BS00]{bib:BonnansShapiro}
J.F. Bonnans and A.~Shapiro.
\newblock {\em {Perturbation Analysis of Optimization Problems}}.
\newblock Springer Series in Operations Research and Financial Engineering.
  Springer New York, 2000.

\bibitem[BV04]{bib:BoydVandenberghe}
S.~Boyd and L.~Vandenberghe.
\newblock {\em {Convex Optimization}}.
\newblock Cambridge University Press, 2004.

\bibitem[CNQ00]{bib:ChenNashedZuhair}
X.~Chen, Z.~Nashed, and L.~Qi.
\newblock Smoothing methods and semismooth methods for nondifferentiable
  operator equations.
\newblock {\em SIAM Journal on Numerical Analysis}, 38(4):1200--1216, 2000.

\bibitem[Deu05]{bib:DeuflhardNewton}
P.~Deuflhard.
\newblock {\em {Newton Methods for Nonlinear Problems: Affine Invariance and
  Adaptive Algorithms}}.
\newblock Springer Series in Computational Mathematics. Springer Berlin
  Heidelberg, 2005.

\bibitem[DF95]{bib:DirkseFerris}
S.~P. Dirkse and M.~C. Ferris.
\newblock The path solver: a nonmonotone stabilization scheme for mixed
  complementarity problems.
\newblock {\em Optimization methods and software}, 5(2):123--156, 1995.

\bibitem[Dut13]{bib:Dutang}
C.~Dutang.
\newblock {Existence theorems for generalized Nash equilibrium problems: An
  analysis of assumptions}.
\newblock {\em Journal of Nonlinear Analysis and Optimization: Theory \&
  Applications}, 4(2):115--126, 2013.

\bibitem[DvKF13]{bib:DrevesHeusingerKanzowFukushima}
A.~Dreves, A:~{von Heusinger}, C~Kanzow, and M.~Fukushima.
\newblock A globalized {N}ewton method for the computation of normalized {N}ash
  equilibria.
\newblock {\em Journal of Global Optimization}, 56(2):327--340, 2013.

\bibitem[Ehr62]{bib:Ehrmann}
H.~H. Ehrmann.
\newblock On implicit function theorems and the existence of solutions of
  non-linear equations.
\newblock Technical report, Wisconsin Univ Madison Mathematics Research Center,
  1962.

\bibitem[FFP09]{bib:FacchineiFischerPiccialli}
F.~Facchinei, A.~Fischer, and V.~Piccialli.
\newblock Generalized {N}ash equilibrium problems and {N}ewton methods.
\newblock {\em Mathematical Programming}, 117(1-2):163--194, 2009.

\bibitem[FK07]{bib:KanzowFacchinei}
F.~Facchinei and C.~Kanzow.
\newblock Generalized {N}ash equilibrium problems.
\newblock {\em 4or}, 5(3):173--210, 2007.

\bibitem[Gri85]{bib:Grisvard}
P.~Grisvard.
\newblock {\em {Elliptic Problems in Nonsmooth Domains}}.
\newblock Classics in Applied Mathematics. Society for Industrial and Applied
  Mathematics, 1985.

\bibitem[HIK02]{bib:HintItoKunisch}
M.~Hinterm{\"u}ller, K.~Ito, and K.~Kunisch.
\newblock {The Primal-Dual Active Set Strategy as a Semismooth {N}ewton
  Method}.
\newblock {\em SIAM Journal on Optimization}, 13(3):865--888, 2002.

\bibitem[HK06a]{bib:HintKunischPathLowMultiplier}
M.~Hinterm{\"u}ller and K.~Kunisch.
\newblock {Feasible and Noninterior Path-Following in Constrained Minimization
  with Low Multiplier Regularity}.
\newblock {\em SIAM Journal on Control and Optimization}, 45(4):1198--1221,
  2006.

\bibitem[HK06b]{bib:HintKunischPath}
M.~Hinterm{\"u}ller and K.~Kunisch.
\newblock {Path-following Methods for a Class of Constrained Minimization
  Problems in Function Space}.
\newblock {\em SIAM J. Optim.}, 17(1):159--187, 2006.

\bibitem[HPUU08]{bib:HinzePinnauMUlbrichSUlbrich}
M.~Hinze, R.~Pinnau, M.~Ulbrich, and S.~Ulbrich.
\newblock {\em {Optimization with PDE Constraints}}.
\newblock Mathematical Modelling: Theory and Applications. Springer
  Netherlands, 2008.

\bibitem[HR15]{bib:HintRasch}
M.~Hinterm{\"u}ller and J.~Rasch.
\newblock Several path-following methods for a class of gradient constrained
  variational inequalities.
\newblock {\em Computers \& Mathematics with Applications}, 69(10):1045--1067,
  2015.

\bibitem[HS13]{bib:HintSurGNEP}
M.~Hinterm\"uller and T.~M. Surowiec.
\newblock {A PDE-constrained generalized {N}ash equilibrium problem with
  pointwise control and state constraints}.
\newblock {\em Pac. J. Optim}, 9(2):251--273, 2013.

\bibitem[HS20]{bib:HintStenglKConvex}
M.~Hinterm\"uller and S.M. Stengl.
\newblock {On the convexity of optimal control problems involving non-linear
  PDEs or VIs and applications to {N}ash games}.
\newblock September 2020.

\bibitem[HS21]{bib:HintStenglEquiGamma}
M.~Hinterm{\"u}ller and S.M. Stengl.
\newblock {A Generalized $\Gamma$-Convergence Concept for a Class of
  Equilibrium Problems}.
\newblock {\em WIAS Preprint No. 2879}, 2021.

\bibitem[HSK15]{bib:HintSurKaemmler}
M.~Hinterm{\"u}ller, T.~M. Surowiec, and A.~K{\"a}mmler.
\newblock {Generalized {N}ash Equilibrium Problems in {B}anach Spaces: Theory,
  {N}ikaido–-{I}soda-Based Path-Following Methods, and Applications}.
\newblock {\em SIAM Journal on Optimization}, pages 1826--1856, 2015.

\bibitem[KKSW19]{bib:KanzowKarlSteckDWachsmuth}
C.~Kanzow, V.~Karl, D.~Steck, and D.~Wachsmuth.
\newblock {The Multiplier-Penalty Method for Generalized {N}ash Equilibrium
  Problems in {B}anach Spaces}.
\newblock {\em SIAM Journal on Optimization}, 29(1):767--793, 2019.

\bibitem[Lla86]{bib:Llavona}
J.G. Llavona.
\newblock {\em {Approximation of Continuously Differentiable Functions}}.
\newblock North-Holland Mathematics Studies. Elsevier Science, 1986.

\bibitem[Nas50]{bib:NashNPerson}
J.F. Nash.
\newblock {Equilibrium points in n-person games}.
\newblock {\em Proceedings of the National Academy of Sciences}, 36(1):48--49,
  1950.

\bibitem[NI55]{bib:NikaidoIsoda}
H.~Nikaid{\^o} and K.~Isoda.
\newblock Note on non-cooperative convex games.
\newblock {\em Pacific Journal of Mathematics}, 5(Suppl. 1):807--815, 1955.

\bibitem[NRTV07]{bib:NisanRoughgardenAlgorithmicGameTheory}
N.~Nisan, T.~Roughgarden, E.~Tardos, and V.V. Vazirani.
\newblock {\em {Algorithmic Game Theory}}.
\newblock Cambridge University Press, 2007.

\bibitem[PF05]{bib:PangFukushima}
J.-S. Pang and M.~Fukushima.
\newblock Quasi-variational inequalities, generalized {N}ash equilibria, and
  multi-leader-follower games.
\newblock {\em Computational Management Science}, 2(1):21--56, 2005.

\bibitem[Ros65]{bib:Rosen}
J.~B. Rosen.
\newblock {Existence and Uniqueness of Equilibrium Points for Concave N-Person
  Games}.
\newblock {\em Econometrica}, 33(3):520--534, 1965.

\bibitem[Sav98]{bib:Savare}
G.~Savar{\'e}.
\newblock {Regularity Results for Elliptic Equations in Lipschitz Domains}.
\newblock {\em Journal of Functional Analysis}, 152(1):176--201, 1998.

\bibitem[vH09]{bib:vonHeusinger}
A.~von Heusinger.
\newblock {\em {Numerical Methods for the Solution of the Generalized {N}ash
  Equilibrium Problem}}.
\newblock Doctoral thesis, Universit{\"a}t W{\"u}rzburg, 2009.

\bibitem[vK09]{bib:KanzowHeusingerOptReformulation}
A.~{von Heusinger} and C.~Kanzow.
\newblock Optimization reformulations of the generalized {N}ash equilibrium
  problem using {N}ikaido-{I}soda-type functions.
\newblock {\em Computational Optimization and Applications}, 43(3):353--377,
  2009.

\bibitem[ZK79]{bib:ZoweKurcyusz}
J.~Zowe and S.~Kurcyusz.
\newblock Regularity and stability for the mathematical programming problem in
  {B}anach spaces.
\newblock {\em Applied Mathematics and Optimization}, 5(1):49--62, Mar 1979.

\end{thebibliography}
